\documentclass[11pt]{article}

\usepackage{jabbrv}
\usepackage{float}
\usepackage{fancyhdr}
\usepackage{amsmath,amsfonts}
\usepackage{graphicx}
\usepackage{color}
\usepackage[footnotesize]{caption}
\usepackage{epstopdf}
\usepackage[Lenny]{fncychap}
\usepackage{algorithm}
\usepackage{algpseudocode}
\usepackage{tocvsec2}
\usepackage{mathtools}
\usepackage{nicefrac}
\mathtoolsset{showonlyrefs}

\makeatletter

\newcommand{\Rmnum}[1]{\expandafter\@slowromancap\romannumeral #1@}
\makeatother

\newtheorem{oss}{Remark}

\def\R{{\mathbb R}}

\def\N{{\mathbb N}}

\let\theta\vartheta

\def\vec#1{{\mathbf{#1}}}

\DeclareMathOperator{\sgn}{sign}

\newcommand{\pad}[2]{\frac{\partial #1}{\partial #2}}

\usepackage{authblk}
\title{A multigrid ghost-point level-set method for incompressible Navier-Stokes equations on moving domains with curved boundaries}
\author[$\star$]{Armando Coco}
\affil[$\star$]{School of Engineering, Computing and Mathematics, Oxford Brookes University, Wheatley campus, Oxford (United Kingdom)}

\begin{document}
\sffamily
\maketitle

\begin{abstract}
In this paper we present a numerical approach to solve the Navier-Stokes equations on moving domains with second-order accuracy.
The space discretization is based on the ghost-point method, which falls under the category of unfitted boundary methods, since the mesh does not adapt to the moving boundary. 
The equations are advanced in time by using Crank-Nicholson. The momentum and continuity equations are solved simultaneously for the velocity and the pressure by adopting a proper multigrid approach. To avoid the checkerboard instability for the pressure, a staggered grid is adopted, where velocities are defined at the sides of the cell and the pressure is defined at the centre. The lack of uniqueness for the pressure is circumvented by the inclusion of an additional scalar unknown, representing the average divergence of the velocity, and an additional equation to set the average pressure to zero.
Several tests are performed to simulate the motion of an incompressible fluid around a moving object, as well as the lid-driven cavity tests around steady objects. The object is implicitly defined by a level-set approach, that allows a natural computation of geometrical properties such as distance from the boundary, normal directions and curvature.
Different shapes are tested: circle, ellipse and flower. Numerical results show the second order accuracy for the velocity and the divergence (that decays to zero with second order) and the efficiency of the multigrid, that is comparable with the tests available in literature for rectangular domains without objects, showing that the presence of a complex-shaped object does not degrade the performance.
\end{abstract}

\paragraph*{Keywords:}
finite-difference; unfitted boundary methods; high-order accuracy; staggered grid; incompressible fluid dynamics; complex-shaped domains

\section{Introduction}
Navier-Stokes equations are a system of partial differential equations (PDEs) describing the motion of an incompressible viscous fluid subject to external forces. In the non-dimensional form, they read
\begin{equation}\label{NSeq}
\begin{aligned}
\frac{\partial \vec{u}}{\partial t} + 
\vec{u} \cdot \nabla \vec{u}
+ \nabla p &= \frac{1}{\text{Re}} \nabla^2 \vec{u} + \vec{f} \\
\nabla \cdot \vec{u} &= 0
\end{aligned}
\end{equation}
where $\vec{u}$ is the fluid velocity, $p$ is the pressure, $\text{Re}$ is the Reynolds number and $\vec{f}$ the external force.
The first equation is obtained from the momentum equation, while the second one is the continuity equation and represents the incompressibility condition.
Navier-Stokes equations have been central to innumerable scientific and engineering applications for several decades, as the necessity to model the motion of incompressible fluids has a tremendous impact on the real-life world. 
Several numerical approaches have been developed so far to provide efficient and accurate methods in different contexts, both for time and space discretization. An overview on traditional methods can be found in~\cite{langtangen2002numerical}.

One of the first approaches for time discretization was the {\it projection method}, proposed by Chorin in~\cite{chorin1997numerical}. It consists of a fractional step, where the velocity $\vec{u}$ is advanced in time by discretising the momentum equation and ignoring the pressure term, obtaining an intermediate velocity $\vec{u}^*$. This velocity is not divergence-free, and therefore is corrected by the pressure term, that is determined by enforcing the divergence-free condition.
Projection methods have been very popular for Navier-Stokes equations and are widely adopted nowadays, as they are computationally efficient and easy to implement (they require the solution of a sequence of elliptic equations for the intermediate velocity and for the pressure). However, the definition of boundary conditions for $\vec{u}^*$ and $p$ is not obvious and it is challenging to achieve an accuracy higher than first order, even if the internal equations are discretized with higher order. There are existing approaches that suggest appropriate boundary conditions to improve the accuracy, as described in~\cite{brown2001accurate}. However, these boundary conditions are fully coupled each other and the resulting elliptic equations for $\vec{u}^*$ and for $p$ cannot be solved in sequence, as the boundary condition on $p$ depends on $\vec{u}^*$.

Another approach consists of discretising both PDEs of Eq.~\eqref{NSeq} and solve them simultaneously in the unknowns $(\vec{u},p)$. This is the {\it monolithic} approach and has the advantage that no boundary conditions are required for the pressure, since the order of the system of PDEs \eqref{NSeq} is the same as the number of components of $\vec{u}$~\cite{langtangen2002numerical}. 
These methods can usually achieve higher order accuracy than projection methods, but they have the disadvantage that extra computational effort is needed to solve both PDEs simultaneously. This aspect can be alleviated if efficient solvers are adopted such as multigrid specifically designed for this class of problems.

Projection and monolithic methods refer only to time discretization and then spatial discretization methods must be addressed regardless of the choice between the methods described above. Standard problems on rectangular domains are usually discretized in space by using a finite volume approach on a staggered grid (MAC grid, see for example~\cite{srinivas2002computational}), where the velocities are defined at the edges of the cell and the pressure is defined at the centre.
The MAC grid is used for example to avoid the checkerboard instability for the pressure term observed in non-staggered grids due to the fact that $p$ appears in the equations only in the form of $\nabla p$.
The extension of MAC grids to non-rectangular domains (curved boundaries) is not straightforward, since a particular treatment is needed close to the boundary. Moreover, there is an increasing scientific interest to Navier-Stokes equations for complex-shaped (possibly moving) domains, for example to model the motion of an incompressible fluids around a deforming object, that may be central to more sophisticated applications such as fluid-structure interactions.
Numerical methods for steady complex-shaped objects are mainly divided between fitted and unfitted boundary methods. In fitted boundary methods, the grid, either structured or unstructured (such as Finite Element Methods (FEM)), is adapted to the curved boundary~\cite{shukla2007very, gross2007extended, girault2012finite, masud1997space, lozovskiy2018quasi, danilov2017finite, court2015fictitious}.
Although they provide an extremely accurate representation of the domain, the generation of the mesh may become cumbersome and significantly expensive from a computational point of view, especially in moving domain problems where a new mesh generation is needed at each time step.
In unfitted boundary methods the moving domain is instead embedded in a fixed grid.
The methods of this class are based either on FEM, such as XFEM~\cite{heimann2013unfitted, peterseim2011finite} or GFEM~\cite{christiansen2018generalized}, on meshless methods~\cite{yang2019solving}, or Cartesian meshes, such as cut cell methods~\cite{gokhale2018dimensionally}.

An alternative approach of finite--difference methods for curved boundaries is represented by the Ghost--Point method.
A pioneering work of this class of methods is the Immersed Boundary Method~\cite{Peskin:IBM} to model blood flows in the heart, where Peskin presented a first-order accurate method, later extended to higher order by the Immersed Interface Methods proposed by LeVeque and Li in~\cite{LeVequeLi:IIM}.
More recent finite-difference methods to discretize PDEs on complex domains are the Ghost-Fluid Methods proposed
in~\cite{Fedkiw:GFM, Gibou:Ghost, Gibou:fourth_order, Gibou:fluid_solid},
where the grid points external to the domain are called ghost points and a fictitious (ghost) value of the numerical solution is formally extrapolated onto these grid points in order to maintain a stable algorithm and achieve the desired accuracy order.
The main feature of these methods is that the internal equations are firstly formally solved for the ghost values, resulting in a final linear system with eliminated boundary conditions. This class of methods may fail to obtain high accuracy for Neumann boundary condition, and anyway it remains first order accurate in the post-processing reconstruction of the gradient of the solution.

A new finite-difference method with an improved accuracy was proposed by the author and Russo in~\cite{CocoRusso:Elliptic, coco2018second, coco2012second} for elliptic equations and extended to thermo-poroelastic Cauchy-Navier equations for elastic deformation~\cite{CCDR:LinearElasticity}, Euler equations for gas dynamics~\cite{chertock2018second, CocoRusso:Hyp2012}, electron transport in semiconductors~\cite{muscato2019hydrodynamic, muscato2018low}, fluid flow in porous media for water-rock interaction \cite{coco2016numerical, coco2016hydro}.
In this approach, the ghost point values are eventually coupled each other, resulting in a bigger linear system with non-eliminated boundary conditions,
solved by a multigrid approach suitably designed for ghost-point methods.
The domain is implicitly described by a level-set function.

Other numerical methods to solve sharp-edge boundary problems are the finite volume methods, the non-symmetric positive definite finite element method, the matched interface and boundary (MIB) method~\cite{Zhou:MIB}, the arbitrary Lagrangian Eulerian method (ALE)~\cite{FormaggiaNobile:ALE, Donea:ALE},
the penalization methods~\cite{lacanette2009eulerian, Angot:penalization, Iollo:penalization}, and the class of Immersed Finite Volume Methods (IFVM)~\cite{ewing1999immersed}.
Other recent advances have been obtained in~\cite{helgadottir2015imposing} to achieve higher accuracy in the presence of Neumann boundary conditions, and in~\cite{bochkov2019solving} to obtain higher accuracy in the gradient of the solution.

Time and space discretizations of Navier-Stokes equations lead to a linear system that must be solved efficiently.
Multigrid method is one of the most efficient iterative approaches to solve the linear systems arising from the discretization of a class of PDEs. Initially designed for elliptic problems, multigrid methods have been widely adopted to solve countless problems~\cite{Trottemberg:MG, hackbusch2013multi, briggs2000multigrid}.
Several approaches have been proposed for boundaries aligned with grid lines. As example, we mention~\cite{alcouffe1981multi, dendy1982black, adams2002immersed, li2001maximum, adams2004new}. These approaches fall under the class of geometric multigrid, since the operators are defined from the geometrical background of the problem. An extension of geometric multigrid is represented by the class of algebraic multigrid methods~\cite{stuben2001review}, where the operators are defined merely from the information gathered from the linear systems, and therefore it is more practical for complex problems such as curved boundaries. However, the inclusion of geometrical information to the definition of the operators usually results in a more efficient solver, as highlighted for example in~\cite{adams2005comparison}, where a comparison between geometric and algebraic multigrid methods for curved boundaries is performed.

In this paper, we provide a finite-difference method to solve the Navier-Stokes equations on moving domains following a monolithic approach for the time discretization and a finite-difference ghost-point method for the space discretization. The singularity of the problem (due to the fact that $p$ is not uniquely defined, since it appears only in terms of $\nabla p$) is circumvented by an additional equation for the pressure that requires that its average is zero. This additional equation is balanced with an additional unknown $\xi$ that measures the average divergence of the velocity (that is therefore not necessarily zero, as in traditional methods, but decays with the same order of accuracy of the method).
The linear system arising from the discretization is solved by a geometric multigrid approach, where the relaxation, restriction and interpolation operators are modified close to the boundary in order to maintain the same efficiency that standard multigrid techniques achieve away from the boundary.

The paper is structured as follows. In Sect.~\ref{sect:disctime} and \ref{sect:space} we describe the time and space discretizations, respectively. The multigrid technique is detailed in Sect.~\ref{sect:mg}, while numerical tests are provided in Sect.~\ref{sect:numtests} and conclusions in Sect.\ref{sect:concl}.

\section{Discretization in time}\label{sect:disctime}
The Navier-Stokes equations \eqref{NSeq} are discretized in time using Crank-Nicholson and solved for $\vec{u}^{n+1}$ in the domain $\Omega^{(n+1)}$
\begin{equation}\label{NSCN}
\begin{aligned}
\frac{\vec{u}^{(n+1)}-\vec{u}^{(n)}}{\Delta t} + \left( \vec{u} \cdot \nabla \vec{u} \right)^{(n+1/2)}
+ \nabla p^{(n+1/2)} &= \frac{1}{2\,\text{Re}} \left( \nabla^2 \vec{u}^{(n)} + \nabla^2 \vec{u}^{(n+1)} \right) + \vec{f}^{(n+1/2)} \text{ in } \Omega^{(n+1)} \\
\nabla \cdot \vec{u}^{(n+1)} &= 0 \text{ in } \Omega^{(n+1)} \\
\left. \vec{u}^{(n+1)} \right|_{\partial \Omega^{(n+1)}} &= \vec{u}_b
\end{aligned}
\end{equation}
Taking the implicit terms on the left hand side and the explicit terms on the right hand side and multiplying by $\Delta t$, Eq.~\eqref{NSCN} reads:
\begin{equation}\label{NSCN1}
\begin{aligned}
\vec{u}^{(n+1)} - \frac{\Delta t}{2\,\text{Re}} \nabla^2 \vec{u}^{(n+1)}
+ \Delta t \, \nabla p^{(n+1/2)} &= \tilde{\vec{f}} - \Delta t \left( \vec{u} \cdot \nabla \vec{u} \right)^{(n+1/2)} \\
\nabla \cdot \vec{u}^{(n+1)} &= 0 \\
\left. \vec{u}^{(n+1)} \right|_{\partial \Omega^{(n+1)}} &= \vec{u}_b
\end{aligned}
\end{equation}
where $\tilde{\vec{f}} = \vec{u}^{(n)} + \displaystyle \frac{\Delta t}{2\,\text{Re}} \nabla^2 \vec{u}^{(n)} + \Delta t \, \vec{f}^{(n+1/2)}$.
The first question to address is the treatment of the convective term $\left( \vec{u} \cdot \nabla \vec{u} \right)^{(n+1/2)}$. A natural approach is to approximate this term fully explicitly $\left( \vec{u} \cdot \nabla \vec{u} \right)^{(n+1/2)} \approx \vec{u}^{(n)} \cdot \nabla \vec{u}^{(n)}$ or with an implicit-explicit technique
$\left( \vec{u} \cdot \nabla \vec{u} \right)^{(n+1/2)} \approx \vec{u}^{(n)} \cdot \nabla \vec{u}^{(n+1)}$, which however are only first order accurate.
A higher order fully explicit approximation can be obtained by extrapolating the convective term from previous time steps at the desired accuracy. For example, for a second order approximation, the extrapolation is
\begin{equation}\label{extrUdU}
\left( \vec{u} \cdot \nabla \vec{u} \right)^{(n+1/2)} \approx \left( \vec{u} \cdot \nabla \vec{u} \right)^{(n)}
+ \frac{1}{2} \left( \vec{u} \cdot \nabla \vec{u} \right)^{(n-1)}
- \frac{1}{2} \left( \vec{u} \cdot \nabla \vec{u} \right)^{(n-2)}.
\end{equation}
This approximation leads to a multistep method using the solution at the three previous time steps. For this reason,
the solution at the first two time steps $\vec{u}^{(i)}$ with $i=1,2$ must be computed with an alternative method (observe that $\vec{u}^{(0)}$ is set up from the initial condition).
Another second-order approximation is provided by the mid-point discretization
\[
\left( \vec{u} \cdot \nabla \vec{u} \right)^{(n+1/2)} \approx \frac{1}{2} \left( \left(\vec{u} \cdot \nabla \vec{u}\right)^{(n)} + \left( \vec{u} \cdot \nabla \vec{u}\right)^{(n+1)} \right).
\]
With this approximation, Eqs.~\eqref{NSCN1} becomes a nonlinear equation due to the presence of the nonlinear convective term $\vec{u}^{(n+1)} \cdot \nabla \vec{u}^{(n+1)}$. The nonlinearity can be circumvented by using, for example, a Newton-Raphson approach, or alternatively an iterative method (with an iterative parameter $k$)
\begin{equation}\label{NSitconv}
\begin{aligned}
\vec{u}^{(n+1,k+1)} - \frac{\Delta t}{2\,\text{Re}} \nabla^2 \vec{u}^{(n+1,k+1)}
+ \Delta t \, \nabla p^{(n+1/2,k+1)} &= \tilde{\vec{f}} - \Delta t \left( \vec{u} \cdot \nabla \vec{u} \right)^{(n+1/2,k)} \\
\nabla \cdot \vec{u}^{(n+1,k+1)} &= 0 \\
\left. \vec{u}^{(n+1,k+1)} \right|_{\partial \Omega^{(n+1)}} &= \vec{u}_b
\end{aligned}
\end{equation}
with 
\[
\left( \vec{u} \cdot \nabla \vec{u} \right)^{(n+1/2,k)} = \frac{ \left(\vec{u} \cdot \nabla \vec{u}\right)^{(n)} + \left( \vec{u} \cdot \nabla \vec{u} \right)^{(n+1,k)} }{2}.
\]
The iteration algorithm is performed until $\vec{u}^{(n+1,k+1)}$ satisfies Eqs.~\eqref{NSCN1} within a certain tolerance $\varepsilon$ (say $\varepsilon = 10^{-6}$), i.e.
\begin{equation}\label{iterativeNS}
\left\| \tilde{\vec{f}} - \Delta t \left( \vec{u} \cdot \nabla \vec{u} \right)^{(n+1/2,k+1)} - \vec{u}^{(n+1,k+1)} + \frac{\Delta t}{2\,\text{Re}} \nabla^2 \vec{u}^{(n+1,k+1)}
- \Delta t \, \nabla p^{(n+1/2,k+1)} \right\|
< \max \left\{ \varepsilon, \varepsilon \left\| \tilde{\vec{f}} \right\|  \right\}.
\end{equation}
Then, $\vec{u}^{(n+1)} \coloneqq \vec{u}^{(n+1,k+1)}$ and the scheme is moved to the next time step.

Other approaches can be adopted, for example approximation techniques based on the Godunov method~\cite{bell1989second}, or Runge-Kutta time discretization~\cite{nikitin2006third}. A detailed treatment is beyond the scope of this paper, as we mainly focus on the treatment of curved boundaries and the method proposed for that purpose do not rely on the approach adopted to approximate the nonlinear term.

In this paper, we adopt a hybrid approach: we use the iteration technique \eqref{NSitconv} for the first two time steps and the extrapolation \eqref{extrUdU} from the third time step.

Either if we are using the iteration technique \eqref{NSitconv} or the extrapolation approach \eqref{extrUdU}, the problem that we have to solve at each iteration can be formulated as follows: find $(\vec{u},p)$ such that
\begin{equation}\label{spaceprob}
\begin{aligned}
\vec{u} - \alpha \nabla^2 \vec{u}
+ \nabla p &= \tilde{\vec{f}} \quad \text{ in } \Omega \\
\nabla \cdot \vec{u} &= 0 \quad \text{ in } \Omega \\
\vec{u} &= \vec{u}_b \quad \text{ on } \partial \Omega \\
\end{aligned}
\end{equation}
where $\Omega = \Omega^{(n+1)}$,  $\alpha = \Delta t/(2\,\text{Re})$, and $p$ and $\tilde{\vec{f}}$ are redefined as $p \coloneqq \Delta t \, p$ and 
$\tilde{\vec{f}} \coloneqq  \tilde{\vec{f}} - \Delta t \left( \vec{u} \cdot \nabla \vec{u} \right)^{(n+1/2)}$.

The approach suggested by the projection or gauge methods consists of solving an intermediate step 
\begin{equation}\label{proj}
\vec{u}^* - \alpha \nabla^2 \vec{u}^* = \tilde{\vec{f}}
\end{equation}
and then projecting the intermediate solution $\vec{u}^*$ onto the space of divergence-free functions. In detail, $\vec{u}^*$ is split as $\vec{u}^* = \vec{u} + \nabla \chi $, where $\vec{u}$ is divergence-free. The auxiliary function $\chi$ is determined by taking the divergence of both sides, leading to an elliptic problem $\nabla^2 \chi = \nabla \cdot \vec{u}^*$. Finally, the solution $\vec{u}$ is recovered from the splitting operation: $\vec{u} = \vec{u}^*- \nabla \chi$.
Substituting this latter expression into Eq.~\eqref{spaceprob}, we obtain the expression to compute the pressure term $p$ from the auxiliary function: $  p = \chi - \alpha \nabla^2 \chi$. 

One of the main aspects to consider in the projection method is the choice of boundary conditions for $\vec{u}^*$ in Eq.~\eqref{proj} and $\chi$ in the elliptic equation, in order to maintain the desired accuracy. The na\"{i}ve choice $\left. \vec{u}^* \right|_{\partial \Omega} =  \vec{u}_b $
does degrade the overall method to first order accuracy, regardless of space/time discretization accuracy, because it is not guaranteed that $\left. \nabla \chi \right|_{\partial \Omega} = 0$. By projecting $\vec{u}^* = \vec{u} + \nabla \chi $ to the boundary $\partial \Omega$ along the normal and tangential directions $\vec{n}$ and $\vec{\tau}$, respectively, we obtain the following boundary conditions for $\vec{u}^*$:
\begin{equation}\label{ustarbc}
 \vec{u}^* \cdot \vec{n} = \vec{u}_b \cdot \vec{n} + \pad{\chi}{\vec{n}},
\qquad
 \vec{u}^* \cdot \vec{\tau} =\vec{u}_b \cdot \vec{\tau}+ \pad{\chi}{\vec{\tau}}.
\end{equation}
Although we can set an arbitrary value for $\displaystyle \pad{\chi}{\vec{n}}$ (because this value will then be prescribed as boundary condition in the elliptic equation for $\chi$),
the value $\displaystyle \pad{\chi}{\vec{\tau}}$ is not known in \eqref{ustarbc} as $\chi$ is computed later.
One possibility is to use an approximation for $\displaystyle \pad{\chi}{\vec{\tau}}$ in \eqref{ustarbc} by extrapolating $\chi$ from previous time steps: $\chi^{(n+1)} = 2 \chi^{(n)} - \chi^{(n-1)}$. 
A detailed discussion on projection and gauge methods can be found in~\cite{brown2001accurate}.

In this paper we adopt a {\it monolithic} approach, i.e.~we solve the simultaneous equations \eqref{spaceprob} as a single block of equations for both unknowns $\vec{u}$ and $p$.
The pressure term $p$ of Eq.~\eqref{spaceprob} is not uniquely defined, as it appears as $\nabla p$, leading to a singular linear system when discretising the equations in space.
The issue is circumvented by augmenting the problem with an additional scalar unknown $\xi \in \R$ and an additional equation for $p$. The problem consists then of finding
$(\vec{u},p,\xi)$ such that

\begin{equation}\label{spaceprobaug}
\begin{aligned}
\vec{u} - \alpha \nabla^2 \vec{u}
+ \nabla p &= \tilde{\vec{f}} \quad \text{ in } \Omega \\
\nabla \cdot \vec{u} - \xi &= 0 \quad \text{ in } \Omega \\
\int_\Omega p \, d \Omega &= 0 \\
\vec{u} &= \vec{u}_b \quad \text{ on } \partial \Omega
\end{aligned}
\end{equation}

\section{Discretization in space}\label{sect:space}
In this section we describe the space discretization of the equations, with particular attention to the treatment of the curved boundary using a fixed grid and a ghost-point approach.
We describe the method in 2D for simplicity, although the approach can be straightforwardly extended to 3D as it does not rely on two-dimensional aspects.

An arbitrary domain $\Omega$ can be represented in different ways. For example, it can be represented explicitly by the parametric equations that describe the boundary as a closed curve $\partial \Omega = \left\{ (x(\gamma),y(\gamma)) \colon \gamma \in (\gamma_m,\gamma_M) \right\}$, or implicitly by a level-set function $\phi(x,y)$ that is negative inside the domain, positive outside, and vanished on the boundary.
In this paper we adopt a level-set approach, which is more convenient to recover additional geometric information that are required from the numerical method proposed in this paper as we will see later, such as normal direction, curvature, or to easily detect whether a point is inside or outside the domain or to approximate its distance to the boundary. In addition, the level-set function can be used in the framework of level-set methods~\cite{osher2006level}, that is widely adopted in the context of PDEs in moving domains as it allows the boundary of a domain to move according to a velocity field $\vec{u}$. Therefore, although the numerical tests provided in this paper model fluid dynamics around moving objects whose motion is assigned {\it a priori},
the numerical method proposed can be extended to the case of moving objects whose motion is implicitly generated by the fluid dynamics, as for example in fluid-structure interaction, bubble dynamics, wind turbines.

\subsection{Level-set function}\label{sec:levelset}
In this paper we simulate the motion of an incompressible viscous fluid around a (possibly moving) object $\mathcal{R}(t)$. The fluid domain is then $\Omega(t) = \mathcal{D} \backslash \mathcal{R}(t)$, where $\mathcal{D} = (-1,1)^2$ is the computational domain.
We simplify the notation with $\mathcal{R}$ and $\Omega$ instead of $\mathcal{R} (t)$ and $\Omega (t)$ when the object is steady or when it is not necessary to specify the domains as a function of time.

The object $\mathcal{R}$ is implicitly described by a level-set function $\phi \colon D \to \R$ such that:
\[
\mathcal{R} = \left\{ \vec{x} \in \R^2 \colon \phi(\vec{x}) > 0 \right\}, \quad \Gamma = \partial \mathcal{R} = \left\{ \vec{x} \in \R^2 \colon \phi(\vec{x}) = 0 \right\}.
\]
We observe that we have chosen $\phi$ to be positive inside the object so that it is negative inside the fluid domain $\Omega$.
From the level-set function it is possible to infer the geometric properties of the boundary $\Gamma$. For example, the outward unit normal vector $\vec{n}$ and the curvature $\kappa$ can be computed as:
\begin{equation}\label{ncurv}
\vec{n} = \frac{\nabla \phi}{|\nabla \phi|}, \quad \kappa = \nabla \cdot \vec{n}.
\end{equation}
Given an object $\mathcal{R}$, there are infinite level-set functions that can represent that object. A particular level-set function is the signed distance function
$\phi_d(\vec{x}) = \pm \text{dist}(\vec{x},\Gamma)$, where the sign is positive or negative depending on whether $\vec{x}$ is inside or outside $\mathcal{R}$, respectively. The signed distance function is convenient to avoid unstable issues caused by deep ($|\nabla \phi(\vec{x})| \gg 0$) or shallow ($|\nabla \phi(\vec{x})| \ll 0$) gradients that may develop when the level-set evolves over time according to a velocity field (the signed distance function does not present these issues since $|\nabla \phi_d(\vec{x})| = 1$).
At each time step, the signed distance function $\phi_d$ can be recovered from a level-set function $\phi$ by the reinitialization process, that consists of finding the 
steady-state solution (in a fictitious time $\mu$) of the following equation
\begin{equation}\label{sdf}
\frac{\partial \phi_d}{\partial \mu} = \sgn(\phi) \left(1- \left| \nabla \phi_d \right| \right), \quad \phi_d=\phi \; \mbox{ when } \; \mu=0.
\end{equation}
In this paper, we focus for simplicity on predetermined rigid motion of the object $\mathcal{R}$, and therefore we do not need to perform the reinitialization process.
Therefore, we design the numerical methods for a generic level-set function $\phi$, that is not necessarily a signed-distance function.
The extension of the proposed methods to deforming objects would be straightforward.

\subsection{Internal grid points}\label{sect:inteq}
\begin{figure}[H]
 \begin{minipage}[c]{0.49\textwidth}
   	\centering
   	\captionsetup{width=0.80\textwidth}
		\includegraphics[width=0.99\textwidth]{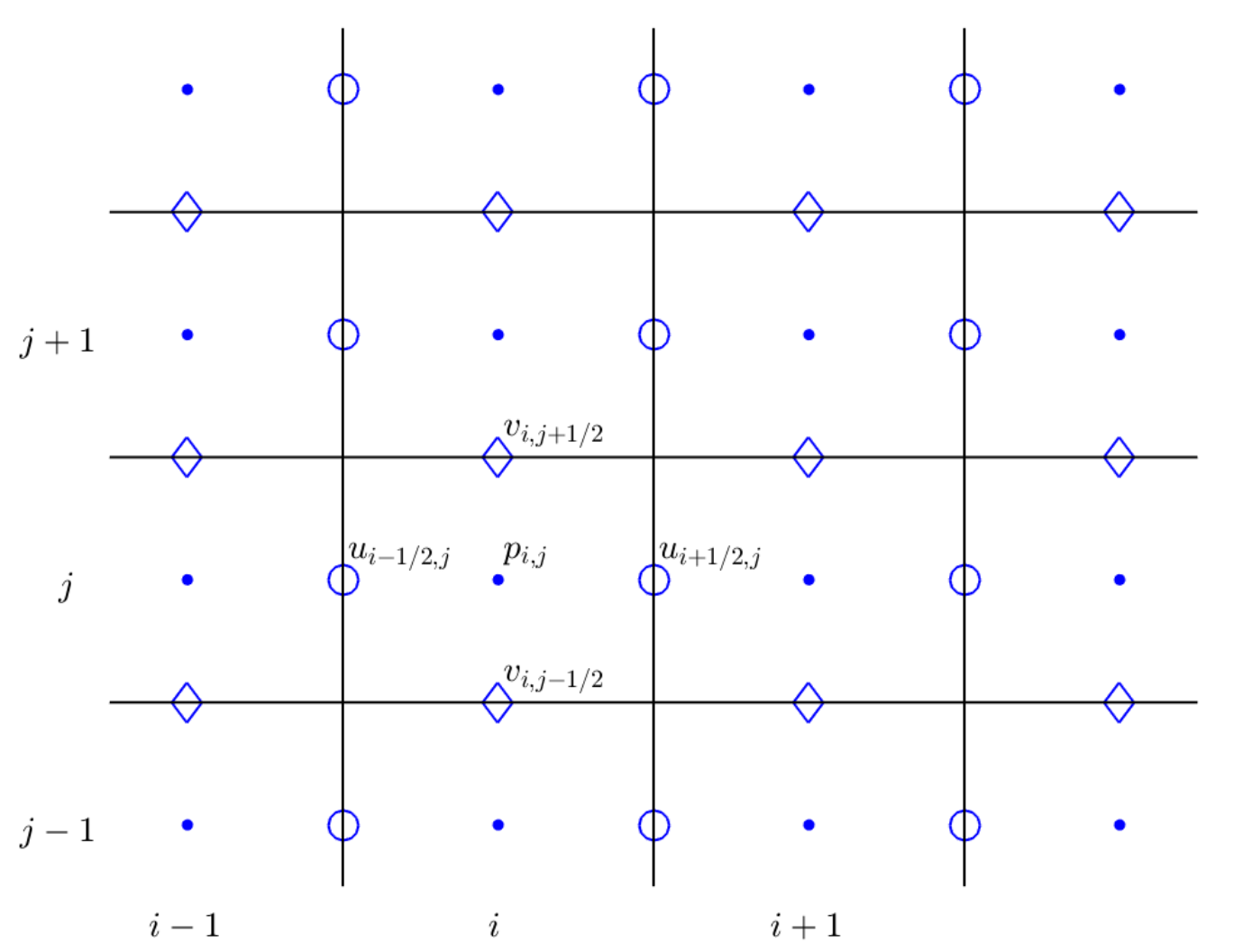}
		\caption{ \footnotesize{ Representation of the staggered grid. $p$ is defined at the centre of each cell ({\it central grid} - small dots), $u$ is defined at the mid points of vertical side cells ({\it horizontally staggered grid} - empty circles) and $v$ at the mid points of horizontal side cells ({\it vertically staggered grid} - diamonds).}}
	\label{fig:staggeredGrid}
 \end{minipage}
 \begin{minipage}[c]{0.49\textwidth}
   	\centering
   	\captionsetup{width=0.80\textwidth}
		\includegraphics[width=0.99\textwidth]{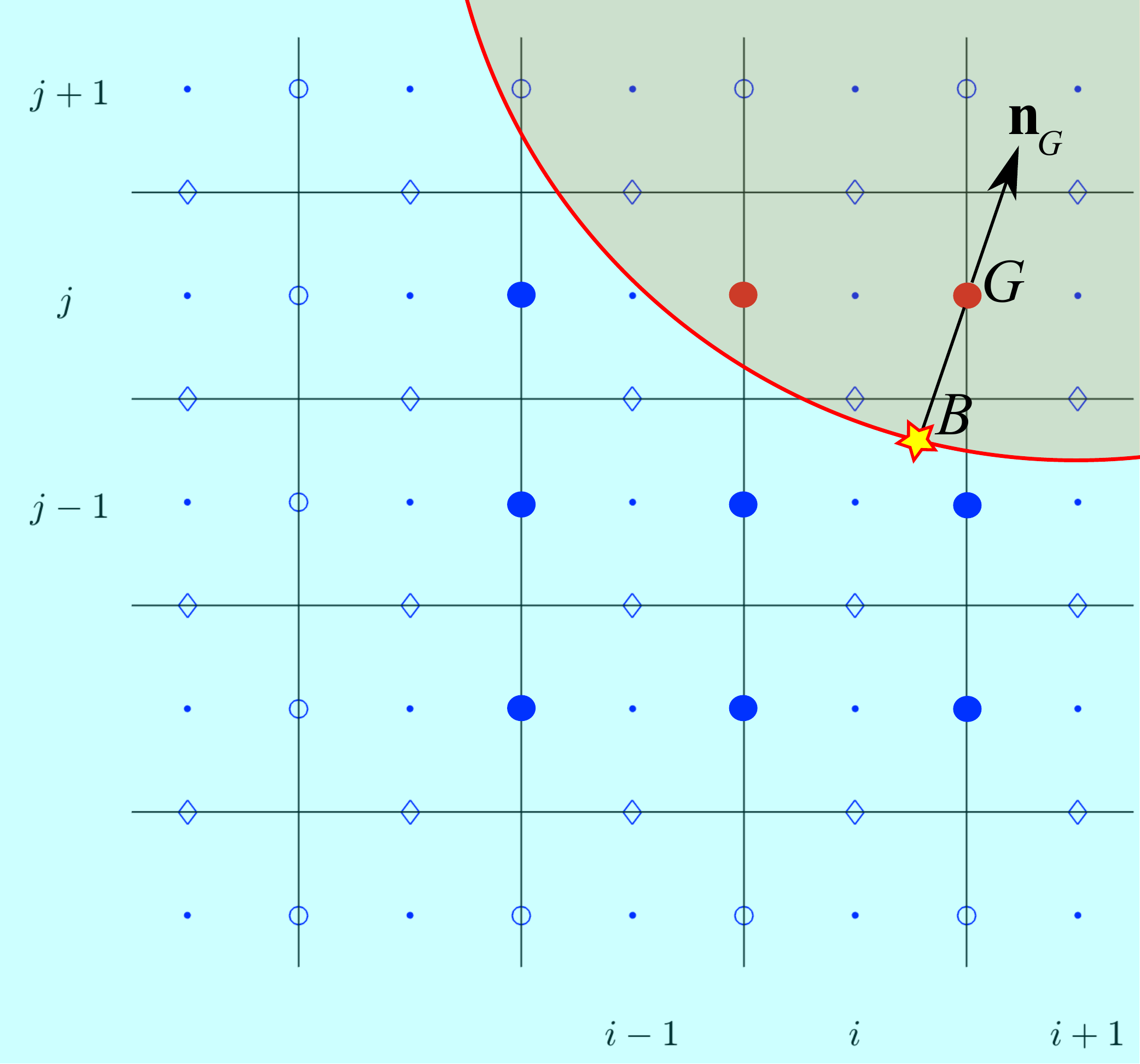}
		\caption{ \footnotesize{ Ghost-point extrapolation. The ghost value $u(G)$ is obtained by enforcing the boundary condition $\tilde{u}(B) = u_b (B)$, where $\tilde{u}(B)$ is the biquadratic interpolation of $u$ on the nine-point stencil represented by the full circles (blue circles inside the fluid and red circles inside the object).}}
	\label{fig:ghost}
 \end{minipage}
\end{figure}

The computational domain $\mathcal{D}$ is discretized through a uniform Cartesian mesh with spatial step $h = \Delta x = \Delta y$. The pressure $p$ and the velocity components $\vec{u}=(u,v)$ are defined on a staggered grid (see Fig.~\ref{fig:staggeredGrid}): $p$ is defined at the centre of each cell ({\it central grid}), $u$ is defined at the mid points of the vertical side cells ({\it horizontally staggered grid}) and $v$ at the mid points of the horizontal side cells ({\it vertically staggered grid}). 
The first two governing equations of~\eqref{spaceprob}, expanded component by component, are:
\begin{equation}\label{spaceprobaugexp}
\begin{aligned}
u - \alpha \nabla^2 u 
+ \pad{p}{x} &= \tilde{{f}}^u \quad \text{ in } \Omega \\
v - \alpha \nabla^2 v
+ \pad{p}{y} &= \tilde{{f}}^v \quad \text{ in } \Omega \\
\nabla \cdot \vec{u} &= \xi \quad \text{ in } \Omega.
\end{aligned}
\end{equation}
We discretize Eqs.~\eqref{spaceprobaugexp} by central finite differences: $\displaystyle \pad{p}{x}$, $ \displaystyle \pad{u}{x}$, $\displaystyle \pad{u}{y}$ and $\nabla^2 u$ are defined on the horizontally staggered grid, $\displaystyle\pad{p}{y}$, $\displaystyle\pad{v}{x}$, $\displaystyle\pad{v}{y}$ and $\nabla^2 v$ are defined on the vertically staggered grid, and $\nabla \cdot \vec{u}$ on the central grid.
\[
\left. \pad{p}{x} \right|_{i+1/2, j} = \frac{p_{i+1,j}-p_{i,j}}{h},
\quad
\left. \pad{p}{y} \right|_{i, j+1/2} = \frac{p_{i,j+1}-p_{i,j}}{h}
\]
\[
\left. \pad{u}{x} \right|_{i+1/2, j} = \frac{u_{i+3/2,j}-u_{i-1/2,j}}{2h},
\quad
\left. \pad{u}{y} \right|_{i+1/2, j} = \frac{u_{i+1/2,j+1}-u_{i+1/2,j-1}}{2h},
\]
\[
\left. \pad{v}{x} \right|_{i, j+1/2} = \frac{v_{i+1,j+1/2}-v_{i-1,j+1/2}}{2h},
\quad
\left. \pad{v}{y} \right|_{i, j+1/2} = \frac{v_{i,j+3/2}-v_{i,j-1/2}}{2h},
\]
\[
\left. \nabla^2 u \right|_{i+1/2, j} = \frac{u_{i+3/2,j}+u_{i-1/2,j}+u_{i+1/2,j+1}+u_{i+1/2,j-1}-4 u_{i+1/2,j}}{h^2},
\]
\[
\left. \nabla^2 v \right|_{i, j+1/2} = \frac{v_{i+1,j+1/2}+v_{i-1,j+1/2}+v_{i,j+3/2}+v_{i,j-1/2}-4 v_{i,j+1/2}}{h^2},
\]
\[
\left. \nabla \cdot \vec{u} \right|_{i,j}=  \frac{u_{i+1/2,j}-u_{i-1/2,j}+v_{i,j+1/2}-v_{i,j-1/2}}{h}.
\]
The discretization of~\eqref{spaceprobaugexp} occurs only on grid points inside of $\Omega$ ({\it internal grid points}). Each grid is treated separately for the identification of internal grid points. In detail, the first equation of~\eqref{spaceprobaugexp} is discretized on internal grid points of the {\it horizontally staggered grid}, the second equation of~\eqref{spaceprobaugexp} is discretized on internal grid points of the {\it vertically staggered grid}, and the third equation of~\eqref{spaceprobaugexp} is discretized on internal grid points of the {\it central grid}.

The third equation of \eqref{spaceprobaug} is discretized by standard mid-point rule, leading to the linear equation
\begin{equation}\label{eq:midpoint}
\sum_{(x_i,y_j) \text{ internal grid point}} p_{i,j} = 0.
\end{equation}

\subsection{Ghost points}\label{sect:ghost}
For internal grid points that are close to the boundary of $\Omega$ it may happen that some of the grid points needed in the discretization stencil may lie outside of the domain. For example, referring to Fig.~\ref{fig:ghost}, the discretization of $\nabla^2 u$ on the internal grid point $(x_{i+1/2},y_{j-1})$ requires the values of $u$ on the external grid point $(x_{i+1/2},y_{j})$, and the discretization of $\nabla^2 v$ on the internal grid point $(x_{i-1},y_{j-1/2})$ requires the values of $v$ on external grid points $(x_{i-1},y_{j+1/2})$ and $(x_{i},y_{j-1/2})$. These external grid points are called {\it ghost points} and a fictitious value (ghost value) of $u$, $v$ and $p$ on these points is required by the discretization of the governing equation on internal grid points. This ghost value is assigned by enforcing the boundary conditions up to the desired accuracy, following a similar approach adopted in~\cite{CocoRusso:Elliptic, coco2018second} for elliptic equation and described below. Let $G=(x_{i+1/2},y_j)$ be a ghost point of the horizontally staggered grid (the approach is similar for any grid), as depicted in Fig.~\ref{fig:ghost}.
Firstly, we compute
the boundary point $B=(x_B,y_B)$ from $G$ by making use of the level-set function as follows. Starting from $G$, we move along the normal direction $\vec{n} = \nicefrac{\nabla \phi}{\left| \nabla \phi \right|}$ by steps of $h/2$, until we obtain an internal point $H$. Let $K$ be the last external point obtained from this process. Then, we perform a bisection method to solve the nonlinear equation $\phi(B) = 0$ in the unknown $B$ along the segment $\overline{HK}$ up to a desired accuracy $\epsilon$ (we choose $\epsilon = 0.01 h$).
For the evaluation of $\phi$ and $\nabla \phi$ we use a bilinear interpolation on $(2 \times 2)$-point stencils containing the point in which the evaluation is required.

The algorithm is summarised in Alg.~\ref{alg:GtoB}.

\begin{algorithm}[H]
\begin{algorithmic}
\State $H \coloneqq G$
\State $\phi_H \coloneqq \phi(G)$
\While{$\phi_H>0$}
\State $K \coloneqq H$
\State let $\tilde{\phi}$ be a bilinear interpolation of $\phi$ on a ($2\times2$)-point stencil containing $H$
\State $\vec{n}_H \coloneqq \displaystyle \left. \left( \nicefrac{ \nabla\tilde{\phi}}{ \left| \nabla \tilde{\phi} \right|} \right) \right|_{H}$
\State $H\coloneqq H-0.5\, h \, \vec{n}_H$
\State $\phi_H \coloneqq \tilde{\phi}(H)$, 
\EndWhile
\For{$\text{iter}=1 \text{  to  } \lceil \log_2 (\text{dist}(H,K)/\epsilon) \rceil$} 
\State $B = \displaystyle \frac{H+K}{2}$
\State $\phi_B \coloneqq \tilde{\phi}(B)$, where $\tilde{\phi}$ is a bilinear interpolation of $\phi$ on a ($2\times2$)-point stencil containing $B$
  \If{$\phi_B \cdot \phi_H < 0$}
    \State $K \coloneqq B$
    \State $\phi_K \coloneqq \phi_B$
  \Else
    \State $H \coloneqq B$
    \State $\phi_H \coloneqq \phi_B$      
\EndIf
\EndFor
\end{algorithmic}
\caption{Bisection method to find the projection boundary point $B$ from a ghost point $G$.}
\label{alg:GtoB}
\end{algorithm}

Once the boundary point $B$ is computed, the ghost value $u_{i+1/2,j}$ is implicitly defined by enforcing the boundary condition $u(B)=u_b(B)$, where $u(B)$ is approximated by interpolating the grid values of $u$ on a suitable stencil. Existing approaches require that the stencil is made only by internal grid points, as the numerical solution is not defined on external points~\cite{costa2019very}. However, we follow the approach of the author and Russo in~\cite{CocoRusso:Elliptic, coco2018second}, where the stencil may contain external grid points, leading to a more flexible scheme.
In detail, we identify the nine-point stencil in the opposite direction to the normal direction $\vec{n}_G$ (upwind direction):
\[
\text{ST}_9 =
\left\{ (x_{i+1/2},y_j) -h\,(s_x \, k_x, s_y \, k_y), (k_x,k_y) \in \left\{ 0,1,2 \right\}^2 \right\},
\]
where $s_x=\text{SIGN} (n_x)$ and $s_y=\text{SIGN} (n_y)$, with $\text{SIGN}(z)=1$ if $z \geq 0$ and $\text{SIGN}(z)=-1$ if $z < 0$ (we take conventionally $\text{SIGN}(0)=1$).
The value $u_{i+1/2,j}$ is implicitly defined by the equation $\tilde{u} (B) = u_b(B)$, 
where $\tilde{u}(B)$ is the biquadratic interpolation of $u$ on the stencil $\text{ST}_9$ evaluated at the boundary point $B$. In detail, let $\theta_x = s_x\left( x_B-x_{i+1/2} \right) / h$ and
$\theta_y = s_y \left( y_B-y_{j} \right) / h$. Then, the equation for $\tilde{u} (B) = u_b(B)$ can be represented using a convenient matrix notation:
\[
 \begin{bmatrix}
u_{i+1/2,j} & u_{i+1/2-s_x,j} & u_{i+1/2-2s_x,j}  \\
u_{i+1/2,j-s_y} & u_{i+1/2-s_x,j-s_y} & u_{i+1/2-2s_x,j-s_y}  \\
u_{i+1/2,j-2s_y} & u_{i+1/2-s_x,j-2s_y} & u_{i+1/2-2s_x,j-2s_y}  
 \end{bmatrix}
 \]
 \[
 \odot \left(
 \begin{bmatrix}
 \displaystyle \frac{(1-\theta_y) (2-\theta_y)}{2} \\
 \\
  \theta_y (2-\theta_y) \\
  \\
  \displaystyle \frac{\theta_y(\theta_y-1)}{2} 
 \end{bmatrix}
 \times 
  \begin{bmatrix}
 \displaystyle \frac{(1-\theta_x) (2-\theta_x)}{2}, \, \theta_x (2-\theta_x), \,  \displaystyle \frac{\theta_x(\theta_x-1)}{2}
 \end{bmatrix}
 \right)
 = u_b(B),
\]
where $A\times B$ is the standard matrix multiplication, and $A \odot B = \sum_{r,s} A_{r,s}B_{r,s}$, or analogously:
\begin{equation}\label{BCu}
\sum_{k_x,k_y=0}^2 \mathcal{T}^{\theta_x}_{k_x} \, \mathcal{T}^{\theta_y}_{k_y} \, u_{i+1/2-k_x s_x,j-k_y s_y} = u_b(B),
\quad \text{ with } 
\mathcal{T}^{\theta} =  \left[ \displaystyle \frac{(1-\theta) (2-\theta)}{2}, \, \theta (2-\theta), \,  \displaystyle \frac{\theta(\theta-1)}{2} \right]^T.
\end{equation}

The reason to take the stencil in the upwind direction is that we want to interpolate values of $u$ from internal grid points and this approach ensures (for smooth boundaries and sufficiently fine grids) that the stencil contains as many internal grid points as possible. However, as in the example of Fig.~\ref{fig:ghost}, the stencil may contain external grid points other than $G$ ($(x_{i-1/2},y_j)$ is another ghost point). These grid points are marked as {\it ghost points} too, and a ghost value for them is defined repeating the entire procedure.

Analogously, we define ghost values for $v$ on ghost points of the vertically staggered grid.

\begin{oss}\label{remark:delta}
The entire set of ghost points needed in the discretization is a subset of all external grid points that are within a predetermined distance $\delta = \sqrt{5} \, h$ from the boundary (referring to the notation of Fig.~\ref{fig:ghost}, the extreme case is when all grid points of the nine-point stencil except $(x_{i+1/2},y_{j-1})$ are outside of the domain and the maximum distance between the external grid points of the stencil and the internal grid point $(x_{i+1/2},y_{j-1})$ is $\delta = \sqrt{5} \, h$). If the grid is sufficiently fine with respect to the curvature of $\Gamma$, this value can be decreased to $\delta=\sqrt{2} \, h$, as farther ghost points are actually not involved in the computation of the numerical solution inside $\Omega$. For example, 
in Fig.~\ref{fig:ghost}, the ghost point $(x_{i+1/2},y_{j+1})$ is not involved in the discretization of any internal equations, neither belongs to the nine-point stencil of any ghost values needed for internal discretizations.
For simplicity, we mark as ghost points all external points within $\delta$ distance from the boundary. Another approach could be to scan all external points and mark as ghost points only those external points whose equation \eqref{BCu} is actually coupled with internal equation, following a similar approach of~\cite{CocoRusso:Elliptic, coco2018second}. We experienced that from a computational cost point of view the two approaches are practically similar.  
\end{oss}

\begin{oss}
Some existing ghost-point methods eliminate the ghost value from the internal equation by substituting its value with the extrapolated value from the boundary condition~\cite{Gibou:Ghost, Shortley-Weller:discretization}. This approach cannot be proposed here, as Eq.~\eqref{BCu} may contain more than one ghost point and therefore the elimination of ghost values from the internal equation is not practical. This approach is called {\it non-eliminated ghost values}.
\end{oss}

\begin{oss}
Existing ghost-point methods allow the definition of more ghost values on the same ghost point $G$ if this is required from more than one internal equation~\cite{Gibou:Ghost, costa2019very}. For example, looking at Fig.~\ref{fig:ghost}, the ghost value $u_{i-1/2,j}$ is required for the discretization of the internal equations at both $(x_{i-1/2},y_{j-1})$ and $(x_{i-3/2},y_{j})$. The approach described in this paper, which is proposed by the author and Russo in~\cite{CocoRusso:Elliptic, coco2018second}, define one single ghost value for each ghost point. 
\end{oss}

Finally, internal equations and ghost interpolations \eqref{BCu} for $u$, $v$ and $p$ constitute a linear system in the unknowns $u_{i+1/2,j}$, $v_{i,j+1/2}$ and $p_{i,j}$, defined on all grid points of the horizontally staggered grid, vertically staggered grid and central grid, respectively.
In summary, the linear system is obtained by writing a linear equation for each grid point $P$ as described in the Algorithm \ref{alg:LS}. 
\begin{algorithm}
\begin{algorithmic}
\For{$P$ in the horizontally-staggered, vertically-staggered and central grids}
  \If{$\phi(P)<0$}
    \State mark $P$ as internal grid point;
    \State the linear equation is obtained by discretizing the three PDEs \eqref{spaceprobaugexp} (first, second or third equation whether $P$ is in the horizontally-staggered, vertically-staggered or central grid, respectively), using the discretization approach described in Sect.~\ref{sect:inteq}
  \ElsIf{$\text{dist}(P,\Gamma) < \delta$}
    \State mark $P$ as ghost point;
    \State the linear equation is obtained from \eqref{BCu} (boundary condition for $u$ or for $v$ whether $P$ is in the horizontally-staggered or vertically-staggered grid, respectively)
  \Else
    \State mark $P$ as inactive grid point;
    \State the linear equation is not needed (or it can be set up fictitiously as, for example, $u_{i+1,j}=0$.)
  \EndIf
\EndFor
\State Write the linear equation \eqref{eq:midpoint}.
\end{algorithmic}
\caption{In this algorithm we summarise the steps to be taken to write the linear system arising from the discretization of \eqref{spaceprobaug}.}
\label{alg:LS}
\end{algorithm}
See Fig.~\ref{fig:ghostBand} for an example of internal, ghost and inactive grid points. We notice that the distance $\text{dist}(P,\Gamma)$ can be easily computed if the signed distance function is available, as $\text{dist}(P,\Gamma) = \phi_d (P)$. Otherwise, $\text{dist}(P,\Gamma) = d(P,B)$, where $B$ is the boundary point obtained by Alg.~\ref{alg:GtoB}.

\subsection{Moving domains}\label{sec:moving}
\begin{figure}[H]
 \begin{minipage}[c]{0.49\textwidth}
   	\centering
   	\captionsetup{width=0.80\textwidth}
		\includegraphics[width=0.79\textwidth]{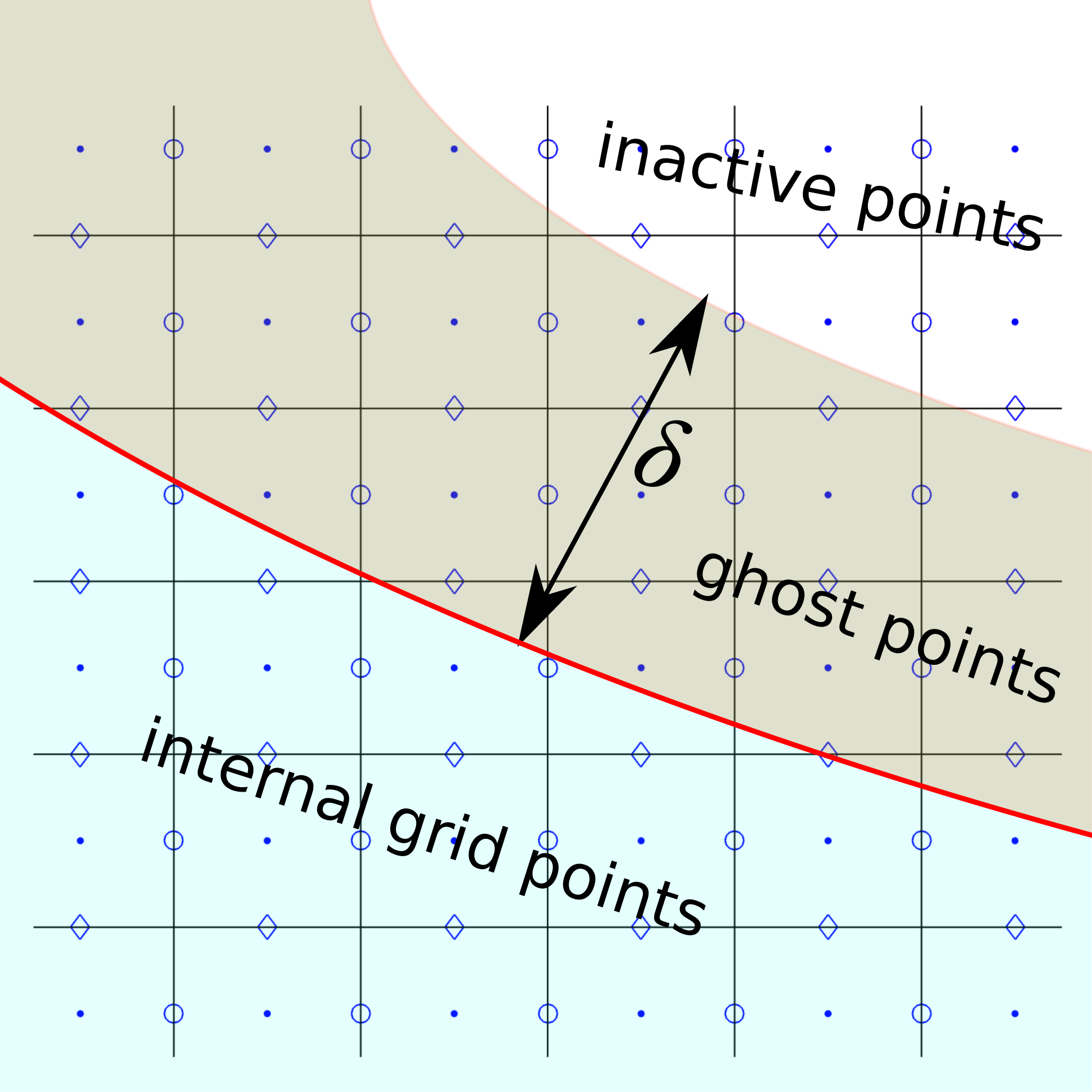}
		\caption{ \footnotesize{ The layer of ghost points (width $\delta$) separate the internal grid points from the inactive grid points. }}
	\label{fig:ghostBand}
 \end{minipage}
 \begin{minipage}[c]{0.49\textwidth}
   	\centering
   	\captionsetup{width=0.80\textwidth}
		\includegraphics[width=0.79\textwidth]{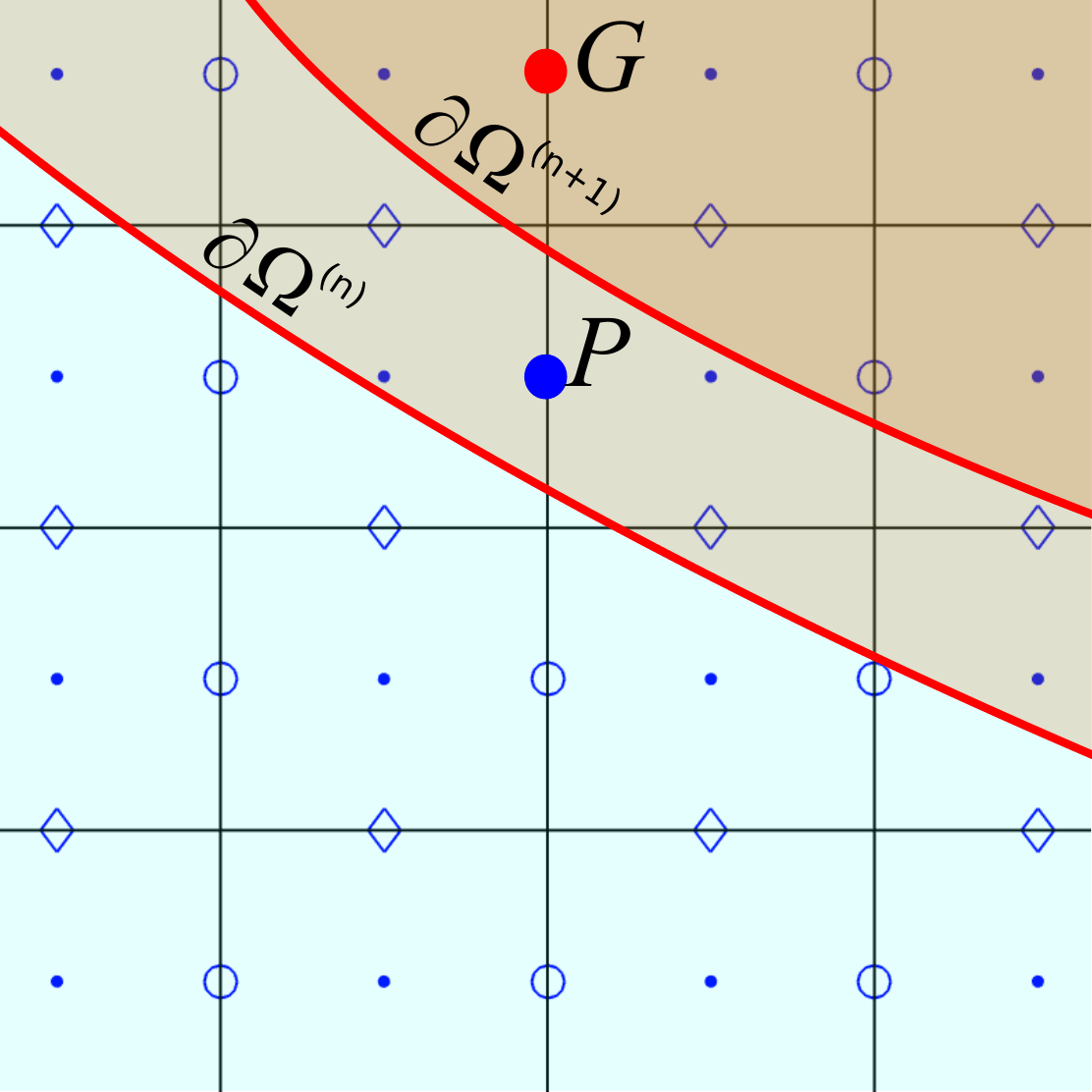}
		\caption{ \footnotesize{ The boundary $\partial \Omega$ moves from $\partial \Omega^{(n)}$ to $\partial \Omega^{(n+1)}$ in one time step. The grid point $P$ is internal to the domain $\Omega^{(n+1)}$ and the discretization of the equations in $P$ requires the value on $G$ at the previous time step. In order to guarantee that a ghost value was computed in $G$ at the previous time step, we have to increase (at each time step) the width of the narrow band of ghost points from $\delta = \sqrt{5} \, h$ to $\delta = (\sqrt{5}+1)h$ (or from $\delta = \sqrt{2} \, h$ to $\delta = (\sqrt{2}+1)h$ if the grid is sufficiently fine with respect to the curvature of the boundary).}}
	\label{fig:ghost_moving}
 \end{minipage}
\end{figure}

The description of the numerical method proposed so far is designed for steady domains $\Omega^{(n)} = \Omega$, $\forall n \in \N$.
In case of moving domain, there are some modifications to be considered. We observe that the explicit term $\tilde{\vec{f}}$ in Eq.~\eqref{spaceprob} 
\[
\tilde{\vec{f}} = \vec{u}^{(n)} + \displaystyle \frac{\Delta t}{2\,\text{Re}} \nabla^2 \vec{u}^{(n)} - \Delta t \left( \vec{u} \cdot \nabla \vec{u} \right)^{(n+1/2)} + \Delta t \, \vec{f}^{(n+1/2)}
\]
must be computed on internal grid points of $\Omega^{(n+1)}$.
On the other hand, $\tilde{\vec{f}}$ contains values of $(\vec{u},p)$ and of its derivatives evaluated at the previous time step, and therefore defined in $\Omega^{(n)}$, that may differ from $\Omega^{(n+1)}$. In particular, an issue arises when an internal grid point for $\Omega^{(n+1)}$ was an external point for $\Omega^{(n)}$, as in Fig.~\ref{fig:ghost_moving}, since central differences centred on ghost points of $\Omega^{(n)}$ may involve grid points that are too far from the boundary and for which there is no value of $(\vec{u},p)$ available. This issue can be circumvented by enlarging the region of ghost points, i.e.~by increasing the value of $\delta$ in the Algorithm $\ref{alg:LS}$ and in the Remark \ref{remark:delta}. In order to mitigate this effect, we ensure that the boundary does not move more than $h$ by providing a suitable restriction for the time step $\Delta t$, that is $\displaystyle \Delta t < h/\max_{\partial \Omega^{(n)}} \left\|\vec{u}^{(n)}_b \right\|_2$, and then we set up $\delta = (\sqrt{5}+1)h$.

\begin{oss}
Following the Remark \ref{remark:delta}, 
if the grid is sufficiently fine with respect to the curvature of the boundary, the value of $\delta$ can be decreased to $\delta = (\sqrt{2}+1) \, h$.
\end{oss}
\begin{oss}
The reverse case in which an internal grid point for $\Omega^{(n)}$ becomes a ghost point for $\Omega^{(n+1)}$ does not give any warnings, as the equation for that grid point is set up normally as described in Sect.\ref{sect:ghost} enforcing boundary conditions on $\Omega^{(n+1)}$.
\end{oss}

\section{Multigrid approach}\label{sect:mg}
The linear system obtained in Sect.~\ref{sect:space} and described in Algorithm~\ref{alg:LS} can be summarised as
\begin{equation}\label{eq:LS}
\begin{pmatrix}
I^u_h-\alpha L^u_h & \vec{0} &  (D_x)^u_h & \vec{0}    \\
B^u_h & \vec{0} &  \vec{0} & \vec{0}    \\
\vec{0}  & I^v_h-\alpha L^v_h &  (D_y)^v_h & \vec{0}  \\
\vec{0}  & B^v_h &  \vec{0} & \vec{0}  \\
(D_x)^p_h  & (D_y)^p_h &  \vec{0} & -{1}  \\
\vec{0}  & \vec{0} &  \vec{1}_h & \vec{0} 
\end{pmatrix}
\cdot 
\begin{pmatrix}
{u}_h   \\
{v}_h  \\
{p}_h  \\
\xi
\end{pmatrix}
=
\begin{pmatrix}
\tilde{{f}}^u_h   \\
u_b  \\
\tilde{{f}}^v_h  \\
v_b  \\
\vec{0}  \\
0  
\end{pmatrix},
\quad
\text{or in a compact form:}
\quad
M_h\vec{u}_h = \vec{b}_h,
\end{equation}
where $I^{u,v}_h$ is the identity matrix, $L^{u,v}_h$ is the discretization of the Laplacian operator, $B^{u,v}_h$ the discrete boundary condition operators, $(D_{x,y})^{u,v,p}_h$ the central difference operartor for the first derivative, and $\vec{1}_h$ a vector of ones.

This is a sparse non-symmetric linear system, whose solution is obtained efficiently by adopting a geometric multigrid approach. Multigrid methods are efficient methods to solve sparse linear systems arising from the discretization of a class of PDEs.
The aim of a multigrd method is to solve the high-frequency components of the solution on a fine grid and the low-frequency components on coarser grids~\cite{Trottemberg:MG, hackbusch2013multi, briggs2000multigrid}.
Initially designed for elliptic equations, their implementation has been extended to PDEs that have similarities with elliptic equations, in particular those (systems of) equations whose discretization presents a high {\it ellipticity} factor, defined as a quantity to establish the qualitative resemblance of the discrete operator to the elliptic operator (see~\cite[Ch.~4]{Trottemberg:MG}.
In the recent decades, they have also been designed to several other classes of PDEs, including parabolic~\cite{gander2016analysis} and hyperbolic equations~\cite{ruggiu2019multigrid}, as well as to general sparse linear systems whose derivation is not (or is partially) known, leading to the class of {\it Algebraic Multigrid Methods}~\cite{stuben2001review} (in contrast with geometric multigrid methods), where geometric information on the grid or on the relation between the linear equations is not available.

In this paper we adopt a geometric multigrid method specifically designed for ghost-point methods, extending the idea proposed in~\cite{CocoRusso:Elliptic, coco2018second} and combining the approach with traditional multigrid methods developed for Navier-Stokes equations for rectangular domains. 
We confine our explanation to the main components of the multigrid approach, i.e.~relaxation operator, restriction operation and interpolation operator, from which the entire multigrid algorithm is uniquely defined. For a comprehensive explanation of multigrid methods see, for example,~\cite{Trottemberg:MG, hackbusch2013multi, briggs2000multigrid}.

\subsection{Relaxation operator}
The most important property that a relaxation operator must satisfy for a multigrid performance is the smoothing property, defined as the ability of the relaxation scheme to smooth (i.e.~quickly dump high-frequency components of) the residual $\vec{r}_h=\vec{b}_h-M_h \vec{\tilde{u}}_h$, where $\vec{\tilde{u}}_h$ is the approximation of the solution obtained after few relaxation steps. The relaxation operator is also called smoother. For elliptic equations discretized by standard central-differences, the Gauss-Seidel scheme is an appropriate smoother, as well as weighted-Jacobi scheme, with weight $\omega_\text{Jac}=4/5$ in 2D, while non-weighted Jacobi does not have the smoothing property. Several types of Gauss-Seidel schemes can be defined, depending on the order chosen to swap the grid points, i.e.~the equations of the linear systems. Common choices are lexicographic order and red-black. Although the latter is more efficient and suitable for parallelization, we focus in this paper on the Gauss-Seidel lexicographic order (GS-LEX) for simplicity, keeping in mind that the ideas can be straightforwardly extended to other Gauss-Seidel schemes.

The linear system~\eqref{eq:LS} is made by internal equations and ghost equations. The smoother will swap onto all grid nodes, acting in a different way whether it is on internal grid points or on ghost grid nodes, as described in the next sections.

\subsubsection{Internal smoother}\label{sec:internalsmoother}
For the internal equations, we follow the approach described in cite~\cite[Sect.~8.7]{Trottemberg:MG} and briefly summarised here. First of all, a standard Gauss-Seidel scheme does not act as an efficient smoother for systems of PDEs. A proper adjustment to retrieve the smoothing property is that for every single node we update simultaneously all variables defined at that node. This idea leads to the solution of small $n\times n$ systems, where $n$ is the number of variables defined in that grid node. This approach is called {\it collective} relaxation. For the linear system~\eqref{eq:LS}, a collective relaxation is not properly defined, as for example the variables defined in the nodes of the central grid, i.e.~$p_{i,j}$, do not appear in the linear equations associated with those nodes, that are
\[
\frac{u_{i+1/2,j}-u_{i-1/2,j}+v_{i,j+1/2}-v_{i,j-1/2}}{h} = \xi,
\]
or, in other terms, the linear system has zero diagonal entries. A modification to overcome this issue is the {\it box relaxation}: for each cell we update simultaneously all variables defined in that cell, including cell borders. This leads to the solution of $5\times5$ systems for each cell, as the variables defined in the cell with centre $(x_i,y_j)$ are $p_{i,j},u_{i\pm1/2,j},v_{i,j\pm1/2}$.
We observe that variables defined on the border cells are updated twice during one relaxation step, as for example $v_{i,j+1/2}$ is updated for the cells having centres $(x_i,y_j)$ and $(x_i,y_{j+1})$.

\subsubsection{Ghost smoother}\label{sec:ghostsmoother}
The relaxation scheme in the ghost points is defined using a similar approach to the one adopted in~\cite{CocoRusso:Elliptic, coco2018second} for scalar equations. Since a na\"{i}ve Gauss-Seidel approach does not converge for ghost equations, we adopt a fictitious-time relaxation on the boundary conditions. In detail, let $B(w)=g$ represent a generic boundary condition for the variable $w$ (e.g.~$B(w)=w$ for Dirichlet boundary conditions and $B(w)=\partial w / \partial n$ for Neumann boundary conditions). The relaxation scheme is obtained by discretising the associated fictitious-time boundary condition:
\begin{equation}\label{fictBC}
\frac{\partial w}{\partial \tau} + B w = g,
\end{equation}
where $\tau$ is a fictitious time. Observe that the steady-state solution of \eqref{fictBC} is the original boundary condition $B(w)=g$: therefore (provided that we have convergence) the solution of \eqref{fictBC} will approach asymptotically a solution $w$ that satisfies $B(w)=g$ on the boundary. Since the accuracy of the method is not affected by the (fictitious) time discretization, we can discretize the time derivative with first order accuracy, and in particular we use forward Euler method. In addition, we can discretize the time derivative of Eq.~\eqref{fictBC} in a convenient point rather than on a boundary point, for example on the grid points in which $w$ is defined, obtaining:
\begin{equation}\label{ghostrelax}
w_{r,s}^{(m+1)} = {w}_{r,s}^{(m)} + \Delta \tau_w \left( g - B w^{(m)} \right),
\end{equation}
where $B w^{(m)}$ is the discretization in space described in Sect.~\ref{sect:ghost}. The fictitious time step $\Delta \tau_w$ is an additional parameter that must be defined in order to guarantee the convergence of the relaxation scheme. It is observed numerically that the following conditions will guarantee convergence (see~\cite{CocoRusso:Elliptic, coco2018second} for a theoretical explanation):
\[
0< \Delta \tau_w < 1 \text{ for Dirichlet boundary conditions,} \quad
0< \Delta \tau_w < h \text{ for Neumann boundary conditions,}
\]
where $h$ is the spatial step. In practice, we choose $\Delta \tau_w = 0.9$ and $\Delta \tau_w = 0.9 \, h$, respectively.

In order to guarantee the smoothing property and cancel possible boundary effects that would affect the multigrid performance, the relaxation over the ghost points is performed a certain number of times, say $\lambda$, against each step of the internal iterations.

In summary, Algorithm \ref{alg:relax} describes how $\nu$ relaxation steps are performed.
\begin{algorithm}
\begin{algorithmic}
\For{  {\tt iter1} from 1 to $\nu$  }
	\For{  $(x_i,y_j) \in \Omega$}
    \State Box-relaxation over the cell centred in $(x_i,y_j)$
    \EndFor
    \For{  {\tt iter2} from 1 to $\lambda$  }
    \For{  each ghost point of the central, horizontally staggered and vertically staggered grid}
        \State ghost relaxation \eqref{ghostrelax}
        \EndFor
    \EndFor        
\EndFor                
\end{algorithmic}
\caption{Summary of $\nu$ iteration steps of the relaxation scheme.}
\label{alg:relax}
\end{algorithm}


\subsection{Restriction operator}
After few relaxations (say $\nu_1$) on the fine grid (with a spatial step $h$), an approximation $\vec{\tilde{u}}_h$ of the final solution is obtained. Then, the residual $\vec{r}_h=\vec{b}_h-M_h \vec{\tilde{u}}_h$ is computed and transferred (restricted) to a coarser grid, i.e.~a staggered mesh with a spatial step $H=2h$, where the residual equation is solved: $M_H \vec{e}_H = \vec{r}_H$. Therefore, proper restriction operators must be defined. Since
the residual $\vec{r}_h = (\vec{r}^u_h,\vec{r}^v_h,\vec{r}^p_h)$ is defined on the horizontally-staggered, vertically-staggered and central grids (for the residuals $\vec{r}^u_h$ of the horizontal velocity, $\vec{r}^v_h$ of the vertical velocity, and $\vec{r}^p_h$ of the continuity equation, respectively), then we have to define three restriction operators $\vec{r}^u_H = \mathcal{I}_H^{u} (\vec{r}^u_h)$, $\vec{r}^v_H = \mathcal{I}_H^{v} (\vec{r}^v_h)$ and $\vec{r}^p_H = \mathcal{I}_H^{p} (\vec{r}^p_h)$. 

\begin{figure}[H]
   	\centering
   	\captionsetup{width=0.80\textwidth}
		\includegraphics[width=0.39\textwidth]{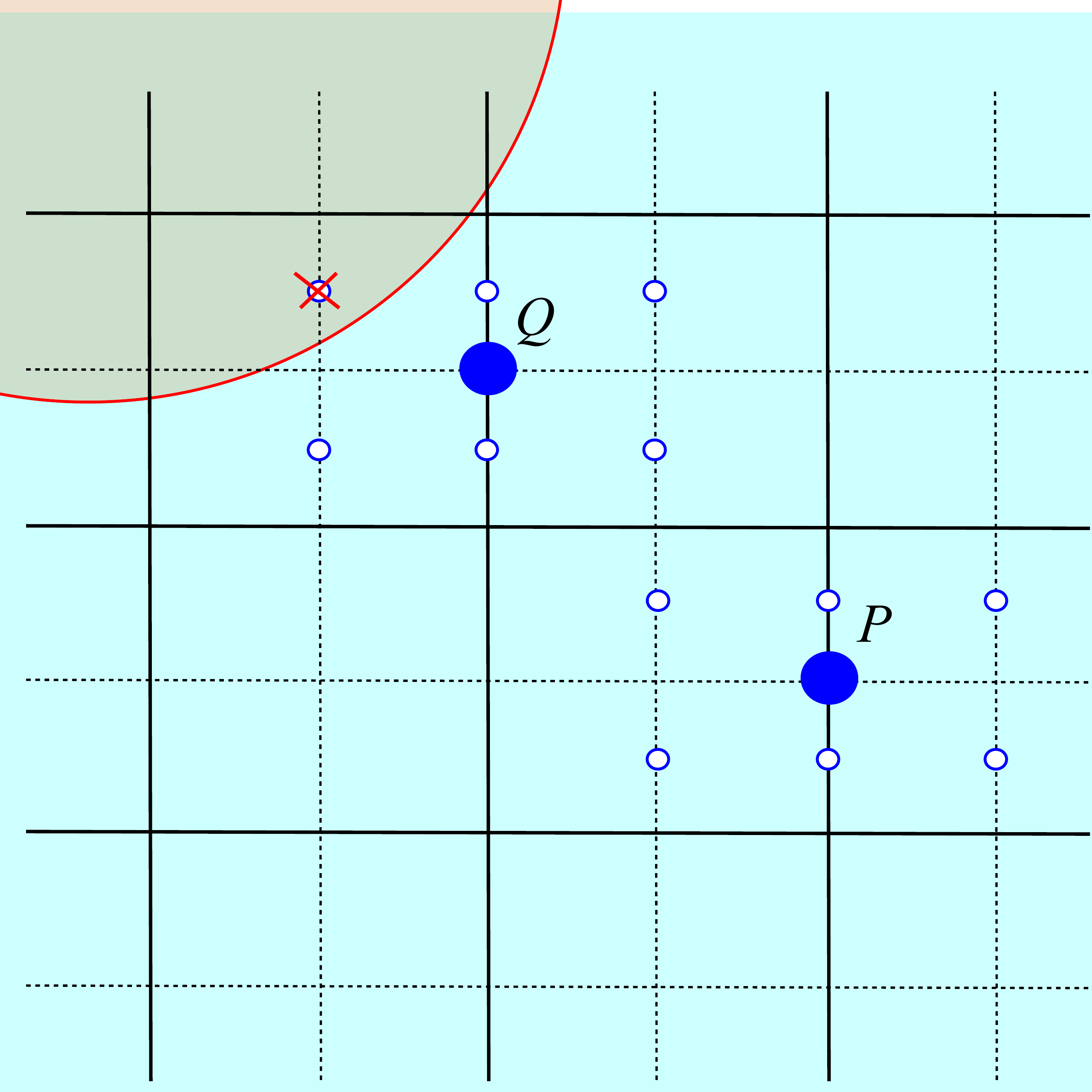}
	\caption{ \footnotesize{Restriction operation $\vec{r}^u_H = \mathcal{I}_H^{u} (\vec{r}^u_h)$ for the horizontally staggered grid. The grid with spatial step $h$ is represented by dashed and bold solid lines, while the grid with spatial step $H=2\,h$ represented by bold solid lines. The restriction value $\vec{r}^u_H (P)$ is obtained by using the standard restriction operator (first stencil of Eq.~\ref{restop}), since the grid points of the stencil are all internal (represented by the six empty small circles surrounding $P$). The restriction value $\vec{r}^u_H (Q)$ is obtained by using the restriction operator of Eq.~\eqref{restop5}, since the grid points of the stencil are all internal except the one at the top-left.}}
	\label{fig:FWrestriction}		
 \end{figure}

The standard restriction operators for staggered grids are (see~\cite[Sect.~8.7.1]{Trottemberg:MG}):
\begin{equation}\label{restop}
\mathcal{I}_H^{u} = \frac{1}{8} 
\begin{bmatrix}
1 & 2 & 1 \\
& \cdot &  \\
1 & 2 & 1 \\
\end{bmatrix}, \quad 
\mathcal{I}_H^{v} = \frac{1}{8} 
\begin{bmatrix}
1 &  & 1 \\
2& \cdot & 2  \\
1 & & 1 \\
\end{bmatrix}, \quad 
\mathcal{I}_H^{p} = \frac{1}{4} 
\begin{bmatrix}
1 &  & 1 \\
& \cdot &  \\
1 &  & 1 \\
\end{bmatrix}.
\end{equation}
An example on the horizontally staggered grid is represented in Fig.~\ref{fig:FWrestriction}.
Since the smoother acting on the internal equations (Sect.~\ref{sec:internalsmoother}) is different from the one acting on the ghost points (Sect.~\ref{sec:ghostsmoother}), the residual that is obtained after few relaxations, although is smooth either internally or along the boundary, might appear discontinuous from inside grid points to ghost points. For this reason, following the idea proposed in~\cite{CocoRusso:Elliptic, coco2018second}, the restriction of the residual for internal grid points is performed without involving the values obtained in the ghost points.

In practice, for each internal grid point of the coarser grid with spatial step $H$, if the grid points of the finer grid with spatial step $h$ that are involved in the restriction operator are all internal grid points, than the usual restriction \eqref{restop} is performed, otherwise the restriction is the average of the residual on the internal grid points of the stencil. For example, assume that a grid point $Q$ of the coarser horizontally-staggered grid is close to the boundary and that the grid points of the finer grid involved in the restriction are all inside $\Omega$ except the point at the top-left (ghost point), then the restriction is (see Fig.~\ref{fig:FWrestriction})
\begin{equation}\label{restop5}
\mathcal{I}_H^{u} = \frac{1}{5} 
\begin{bmatrix}
0 & 1 & 1 \\
& \cdot &  \\
1 & 1 & 1 \\
\end{bmatrix}.
\end{equation}
The restriction of the residual of the ghost equations is performed in the same way, i.e.~only involving ghost values and ignoring the values of the residual of internal equations.

\subsection{Interpolation operator}
Once the restriction of the residual $\vec{r}_H$ is computed and the residual equation $M_H \vec{e}_H = \vec{r}_H$ is solved on the coarser grid, the error $\vec{e}_H$ is interpolated back from the coarser grid to the finer grid $\vec{e}_h$. 
Since the error is the solution of the residual equation, we do not expect any discontinuities from the internal values to the ghost values. Therefore, the standard interpolation operator is adopted, either for internal or ghost values. In practice, for each grid point on the fine grid, we compute the bilinear interpolation on the coarser grid points that surround it (see, for example,~\cite[Sect.~8.7.1]{Trottemberg:MG}).
Therefore, the approximation $\tilde{\vec{u}}_h$ is updated as $\tilde{\vec{u}}_h \colon = \tilde{\vec{u}}_h + \vec{e}_h$, and additional $\nu_2$ relaxations are performed on the fine grid. The final approximation is the result of the Two-Grid Correction Scheme (TGCS). 

When solving the residual equation in the grid with spatial step $H$, one can recursively apply a multigrid iteration again, moving to a coarser grid with spatial step $2H$, and so on until a very coarse grid is reached, on which the (small) linear system is solved exactly by a standard solver. Examples of multigrid approaches based on this strategy are $V-$, $W-$ and $F-$cycles (we refer to~\cite{Trottemberg:MG} for a more comprehensive description).
In this paper, we adopt the $W-$cycle approach, although the technique can be easily extended to other approaches.

\section{Numerical tests}\label{sect:numtests}
\begin{figure}[H]
 \begin{minipage}[c]{0.33\textwidth}
   	\centering
   	\captionsetup{width=0.80\textwidth}
		\includegraphics[width=0.79\textwidth]{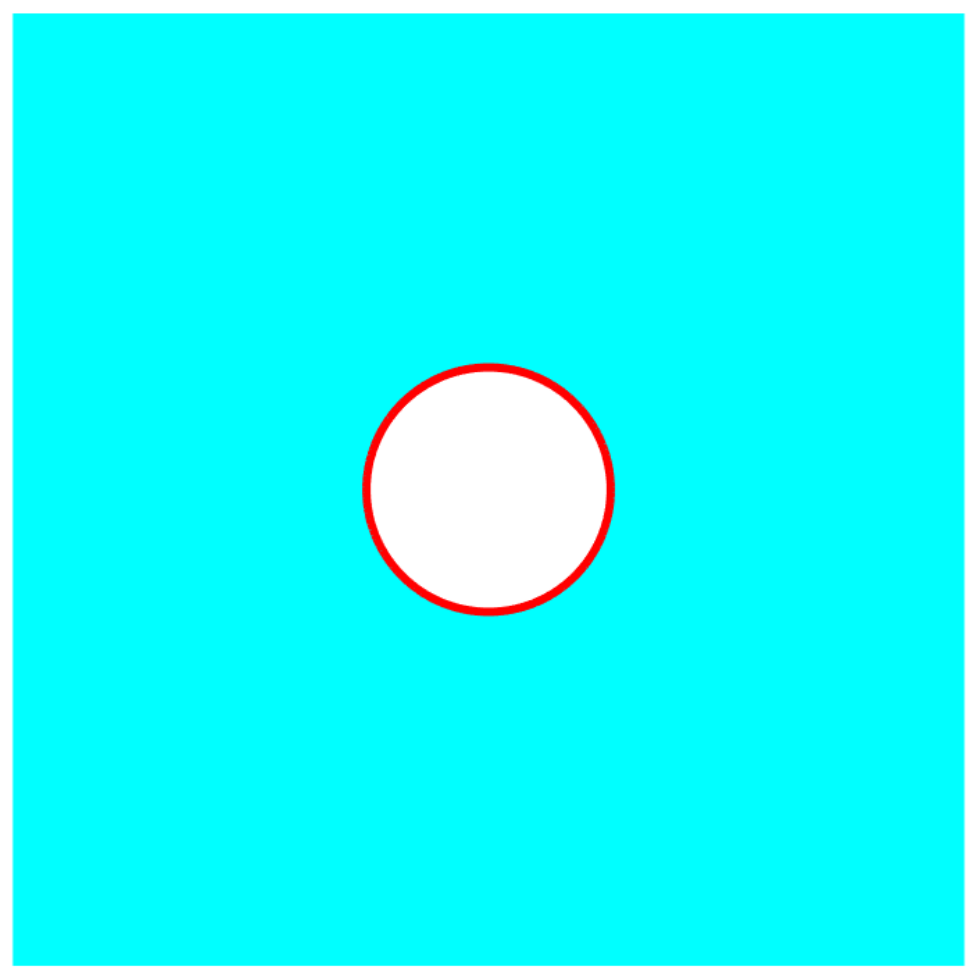}
 \end{minipage}
 \begin{minipage}[c]{0.33\textwidth}
   	\centering
   	\captionsetup{width=0.80\textwidth}
		\includegraphics[width=0.79\textwidth]{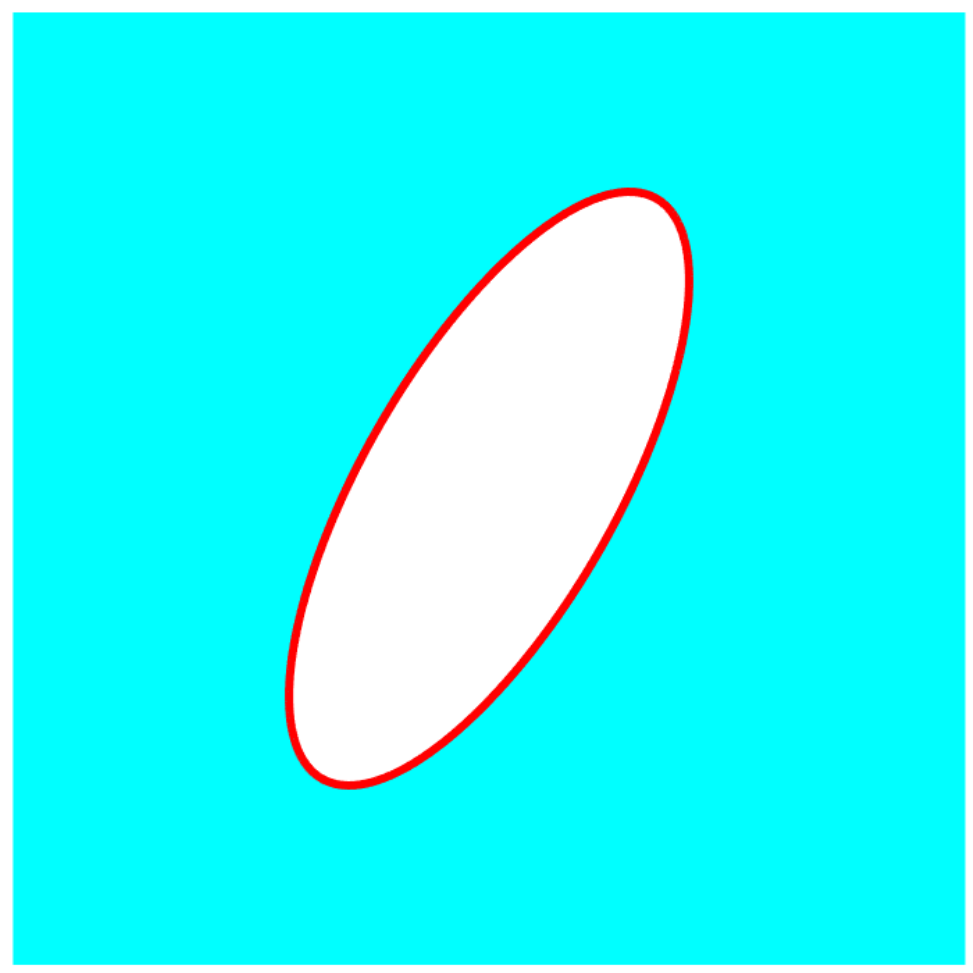}
 \end{minipage}
  \begin{minipage}[c]{0.33\textwidth}
   	\centering
   	\captionsetup{width=0.80\textwidth}
		\includegraphics[width=0.79\textwidth]{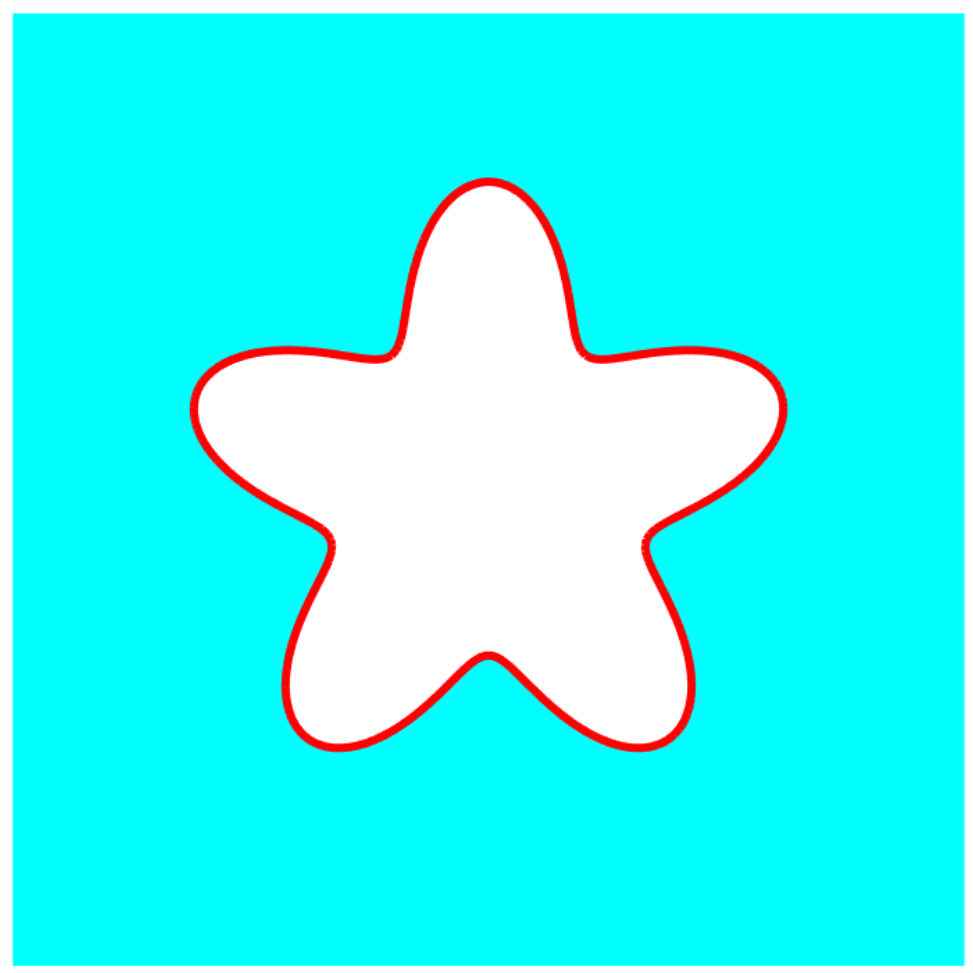}
 \end{minipage}
\caption{ \footnotesize{Different shapes of the object $\mathcal{R}(0)$ adopted in the numerical tests: \textsc{Circle} (left), \textsc{Ellipse} (middle) and \textsc{Flower} (right).} }
	\label{fig:domains}
\end{figure}

In this section we perform several numerical tests to assess the accuracy of the discretization and the efficiency of the multigrid. 
The computational domain is $\mathcal{D} = (-1,1)^2$ and the fluid domain is $\Omega = \mathcal{D} \backslash \mathcal{R}(t)$, where $\mathcal{R}(t)$ represents the object at time $t$, described by a level-set function (Sect.~\ref{sec:levelset}). Different shapes will be tested: circle, ellipse and flower-shaped, defined at time $t=0$ as follows (see Fig.~\ref{fig:domains}):
\begin{equation}\label{eq:domains}
\begin{matrix}
\textsc{Circle:} &
\phi(x,y) = R^2-(x^2+y^2), \quad R = \displaystyle \frac{1}{\sqrt{15}}, \\
\\
\textsc{Ellipse:} &
\phi(x,y) = 1-\left( \displaystyle \frac{(\cos(\theta)x-\sin(\theta)y)^2}{A^2}+ \displaystyle \frac{(\sin(\theta)x+\cos(\theta)y)^2}{B^2} \right), \quad A = \displaystyle \frac{1}{\sqrt{14}}, \; B = \displaystyle \frac{1}{\sqrt{2}}, \;\theta = \displaystyle \frac{\pi}{6},  \\
\\
\textsc{Flower:} &
\phi(x,y) = \sqrt{x^2+y^2}-A-\displaystyle \frac{B (y^5+5x^4y-10x^2y^3)}{\sqrt{x^2+y^2}^5},
\quad A = 0.5, \; B = 0.15.\\
\end{matrix}
\end{equation}
For steady domains, the object will remain steady for the entire simulation, i.e.~$\mathcal{R}(t)=\mathcal{R}(0)$, while for moving domains we rotate the initial object $\mathcal{R}(0)$ around the origin with a constant angular velocity $\omega \in \R$, i.e.~$\mathcal{R}(t)$ is the rotation of $\mathcal{R}(0)$ by an angle $\omega \, t$. The level-set function at time $t$ is computed as $\phi(x,y,t) = \phi(\bar{x},\bar{y},0)$, where $(\bar{x},\bar{y})$ is the point obtained rotating $(x,y)$ by an angle $-\omega \, t$ around the origin, i.e.~
\[
\begin{pmatrix}
\bar{x} \\
\bar{y}
\end{pmatrix}
=
\begin{pmatrix}
\cos (\omega t)  & -\sin (\omega t)\\
\sin (\omega t)  & \cos (\omega t)
\end{pmatrix}
\cdot
\begin{pmatrix}
x\\
y
\end{pmatrix}.
\]

\subsection{Accuracy of the discretization}
In this section we test the accuracy of the method for the different domains presented in \eqref{eq:domains}. We choose an exact solution $\vec{u}_\text{exa}=(u_\text{exa},v_\text{exa},p_\text{exa})$ and adjust the right-hand side of \eqref{spaceprobaug} accordingly. Since the exact solution is not necessary a physical solution, we do not expect in general that it satisfies the continuity equation.
To this purpose, the continuity equation is generalised to $\nabla \cdot \vec{u} - \xi = \vec{g}$, where $\vec{g} = \nabla \cdot \vec{u}_\text{exa}$.
Source term $\vec{f}$ and boundary conditions are defined from the exact solution.

The exact solution that we use in the tests is:
\begin{equation}\label{sol:exa}
u_\text{exa} = \cos(5 x) \cos(6 y \log(t+2)), \;
v_\text{exa} = \sin(4 t) \sin(3 x^2+4y^2+2), \; 
p_\text{exa} = \cos(6 x t) \sin(2 y t) \log(3 t+1)+\sin(5t).
\end{equation}
Reynolds number is $\text{Re}=100$.
The computational domain $\mathcal{D}=(-1,1)^2$ is divided in $N^2$ cells.
As described in Sect.~\ref{sect:disctime}, discretization in time is performed using a Crank-Nicholson approach, which is unconditionally stable. Moreover, for moving domains the time step is chosen in such a way that the object does not move more than the spatial step $h$ at each time step (see the restriction on the time step provided in Sect.~\ref{sec:moving}).
Therefore, the time step chosen for the numerical tests is
\[
\Delta t \leq \displaystyle \min \left\{ h,  \frac{h}{\max \left\| \vec{u} \right\|_2} \right\}.
\]
For each test, we choose several values of $N$ and for each of them we compute the error as the difference between the exact solution and the numerical solution at time $t=10$. We plot the errors against $N$ in logarithmic scale, together with the best fit line to compute an approximation of the accuracy order.
The values are also reported in tables. 
We perform two tests for each geometry of Eq.~\eqref{eq:domains}: one test for the steady case and one test for the moving geometry, with an angular velocity of $\omega = 2 \pi / 5$. We use the notation \textsc{Steady/Moving Geometry}, where \textsc{Geometry} is \textsc{Circle}, \textsc{Ellipse} or \textsc{Flower}.

\begin{figure}[H]
 \begin{minipage}[c]{0.33\textwidth}
   	\centering
   	\captionsetup{width=0.80\textwidth}
		\includegraphics[width=0.99\textwidth]{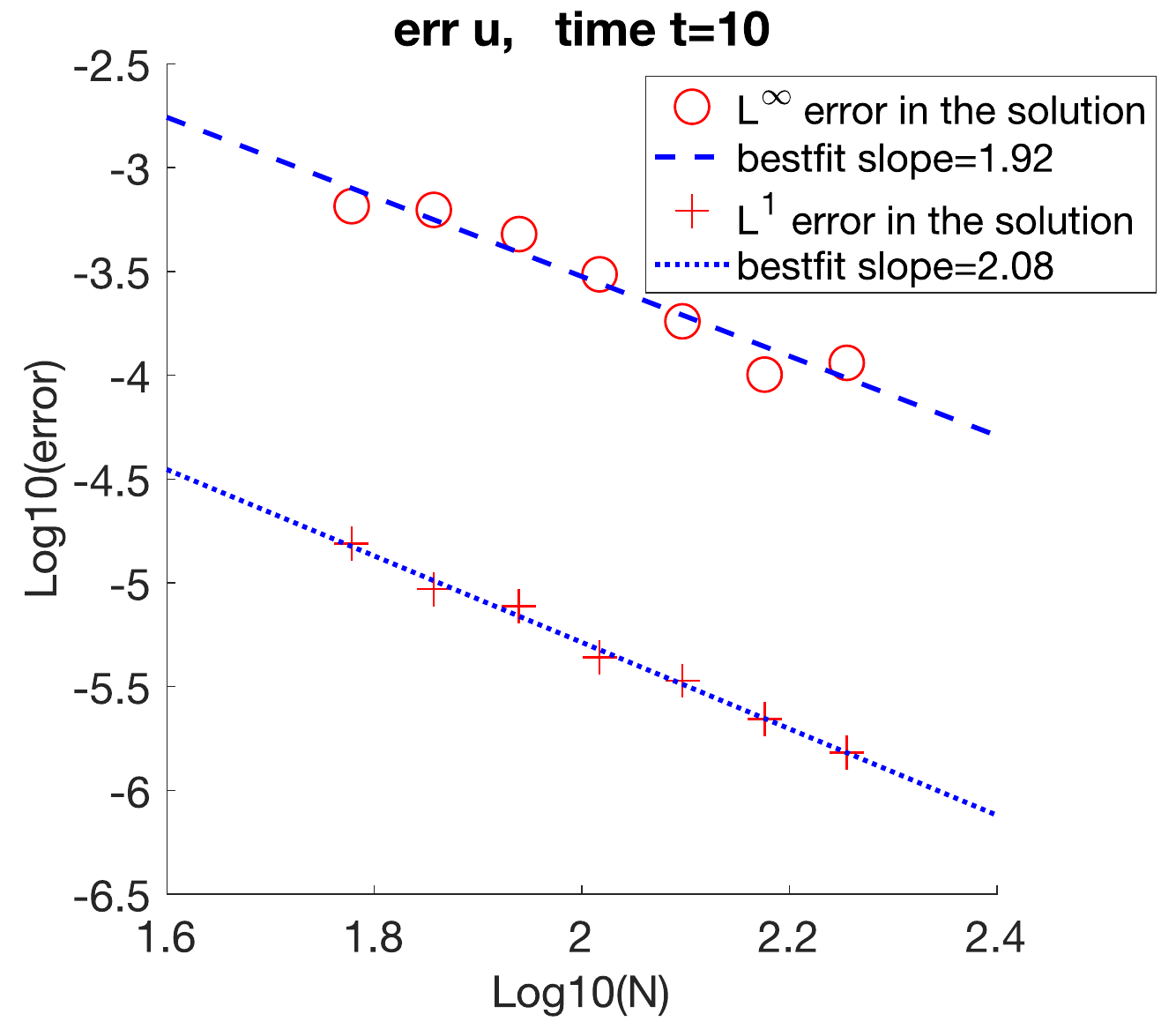}
 \end{minipage}
 \begin{minipage}[c]{0.33\textwidth}
   	\centering
   	\captionsetup{width=0.80\textwidth}
		\includegraphics[width=0.99\textwidth]{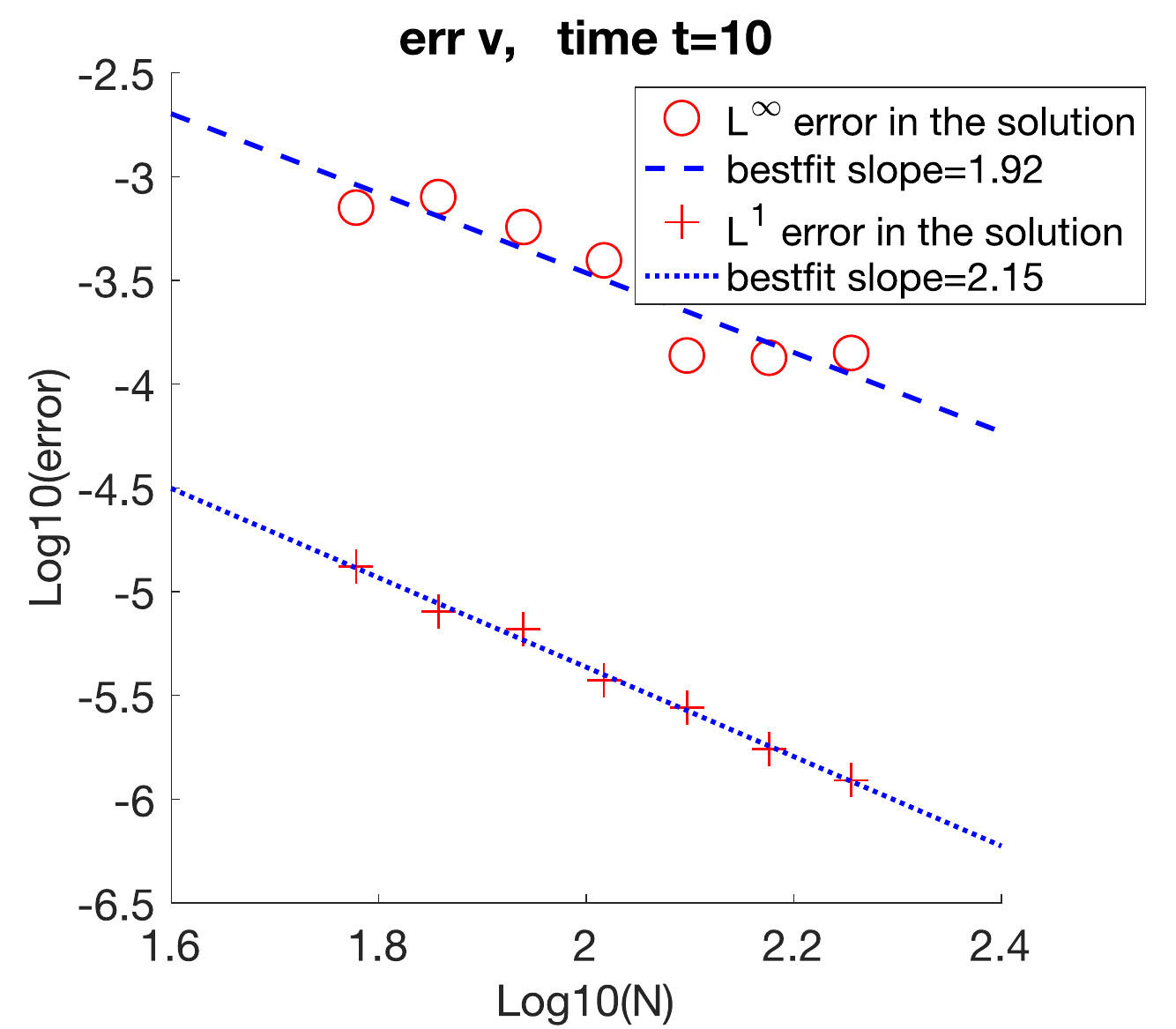}
 \end{minipage}
  \begin{minipage}[c]{0.33\textwidth}
   	\centering
   	\captionsetup{width=0.80\textwidth}
		\includegraphics[width=0.99\textwidth]{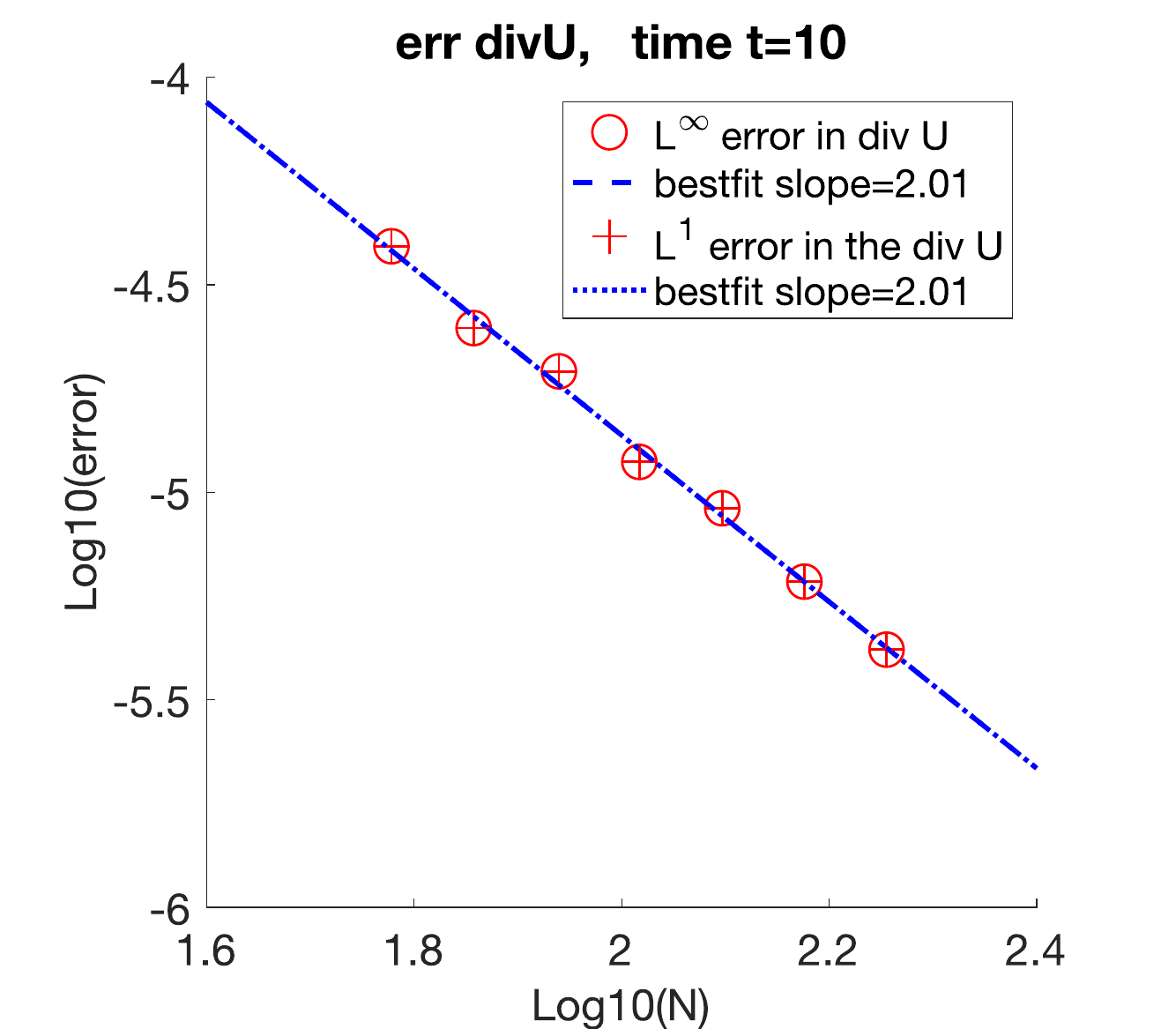}
 \end{minipage}
\caption{ \footnotesize{Test \textsc{Steady Circle}: best fit accuracy order at time $t=10$. The exact solution is \eqref{sol:exa}.} }
	\label{fig:circle_steady}
\end{figure}

\begin{figure}[H]
 \begin{minipage}[c]{0.33\textwidth}
   	\centering
   	\captionsetup{width=0.80\textwidth}
		\includegraphics[width=0.99\textwidth]{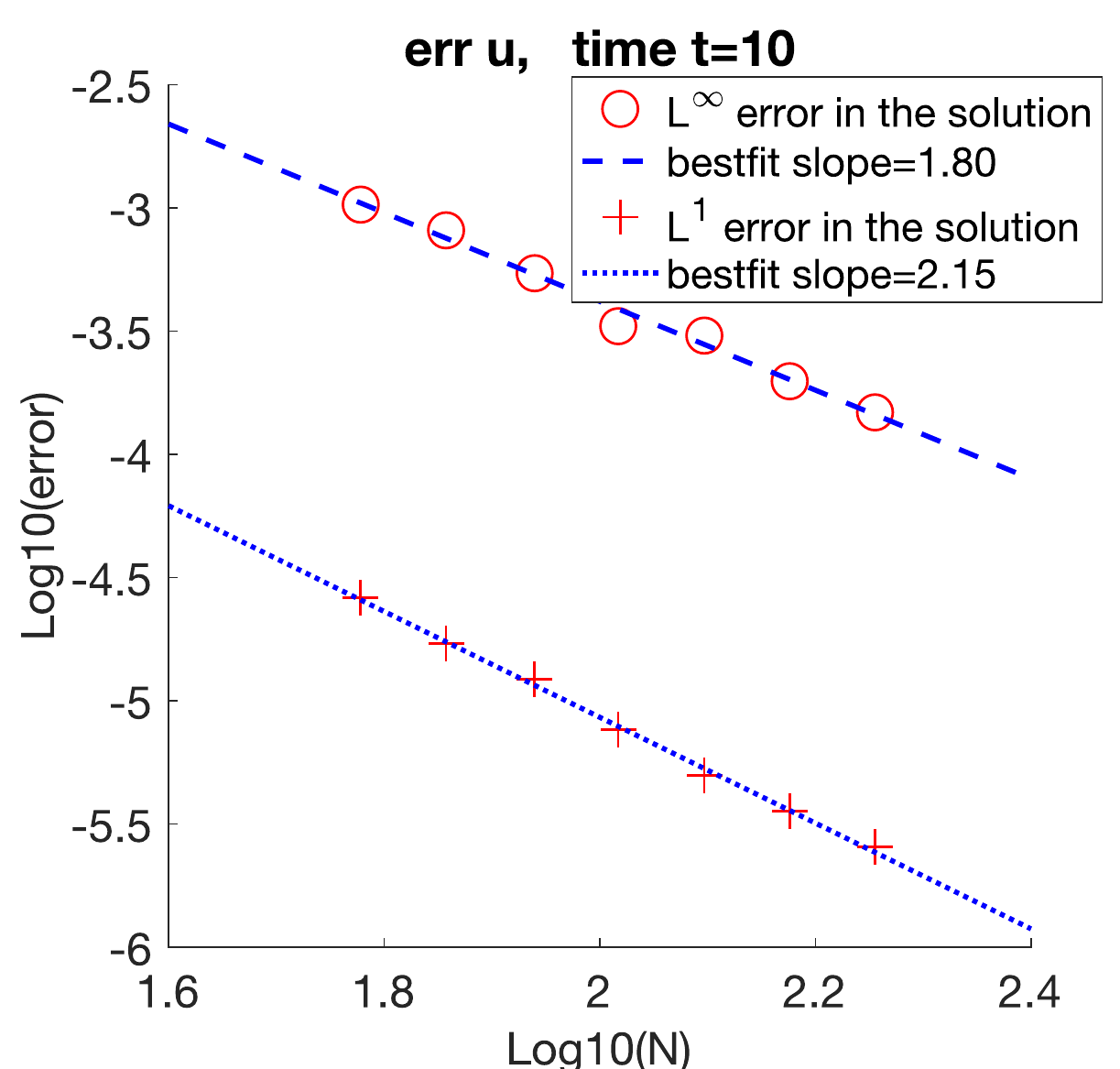}
 \end{minipage}
 \begin{minipage}[c]{0.33\textwidth}
   	\centering
   	\captionsetup{width=0.80\textwidth}
		\includegraphics[width=0.99\textwidth]{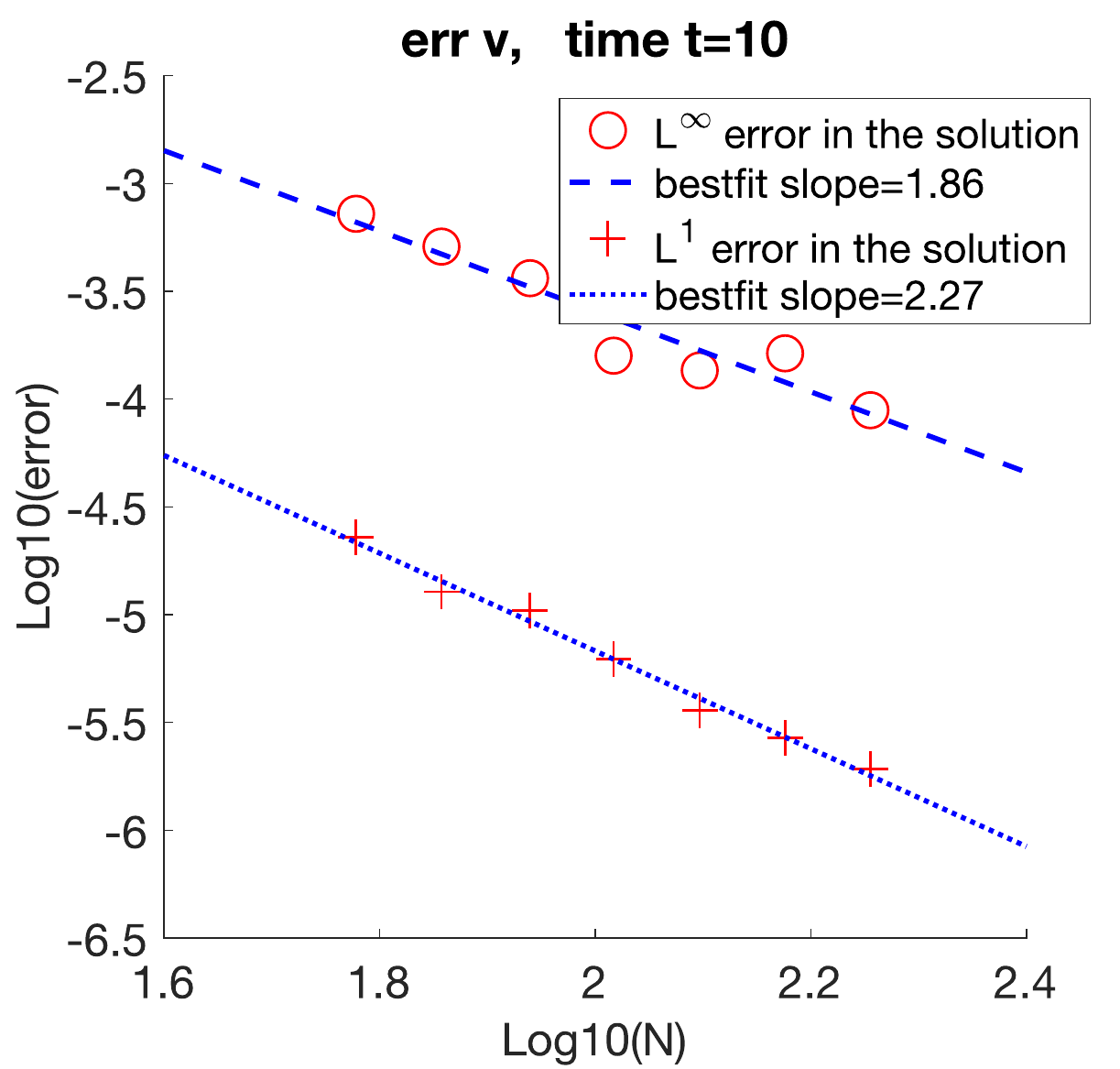}
 \end{minipage}
  \begin{minipage}[c]{0.33\textwidth}
   	\centering
   	\captionsetup{width=0.80\textwidth}
		\includegraphics[width=0.99\textwidth]{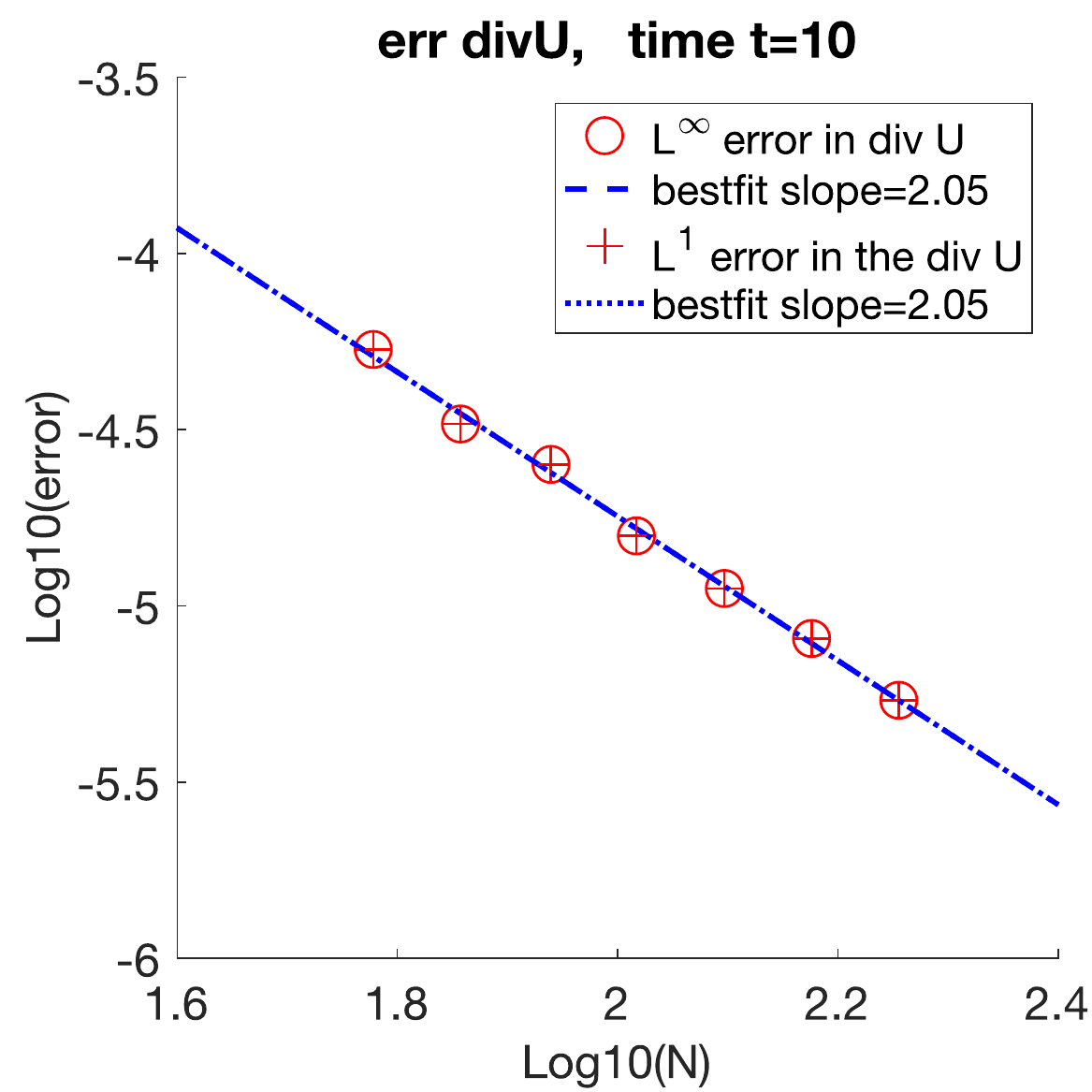}
 \end{minipage}
\caption{ \footnotesize{Test \textsc{Steady Ellipse}: best fit accuracy order at time $t=10$. The exact solution is \eqref{sol:exa}.} }
	\label{fig:ellipse_steady}
\end{figure}

\begin{figure}[H]
 \begin{minipage}[c]{0.33\textwidth}
   	\centering
   	\captionsetup{width=0.80\textwidth}
		\includegraphics[width=0.99\textwidth]{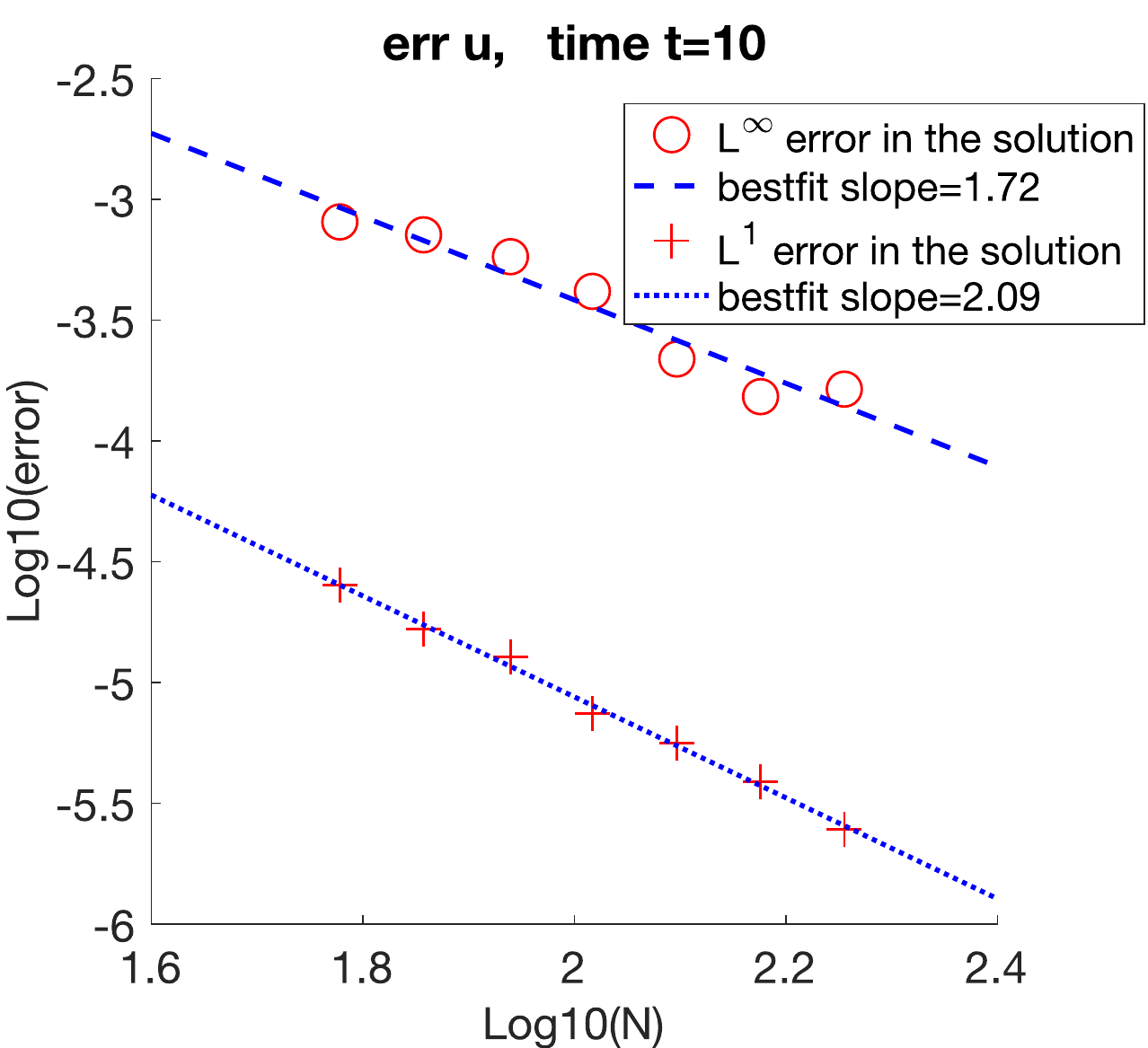}
 \end{minipage}
 \begin{minipage}[c]{0.33\textwidth}
   	\centering
   	\captionsetup{width=0.80\textwidth}
		\includegraphics[width=0.99\textwidth]{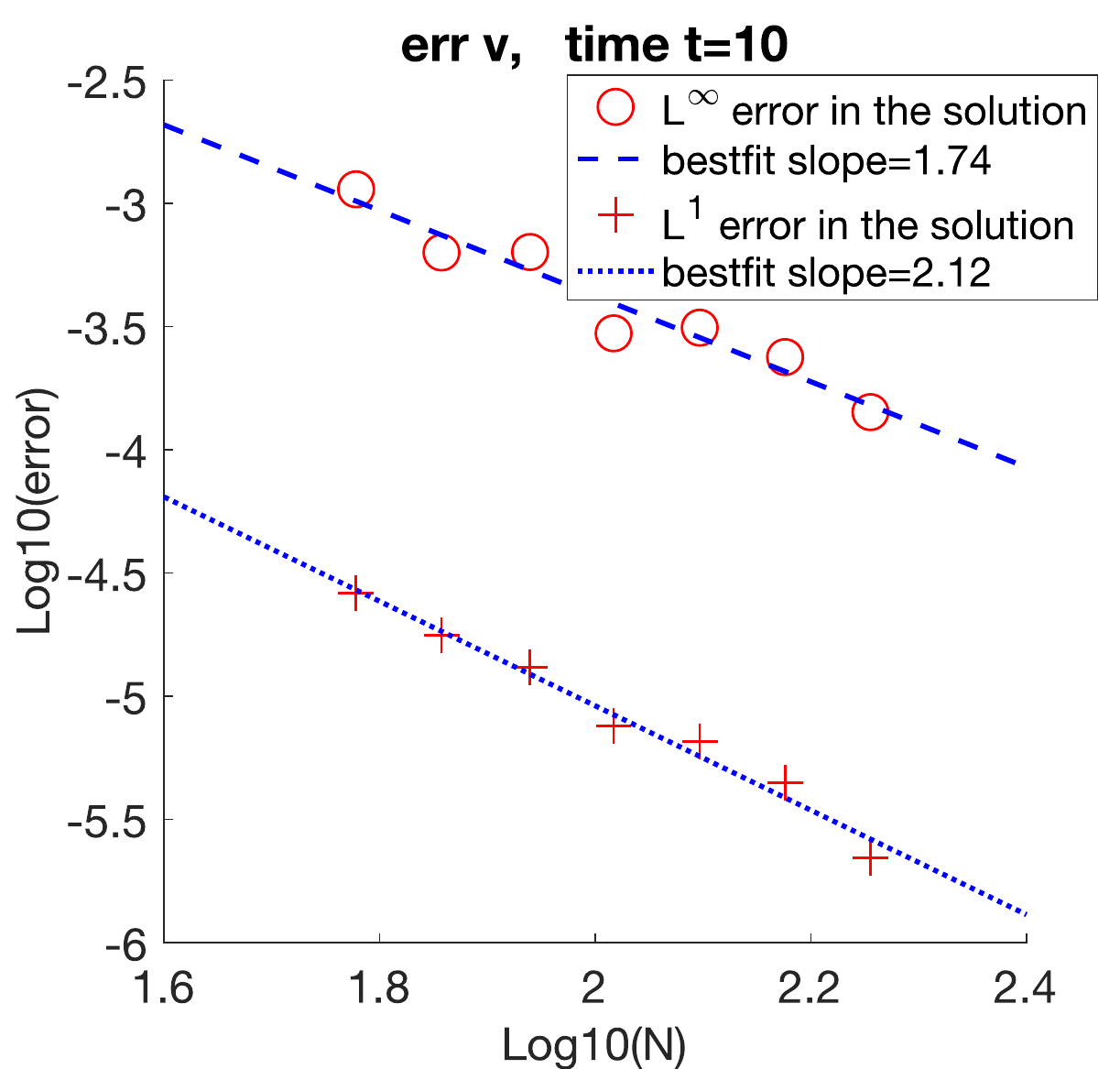}
 \end{minipage}
  \begin{minipage}[c]{0.33\textwidth}
   	\centering
   	\captionsetup{width=0.80\textwidth}
		\includegraphics[width=0.99\textwidth]{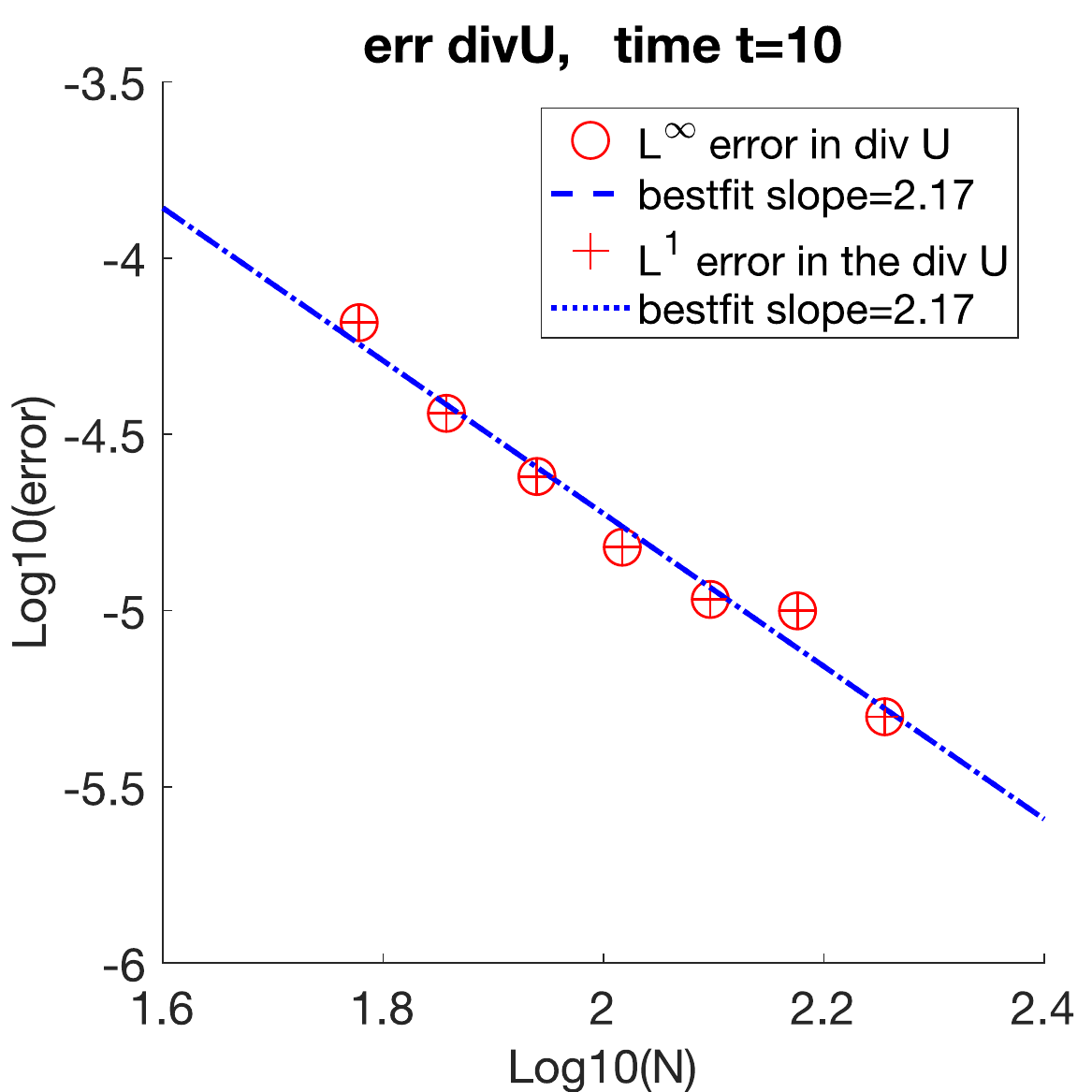}
 \end{minipage}
\caption{ \footnotesize{Test \textsc{Steady Flower}: best fit accuracy order at time $t=10$. The exact solution is \eqref{sol:exa}.} }
	\label{fig:flower_steady}
\end{figure}

\begin{figure}[H]
 \begin{minipage}[c]{0.33\textwidth}
   	\centering
   	\captionsetup{width=0.80\textwidth}
		\includegraphics[width=0.99\textwidth]{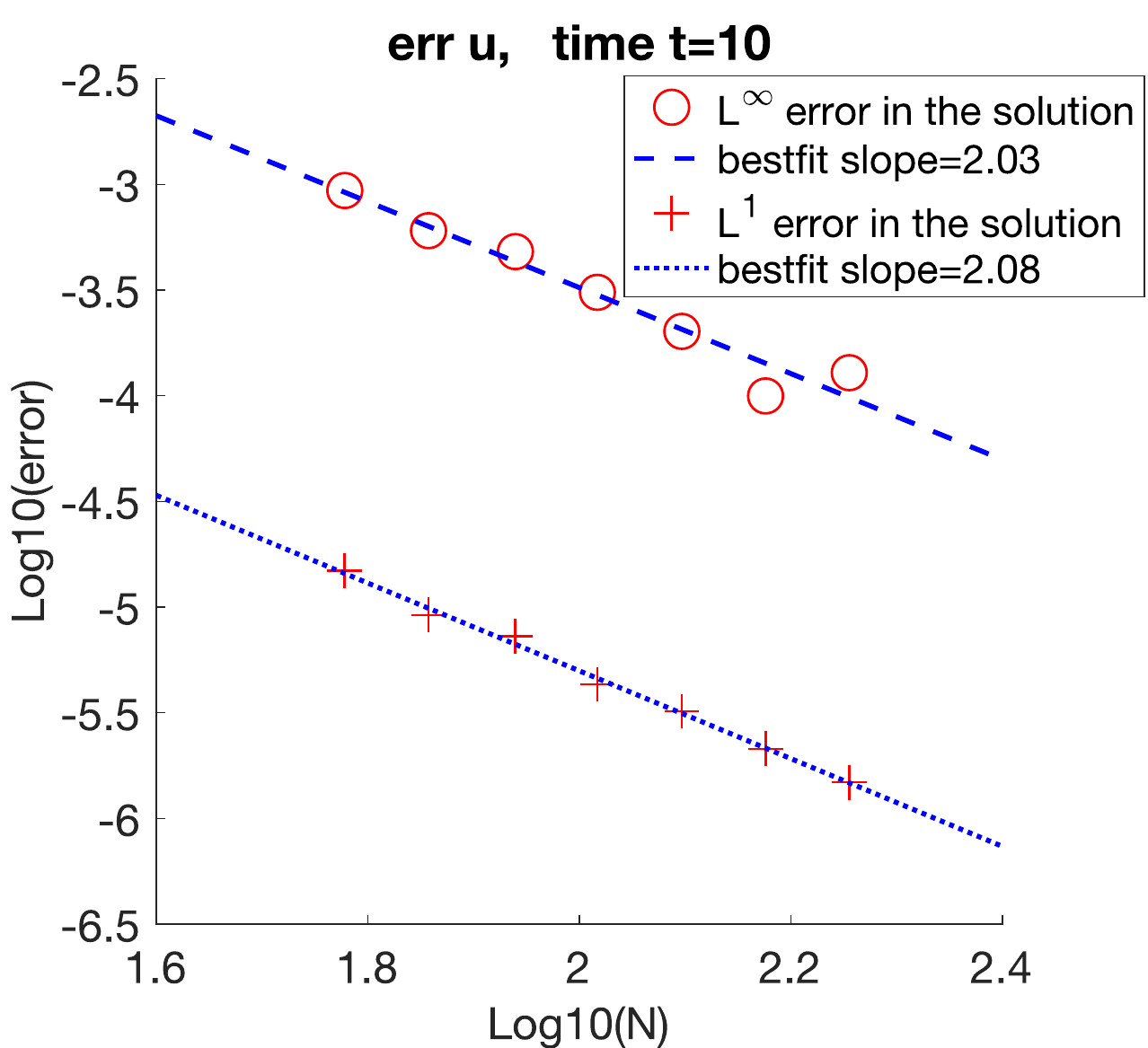}
 \end{minipage}
 \begin{minipage}[c]{0.33\textwidth}
   	\centering
   	\captionsetup{width=0.80\textwidth}
		\includegraphics[width=0.99\textwidth]{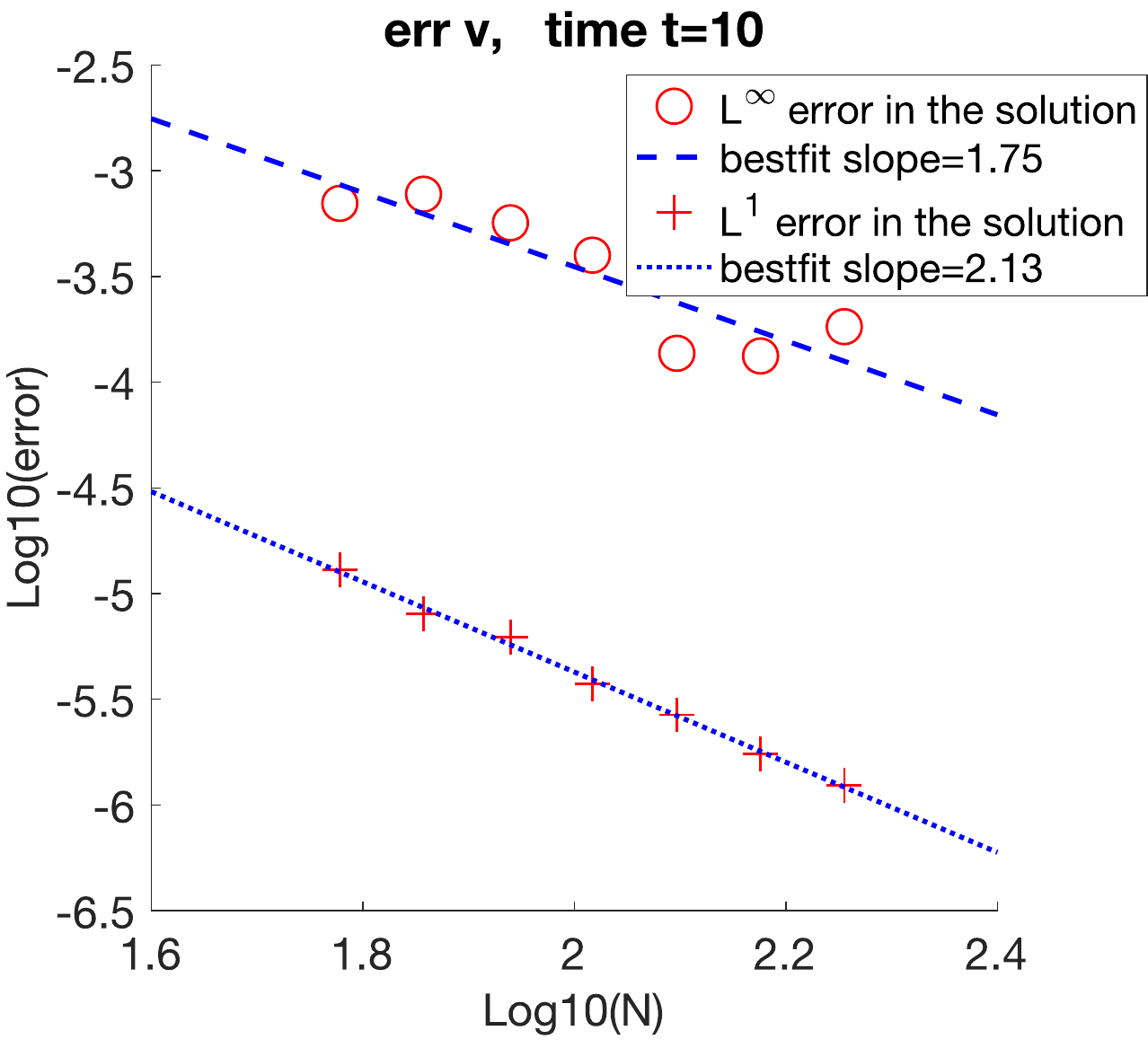}
 \end{minipage}
  \begin{minipage}[c]{0.33\textwidth}
   	\centering
   	\captionsetup{width=0.80\textwidth}
		\includegraphics[width=0.99\textwidth]{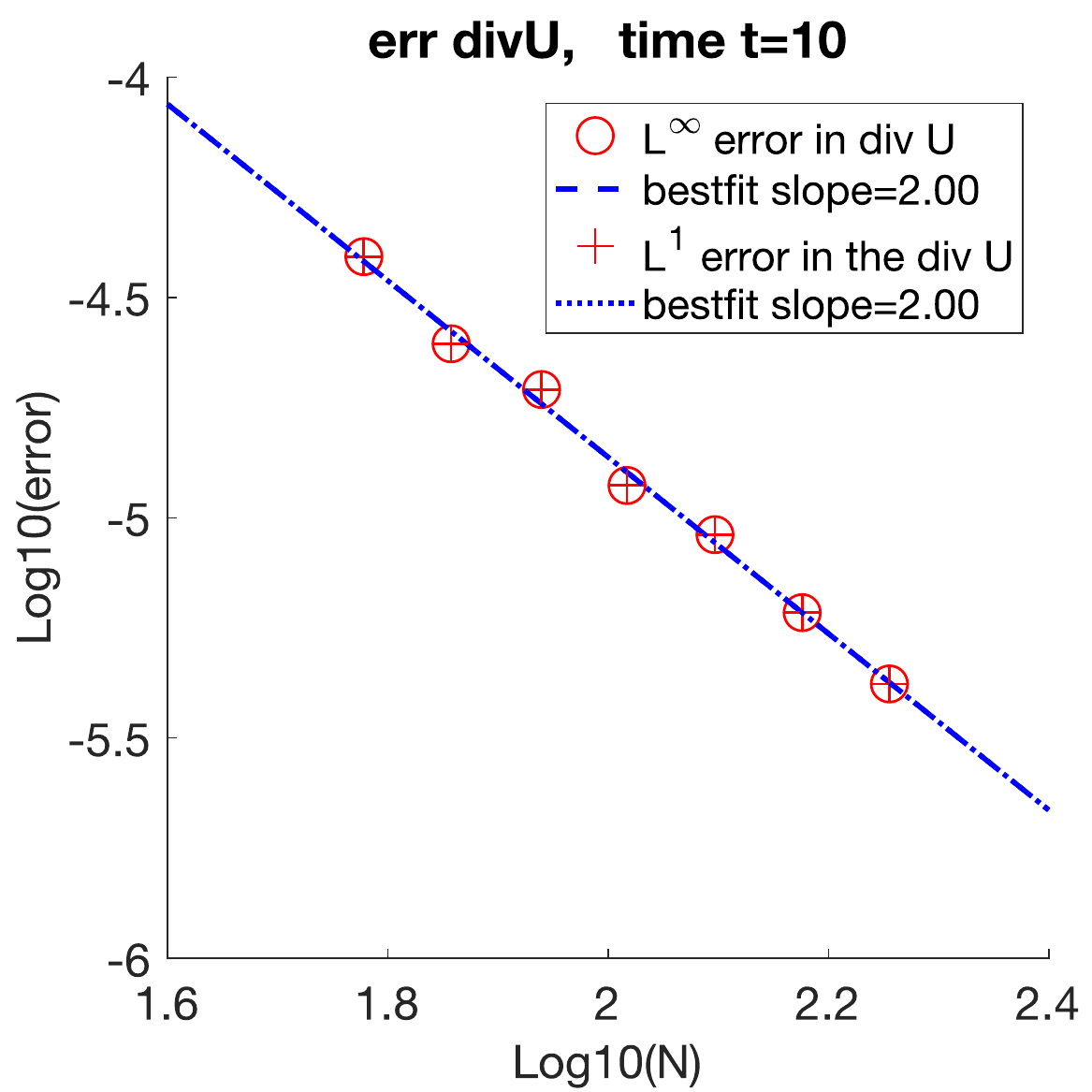}
 \end{minipage}
\caption{ \footnotesize{Test \textsc{Moving Circle}: best fit accuracy order at time $t=10$. The exact solution is \eqref{sol:exa}; the angular velocity is $\omega = 2 \pi / 5$.} }
	\label{fig:circle_moving}
\end{figure}

\begin{figure}[H]
 \begin{minipage}[c]{0.33\textwidth}
   	\centering
   	\captionsetup{width=0.80\textwidth}
		\includegraphics[width=0.99\textwidth]{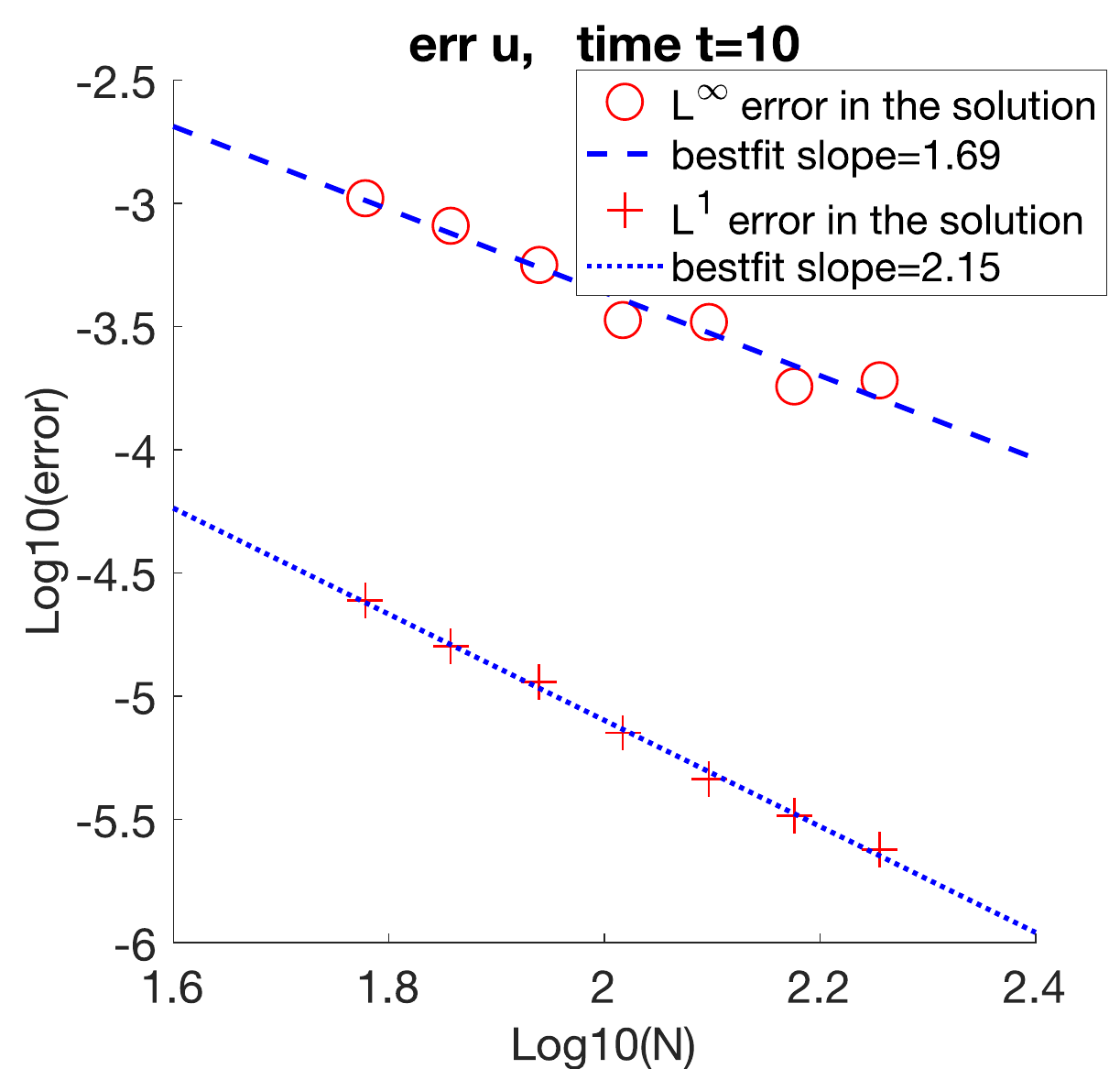}
 \end{minipage}
 \begin{minipage}[c]{0.33\textwidth}
   	\centering
   	\captionsetup{width=0.80\textwidth}
		\includegraphics[width=0.99\textwidth]{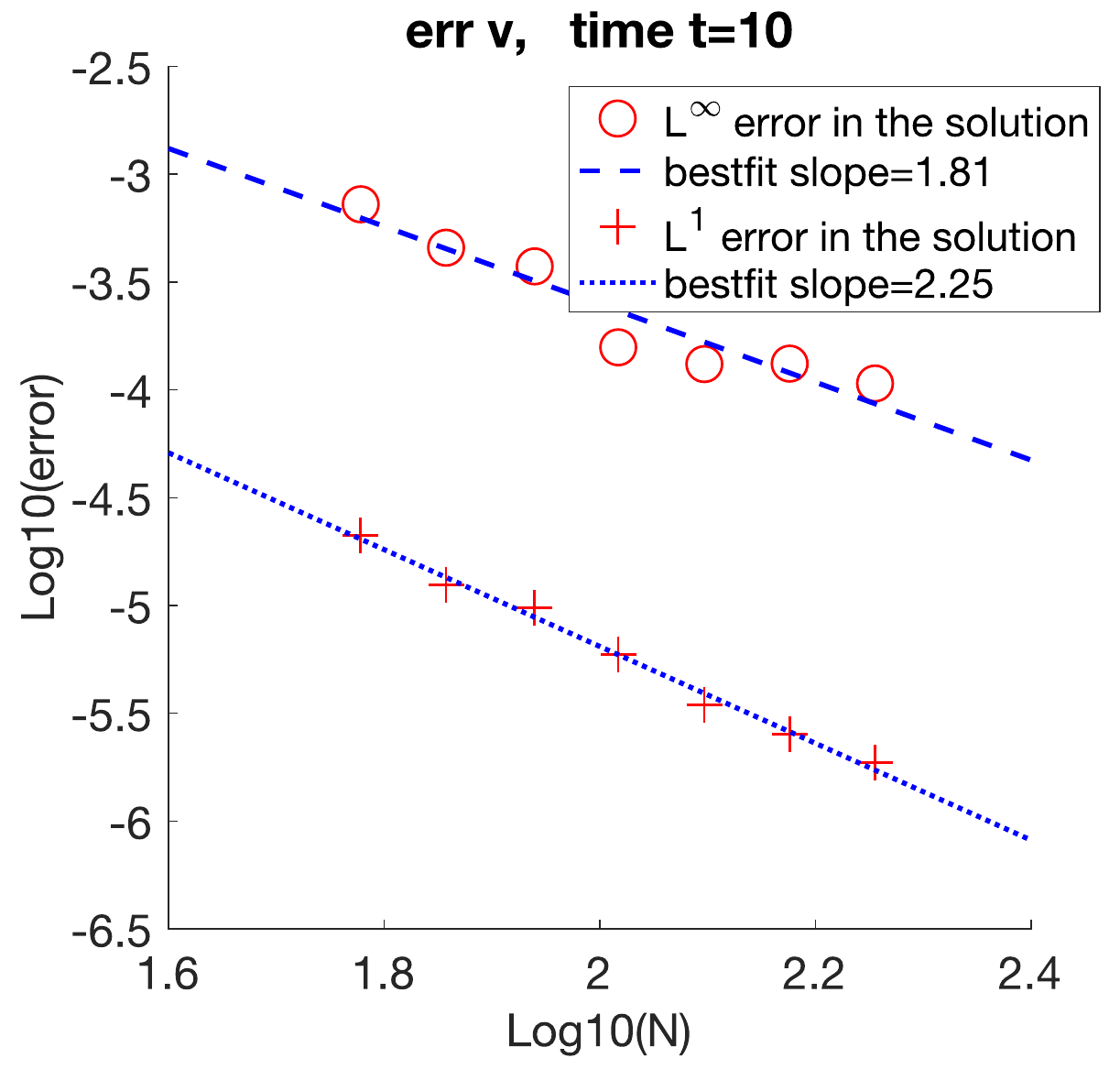}
 \end{minipage}
  \begin{minipage}[c]{0.33\textwidth}
   	\centering
   	\captionsetup{width=0.80\textwidth}
		\includegraphics[width=0.99\textwidth]{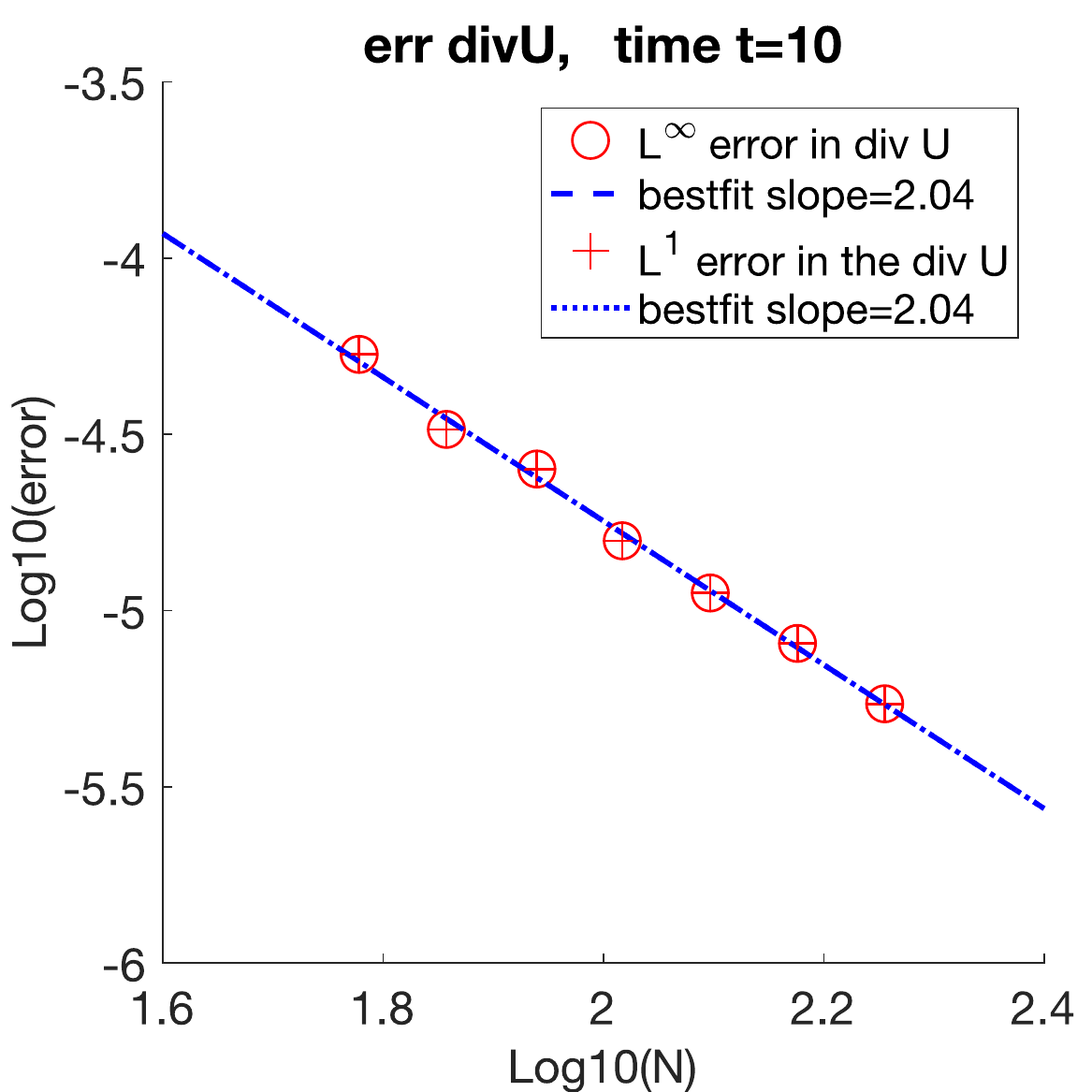}
 \end{minipage}
\caption{ \footnotesize{Test \textsc{Moving Ellipse}: best fit accuracy order at time $t=10$. The exact solution is \eqref{sol:exa}; the angular velocity is $\omega = 2 \pi / 5$.} }
	\label{fig:ellipse_moving}
\end{figure}

\begin{figure}[H]
 \begin{minipage}[c]{0.33\textwidth}
   	\centering
   	\captionsetup{width=0.80\textwidth}
		\includegraphics[width=0.99\textwidth]{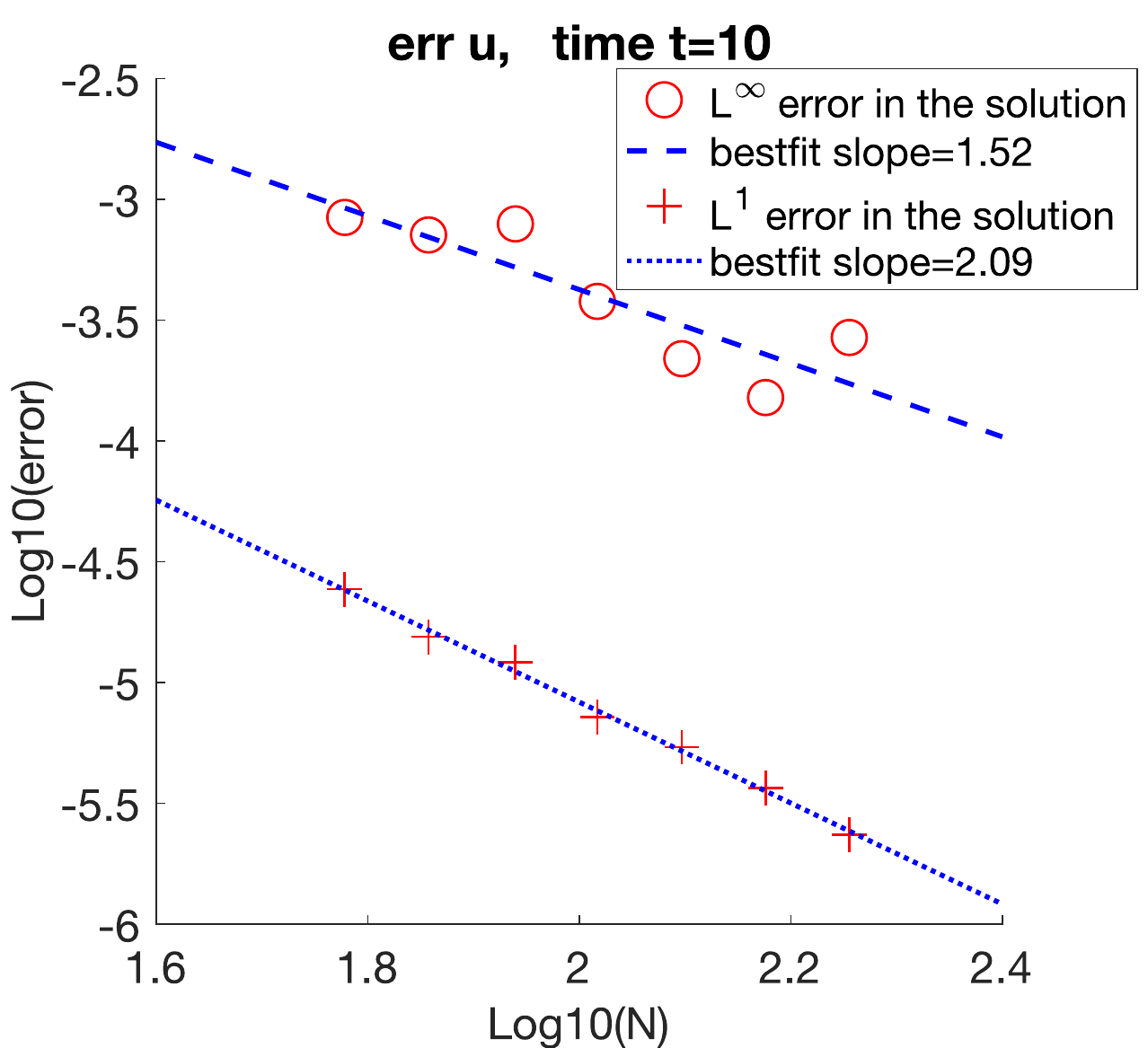}
 \end{minipage}
 \begin{minipage}[c]{0.33\textwidth}
   	\centering
   	\captionsetup{width=0.80\textwidth}
		\includegraphics[width=0.99\textwidth]{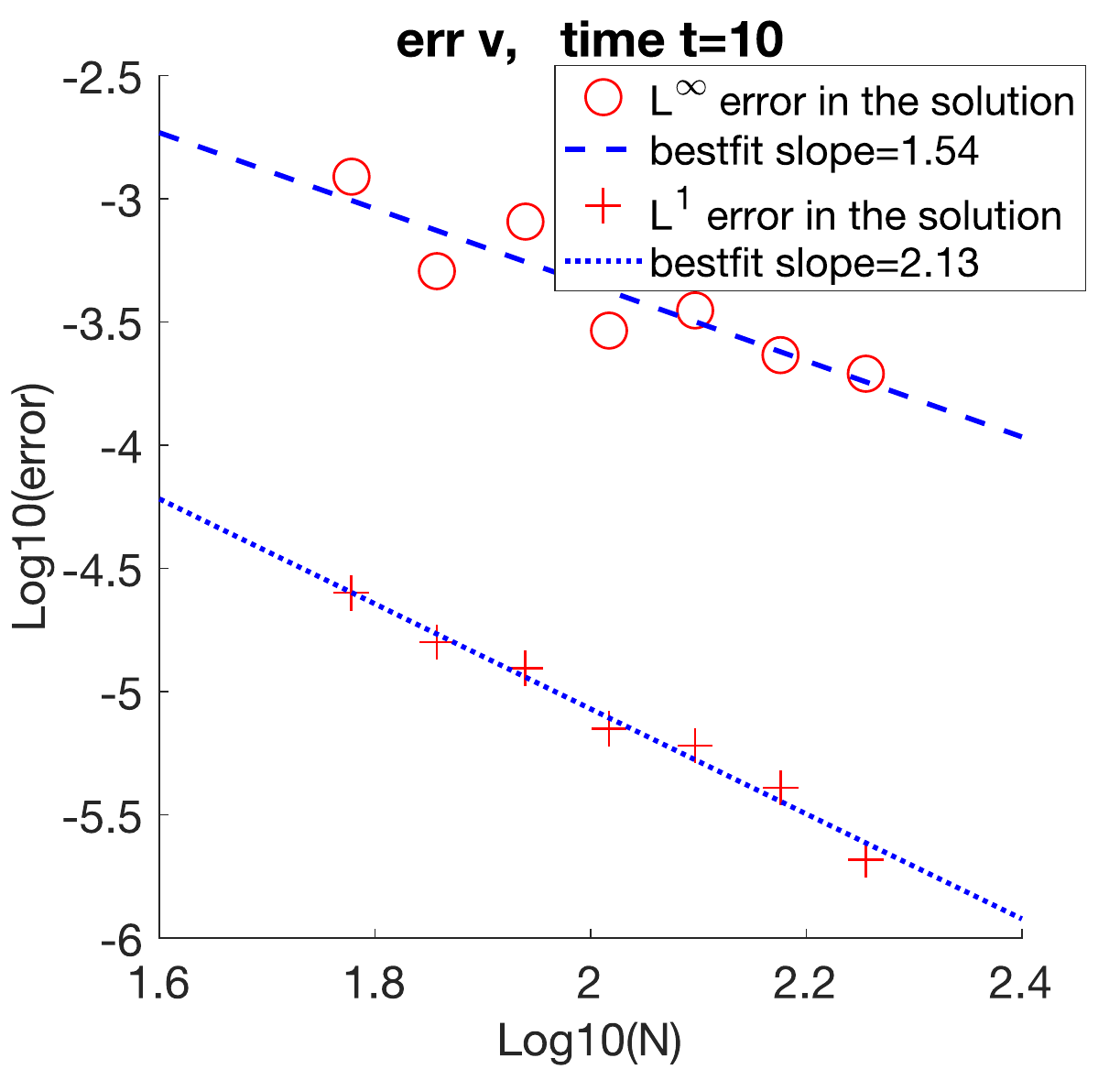}
 \end{minipage}
  \begin{minipage}[c]{0.33\textwidth}
   	\centering
   	\captionsetup{width=0.80\textwidth}
		\includegraphics[width=0.99\textwidth]{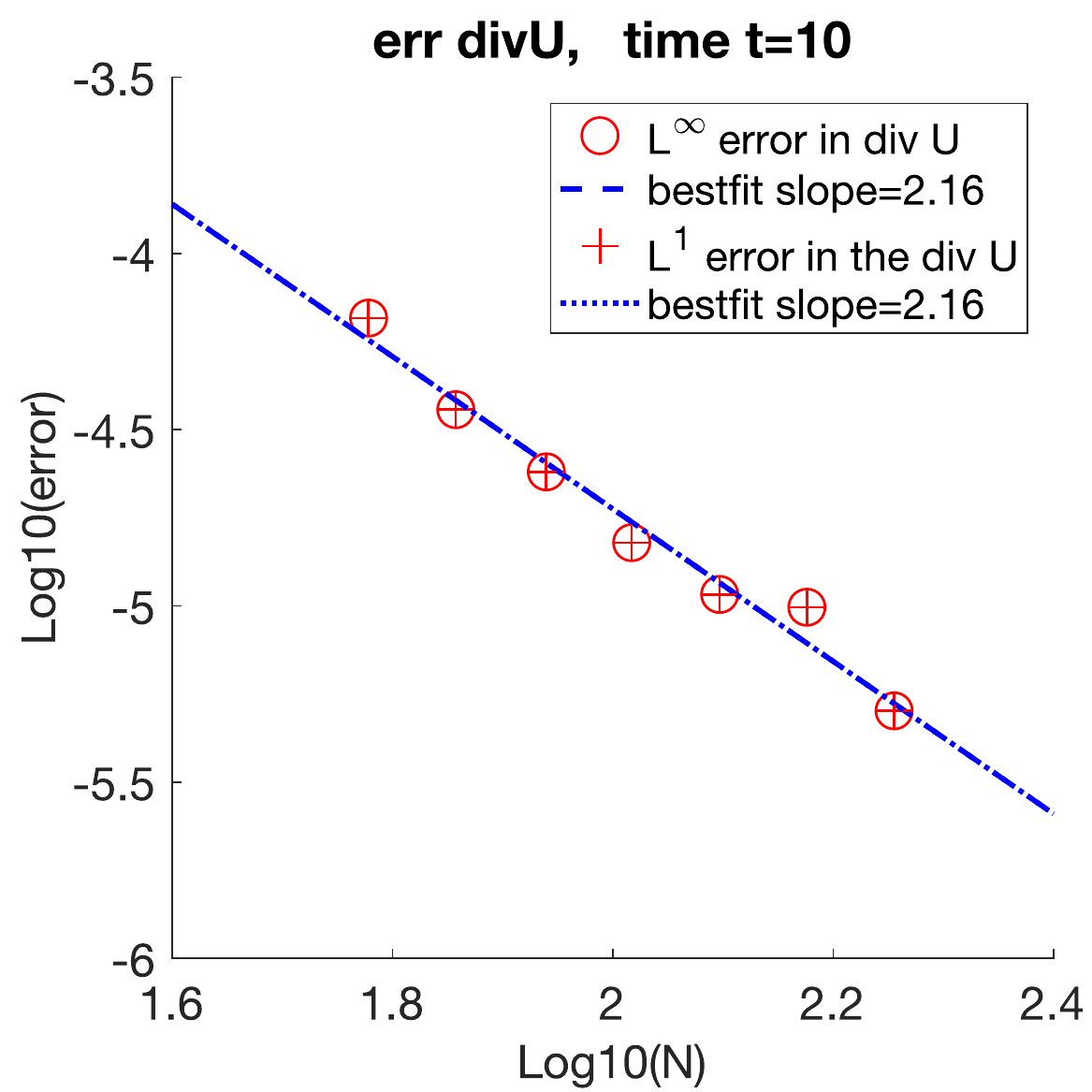}
 \end{minipage}
\caption{ \footnotesize{Test \textsc{Moving Flower}: best fit accuracy order at time $t=10$. The exact solution is \eqref{sol:exa}; the angular velocity is $\omega = 2 \pi / 5$.} }
	\label{fig:flower_moving}
\end{figure}

In Figs.~\ref{fig:circle_steady}, \ref{fig:ellipse_steady} and \ref{fig:flower_steady} we plot the best fit line for the tests \textsc{Steady Circle}, \textsc{Steady Ellipse} and \textsc{Steady Flower}, respectively. The error values are reported in Tables \ref{table:circle_steady}, \ref{table:ellipse_steady} and \ref{table:flower_steady}.
The tests with moving domains, i.e.~\textsc{Moving Circle}, \textsc{Moving Ellipse} and \textsc{Moving Flower}, are plotted in Figs.~\ref{fig:circle_moving}, \ref{fig:ellipse_moving} and \ref{fig:flower_moving}, respectively, with the error values reported in Tables \ref{table:circle_moving}, \ref{table:ellipse_moving} and \ref{table:flower_moving}.

The best fit line shows second order accuracy for all tests in the $L^1$ norm, and for most tests in the $L^\infty$ norm. 
The only test that does not show a clear second order accuracy in the $L^\infty$ norm is the \textsc{Moving Flower} test, although we expect to obtain a second order accuracy asymptotically. 
In fact, we observe (not shown here) that the numerical error is mainly concentrated close to the boundary of the object and we speculate that a more stable second order decay is observed when the ratio $\kappa / N$ is sufficiently small, where $\kappa$ is the curvature of the boundary. Since the curvature of the flower is in some regions higher than the one of the circle, we expect to observe a clear decay for the flower domain with a higher value of $N$.
In general, the norms of the numerical errors fluctuate in all tests around the best fit line, especially in the $L^\infty$ norm. This is mainly due to interpolation errors of ghost values. In fact, from one value of $N$ to the other, the structure of the grid and then of the ghost points around the boundary can significantly change. The effect is more evident for more complex geometries, as the flower-shaped object, and for moving domains. 
The error of $\nabla \cdot \vec{u}$ is computed as $\left\| \vec{g} - \nabla \cdot \vec{u} \right\|$. 
Due to the presence of $\xi$ in the continuity equation (see Eq.~\ref{spaceprobaug}), we do not expect to have a divergence-free solution (or, in this case, a solution that satisfies the equation $\nabla \cdot \vec{u} = \vec{g}$). What we observe is however a second order decay of the error. 
We observe that the $L^1$ and $L^\infty$ norms of the error of $\nabla \cdot \vec{u}$ are the same, since $ \nabla \cdot \vec{u} - \vec{g} = \xi$ is constant in $\Omega$.

 \begin{table}[H]
\captionsetup{width=0.80\textwidth}
\centering      
\begin{tabular}{|| c || c | c || c | c || c | c ||} 
\hline \hline 
No.~of points & $L^1$ error of $u$ & order & $L^1$ error of $v$ & order & $L^1$ error of $\nabla \cdot \vec{u}$ & order \\ 
\hline 
60 $\times$ 60 & 1.55 $\cdot 10^{-5}$ & - & 1.32 $\cdot 10^{-5}$ & - & 3.91 $\cdot 10^{-5}$ & - \\ 
72 $\times$ 72 & 9.33 $\cdot 10^{-6}$ & 2.79 & 8.05 $\cdot 10^{-6}$ & 2.73 & 2.49 $\cdot 10^{-5}$ & 2.49 \\ 
87 $\times$ 87 & 7.69 $\cdot 10^{-6}$ & 1.02 & 6.60 $\cdot 10^{-6}$ & 1.05 & 1.96 $\cdot 10^{-5}$ & 1.27 \\ 
104 $\times$ 104 & 4.38 $\cdot 10^{-6}$ & 3.15 & 3.75 $\cdot 10^{-6}$ & 3.16 & 1.18 $\cdot 10^{-5}$ & 2.82 \\ 
125 $\times$ 125 & 3.38 $\cdot 10^{-6}$ & 1.42 & 2.76 $\cdot 10^{-6}$ & 1.68 & 9.16 $\cdot 10^{-6}$ & 1.40 \\ 
150 $\times$ 150 & 2.20 $\cdot 10^{-6}$ & 2.34 & 1.75 $\cdot 10^{-6}$ & 2.51 & 6.09 $\cdot 10^{-6}$ & 2.23 \\ 
180 $\times$ 180 & 1.52 $\cdot 10^{-6}$ & 2.03 & 1.23 $\cdot 10^{-6}$ & 1.90 & 4.18 $\cdot 10^{-6}$ & 2.07 \\ 
\hline \hline 
No.~of points & $L^\infty$ error of $u$ & order & $L^\infty$ error of $v$ & order & $L^\infty$ error of $\nabla \cdot \vec{u}$ & order \\ 
\hline 
60 $\times$ 60 & 6.52 $\cdot 10^{-4}$ & - & 7.08 $\cdot 10^{-4}$ & - & 3.91 $\cdot 10^{-5}$ & - \\ 
72 $\times$ 72 & 6.26 $\cdot 10^{-4}$ & 0.22 & 7.96 $\cdot 10^{-4}$ & -0.64 & 2.49 $\cdot 10^{-5}$ & 2.49 \\ 
87 $\times$ 87 & 4.79 $\cdot 10^{-4}$ & 1.42 & 5.71 $\cdot 10^{-4}$ & 1.75 & 1.96 $\cdot 10^{-5}$ & 1.27 \\ 
104 $\times$ 104 & 3.06 $\cdot 10^{-4}$ & 2.50 & 3.95 $\cdot 10^{-4}$ & 2.07 & 1.18 $\cdot 10^{-5}$ & 2.82 \\ 
125 $\times$ 125 & 1.82 $\cdot 10^{-4}$ & 2.83 & 1.37 $\cdot 10^{-4}$ & 5.74 & 9.16 $\cdot 10^{-6}$ & 1.40 \\ 
150 $\times$ 150 & 1.01 $\cdot 10^{-4}$ & 3.24 & 1.34 $\cdot 10^{-4}$ & 0.14 & 6.09 $\cdot 10^{-6}$ & 2.23 \\ 
180 $\times$ 180 & 1.15 $\cdot 10^{-4}$ & -0.72 & 1.41 $\cdot 10^{-4}$ & -0.29 & 4.18 $\cdot 10^{-6}$ & 2.07 \\ 
\hline \hline 
\end{tabular} 
\caption{ \footnotesize{Test \textsc{Steady Circle}: errors and accuracy orders in the $L^1$ (top) and $L^\infty$ (bottom) norms, for $u$ (left), $v$ (middle) and $\nabla \cdot \vec{u}$ (right). The exact solution is \eqref{sol:exa}.}}
    \label{table:circle_steady}  
 \end{table}

  \begin{table}[H]
\captionsetup{width=0.80\textwidth}
\centering      
\begin{tabular}{|| c || c | c || c | c || c | c ||} 
\hline \hline 
No.~of points & $L^1$ error of $u$ & order & $L^1$ error of $v$ & order & $L^1$ error of $\nabla \cdot \vec{u}$ & order \\ 
\hline 
60 $\times$ 60 & 2.61 $\cdot 10^{-5}$ & - & 2.29 $\cdot 10^{-5}$ & - & 5.34 $\cdot 10^{-5}$ & - \\ 
72 $\times$ 72 & 1.70 $\cdot 10^{-5}$ & 2.35 & 1.28 $\cdot 10^{-5}$ & 3.21 & 3.28 $\cdot 10^{-5}$ & 2.67 \\ 
87 $\times$ 87 & 1.22 $\cdot 10^{-5}$ & 1.75 & 1.05 $\cdot 10^{-5}$ & 1.04 & 2.52 $\cdot 10^{-5}$ & 1.39 \\ 
104 $\times$ 104 & 7.63 $\cdot 10^{-6}$ & 2.65 & 6.24 $\cdot 10^{-6}$ & 2.91 & 1.58 $\cdot 10^{-5}$ & 2.62 \\ 
125 $\times$ 125 & 4.98 $\cdot 10^{-6}$ & 2.31 & 3.59 $\cdot 10^{-6}$ & 3.01 & 1.12 $\cdot 10^{-5}$ & 1.87 \\ 
150 $\times$ 150 & 3.57 $\cdot 10^{-6}$ & 1.83 & 2.69 $\cdot 10^{-6}$ & 1.59 & 8.08 $\cdot 10^{-6}$ & 1.79 \\ 
180 $\times$ 180 & 2.55 $\cdot 10^{-6}$ & 1.85 & 1.92 $\cdot 10^{-6}$ & 1.83 & 5.40 $\cdot 10^{-6}$ & 2.21 \\ 
\hline \hline 
No.~of points & $L^\infty$ error of $u$ & order & $L^\infty$ error of $v$ & order & $L^\infty$ error of $\nabla \cdot \vec{u}$ & order \\ 
\hline 
60 $\times$ 60 & 1.03 $\cdot 10^{-3}$ & - & 7.23 $\cdot 10^{-4}$ & - & 5.34 $\cdot 10^{-5}$ & - \\ 
72 $\times$ 72 & 8.12 $\cdot 10^{-4}$ & 1.32 & 5.09 $\cdot 10^{-4}$ & 1.92 & 3.28 $\cdot 10^{-5}$ & 2.67 \\ 
87 $\times$ 87 & 5.44 $\cdot 10^{-4}$ & 2.12 & 3.63 $\cdot 10^{-4}$ & 1.79 & 2.52 $\cdot 10^{-5}$ & 1.39 \\ 
104 $\times$ 104 & 3.31 $\cdot 10^{-4}$ & 2.77 & 1.59 $\cdot 10^{-4}$ & 4.63 & 1.58 $\cdot 10^{-5}$ & 2.62 \\ 
125 $\times$ 125 & 3.04 $\cdot 10^{-4}$ & 0.48 & 1.36 $\cdot 10^{-4}$ & 0.86 & 1.12 $\cdot 10^{-5}$ & 1.87 \\ 
150 $\times$ 150 & 1.98 $\cdot 10^{-4}$ & 2.34 & 1.63 $\cdot 10^{-4}$ & -1.00 & 8.08 $\cdot 10^{-6}$ & 1.79 \\ 
180 $\times$ 180 & 1.48 $\cdot 10^{-4}$ & 1.59 & 8.87 $\cdot 10^{-5}$ & 3.33 & 5.40 $\cdot 10^{-6}$ & 2.21 \\ 
\hline \hline 
\end{tabular} 
\caption{ \footnotesize{Test \textsc{Steady Ellipse}: errors and accuracy orders in the $L^1$ (top) and $L^\infty$ (bottom) norms, for $u$ (left), $v$ (middle) and $\nabla \cdot \vec{u}$ (right). The exact solution is \eqref{sol:exa}.}}
    \label{table:ellipse_steady}  
 \end{table}

  \begin{table}[H]
\captionsetup{width=0.80\textwidth}
\centering      
\begin{tabular}{|| c || c | c || c | c || c | c ||} 
\hline \hline 
No.~of points & $L^1$ error of $u$ & order & $L^1$ error of $v$ & order & $L^1$ error of $\nabla \cdot \vec{u}$ & order \\ 
\hline 
60 $\times$ 60 & 2.53 $\cdot 10^{-5}$ & - & 2.62 $\cdot 10^{-5}$ & - & 6.56 $\cdot 10^{-5}$ & - \\ 
72 $\times$ 72 & 1.66 $\cdot 10^{-5}$ & 2.31 & 1.77 $\cdot 10^{-5}$ & 2.14 & 3.63 $\cdot 10^{-5}$ & 3.25 \\ 
87 $\times$ 87 & 1.28 $\cdot 10^{-5}$ & 1.40 & 1.31 $\cdot 10^{-5}$ & 1.59 & 2.40 $\cdot 10^{-5}$ & 2.18 \\ 
104 $\times$ 104 & 7.47 $\cdot 10^{-6}$ & 3.00 & 7.58 $\cdot 10^{-6}$ & 3.08 & 1.51 $\cdot 10^{-5}$ & 2.59 \\ 
125 $\times$ 125 & 5.60 $\cdot 10^{-6}$ & 1.56 & 6.57 $\cdot 10^{-6}$ & 0.77 & 1.08 $\cdot 10^{-5}$ & 1.86 \\ 
150 $\times$ 150 & 3.89 $\cdot 10^{-6}$ & 2.00 & 4.45 $\cdot 10^{-6}$ & 2.14 & 9.98 $\cdot 10^{-6}$ & 0.41 \\ 
180 $\times$ 180 & 2.47 $\cdot 10^{-6}$ & 2.48 & 2.21 $\cdot 10^{-6}$ & 3.83 & 5.00 $\cdot 10^{-6}$ & 3.79 \\ 
\hline \hline 
No.~of points & $L^\infty$ error of $u$ & order & $L^\infty$ error of $v$ & order & $L^\infty$ error of $\nabla \cdot \vec{u}$ & order \\ 
\hline 
60 $\times$ 60 & 8.07 $\cdot 10^{-4}$ & - & 1.14 $\cdot 10^{-3}$ & - & 6.56 $\cdot 10^{-5}$ & - \\ 
72 $\times$ 72 & 7.13 $\cdot 10^{-4}$ & 0.68 & 6.32 $\cdot 10^{-4}$ & 3.24 & 3.63 $\cdot 10^{-5}$ & 3.25 \\ 
87 $\times$ 87 & 5.79 $\cdot 10^{-4}$ & 1.10 & 6.35 $\cdot 10^{-4}$ & -0.02 & 2.40 $\cdot 10^{-5}$ & 2.18 \\ 
104 $\times$ 104 & 4.16 $\cdot 10^{-4}$ & 1.85 & 2.97 $\cdot 10^{-4}$ & 4.26 & 1.51 $\cdot 10^{-5}$ & 2.59 \\ 
125 $\times$ 125 & 2.18 $\cdot 10^{-4}$ & 3.50 & 3.13 $\cdot 10^{-4}$ & -0.29 & 1.08 $\cdot 10^{-5}$ & 1.86 \\ 
150 $\times$ 150 & 1.53 $\cdot 10^{-4}$ & 1.97 & 2.38 $\cdot 10^{-4}$ & 1.51 & 9.98 $\cdot 10^{-6}$ & 0.41 \\ 
180 $\times$ 180 & 1.64 $\cdot 10^{-4}$ & -0.39 & 1.42 $\cdot 10^{-4}$ & 2.82 & 5.00 $\cdot 10^{-6}$ & 3.79 \\ 
\hline \hline 
\end{tabular}
\caption{ \footnotesize{Test \textsc{Steady Flower}: errors and accuracy orders in the $L^1$ (top) and $L^\infty$ (bottom) norms, for $u$ (left), $v$ (middle) and $\nabla \cdot \vec{u}$ (right). The exact solution is \eqref{sol:exa}.}}
    \label{table:flower_steady}  
 \end{table}

  \begin{table}[H]
\captionsetup{width=0.80\textwidth}
\centering      
\begin{tabular}{|| c || c | c || c | c || c | c ||} 
\hline \hline 
No.~of points & $L^1$ error of $u$ & order & $L^1$ error of $v$ & order & $L^1$ error of $\nabla \cdot \vec{u}$ & order \\ 
\hline 
60 $\times$ 60 & 1.49 $\cdot 10^{-5}$ & - & 1.29 $\cdot 10^{-5}$ & - & 3.91 $\cdot 10^{-5}$ & - \\ 
72 $\times$ 72 & 9.17 $\cdot 10^{-6}$ & 2.65 & 8.05 $\cdot 10^{-6}$ & 2.60 & 2.48 $\cdot 10^{-5}$ & 2.50 \\ 
87 $\times$ 87 & 7.28 $\cdot 10^{-6}$ & 1.22 & 6.23 $\cdot 10^{-6}$ & 1.35 & 1.95 $\cdot 10^{-5}$ & 1.26 \\ 
104 $\times$ 104 & 4.31 $\cdot 10^{-6}$ & 2.93 & 3.75 $\cdot 10^{-6}$ & 2.85 & 1.18 $\cdot 10^{-5}$ & 2.81 \\ 
125 $\times$ 125 & 3.22 $\cdot 10^{-6}$ & 1.59 & 2.67 $\cdot 10^{-6}$ & 1.85 & 9.16 $\cdot 10^{-6}$ & 1.40 \\ 
150 $\times$ 150 & 2.14 $\cdot 10^{-6}$ & 2.25 & 1.74 $\cdot 10^{-6}$ & 2.34 & 6.09 $\cdot 10^{-6}$ & 2.24 \\ 
180 $\times$ 180 & 1.48 $\cdot 10^{-6}$ & 1.99 & 1.24 $\cdot 10^{-6}$ & 1.86 & 4.20 $\cdot 10^{-6}$ & 2.04 \\ 
\hline \hline 
No.~of points & $L^\infty$ error of $u$ & order & $L^\infty$ error of $v$ & order & $L^\infty$ error of $\nabla \cdot \vec{u}$ & order \\ 
\hline 
60 $\times$ 60 & 9.34 $\cdot 10^{-4}$ & - & 7.01 $\cdot 10^{-4}$ & - & 3.91 $\cdot 10^{-5}$ & - \\ 
72 $\times$ 72 & 6.03 $\cdot 10^{-4}$ & 2.40 & 7.75 $\cdot 10^{-4}$ & -0.55 & 2.48 $\cdot 10^{-5}$ & 2.50 \\ 
87 $\times$ 87 & 4.80 $\cdot 10^{-4}$ & 1.20 & 5.67 $\cdot 10^{-4}$ & 1.65 & 1.95 $\cdot 10^{-5}$ & 1.26 \\ 
104 $\times$ 104 & 3.08 $\cdot 10^{-4}$ & 2.48 & 3.99 $\cdot 10^{-4}$ & 1.98 & 1.18 $\cdot 10^{-5}$ & 2.81 \\ 
125 $\times$ 125 & 2.01 $\cdot 10^{-4}$ & 2.31 & 1.37 $\cdot 10^{-4}$ & 5.81 & 9.16 $\cdot 10^{-6}$ & 1.40 \\ 
150 $\times$ 150 & 9.97 $\cdot 10^{-5}$ & 3.86 & 1.33 $\cdot 10^{-4}$ & 0.16 & 6.09 $\cdot 10^{-6}$ & 2.24 \\ 
180 $\times$ 180 & 1.28 $\cdot 10^{-4}$ & -1.39 & 1.83 $\cdot 10^{-4}$ & -1.75 & 4.20 $\cdot 10^{-6}$ & 2.04 \\ 
\hline \hline 
\end{tabular}
\caption{ \footnotesize{Test \textsc{Moving Circle}: errors and accuracy orders in the $L^1$ (top) and $L^\infty$ (bottom) norms, for $u$ (left), $v$ (middle) and $\nabla \cdot \vec{u}$ (right). The exact solution is \eqref{sol:exa}; the angular velocity is $\omega = 2 \pi / 5$.}}
    \label{table:circle_moving}  
 \end{table}
 
  \begin{table}[H]
\captionsetup{width=0.80\textwidth}
\centering      
\begin{tabular}{|| c || c | c || c | c || c | c ||} 
\hline \hline 
No.~of points & $L^1$ error of $u$ & order & $L^1$ error of $v$ & order & $L^1$ error of $\nabla \cdot \vec{u}$ & order \\ 
\hline 
60 $\times$ 60 & 2.44 $\cdot 10^{-5}$ & - & 2.12 $\cdot 10^{-5}$ & - & 5.33 $\cdot 10^{-5}$ & - \\ 
72 $\times$ 72 & 1.59 $\cdot 10^{-5}$ & 2.35 & 1.25 $\cdot 10^{-5}$ & 2.90 & 3.27 $\cdot 10^{-5}$ & 2.69 \\ 
87 $\times$ 87 & 1.14 $\cdot 10^{-5}$ & 1.77 & 9.78 $\cdot 10^{-6}$ & 1.29 & 2.52 $\cdot 10^{-5}$ & 1.37 \\ 
104 $\times$ 104 & 7.11 $\cdot 10^{-6}$ & 2.65 & 5.93 $\cdot 10^{-6}$ & 2.81 & 1.58 $\cdot 10^{-5}$ & 2.63 \\ 
125 $\times$ 125 & 4.61 $\cdot 10^{-6}$ & 2.35 & 3.45 $\cdot 10^{-6}$ & 2.95 & 1.12 $\cdot 10^{-5}$ & 1.85 \\ 
150 $\times$ 150 & 3.28 $\cdot 10^{-6}$ & 1.87 & 2.54 $\cdot 10^{-6}$ & 1.68 & 8.07 $\cdot 10^{-6}$ & 1.81 \\ 
180 $\times$ 180 & 2.39 $\cdot 10^{-6}$ & 1.73 & 1.87 $\cdot 10^{-6}$ & 1.68 & 5.44 $\cdot 10^{-6}$ & 2.16 \\ 
\hline \hline 
No.~of points & $L^\infty$ error of $u$ & order & $L^\infty$ error of $v$ & order & $L^\infty$ error of $\nabla \cdot \vec{u}$ & order \\ 
\hline 
60 $\times$ 60 & 1.05 $\cdot 10^{-3}$ & - & 7.25 $\cdot 10^{-4}$ & - & 5.33 $\cdot 10^{-5}$ & - \\ 
72 $\times$ 72 & 8.11 $\cdot 10^{-4}$ & 1.41 & 4.56 $\cdot 10^{-4}$ & 2.54 & 3.27 $\cdot 10^{-5}$ & 2.69 \\ 
87 $\times$ 87 & 5.63 $\cdot 10^{-4}$ & 1.93 & 3.74 $\cdot 10^{-4}$ & 1.05 & 2.52 $\cdot 10^{-5}$ & 1.37 \\ 
104 $\times$ 104 & 3.36 $\cdot 10^{-4}$ & 2.89 & 1.57 $\cdot 10^{-4}$ & 4.84 & 1.58 $\cdot 10^{-5}$ & 2.63 \\ 
125 $\times$ 125 & 3.30 $\cdot 10^{-4}$ & 0.10 & 1.32 $\cdot 10^{-4}$ & 0.97 & 1.12 $\cdot 10^{-5}$ & 1.85 \\ 
150 $\times$ 150 & 1.81 $\cdot 10^{-4}$ & 3.30 & 1.32 $\cdot 10^{-4}$ & -0.03 & 8.07 $\cdot 10^{-6}$ & 1.81 \\ 
180 $\times$ 180 & 1.91 $\cdot 10^{-4}$ & -0.32 & 1.07 $\cdot 10^{-4}$ & 1.16 & 5.44 $\cdot 10^{-6}$ & 2.16 \\ 
\hline \hline 
\end{tabular}
\caption{ \footnotesize{Test \textsc{Moving Ellipse}: errors and accuracy orders in the $L^1$ (top) and $L^\infty$ (bottom) norms, for $u$ (left), $v$ (middle) and $\nabla \cdot \vec{u}$ (right). The exact solution is \eqref{sol:exa}; the angular velocity is $\omega = 2 \pi / 5$.}}
    \label{table:ellipse_moving}  
 \end{table}

  \begin{table}[H]
\captionsetup{width=0.80\textwidth}
\centering      
\begin{tabular}{|| c || c | c || c | c || c | c ||} 
\hline \hline 
No.~of points & $L^1$ error of $u$ & order & $L^1$ error of $v$ & order & $L^1$ error of $\nabla \cdot \vec{u}$ & order \\ 
\hline 
60 $\times$ 60 & 2.44 $\cdot 10^{-5}$ & - & 2.52 $\cdot 10^{-5}$ & - & 6.55 $\cdot 10^{-5}$ & - \\ 
72 $\times$ 72 & 1.55 $\cdot 10^{-5}$ & 2.49 & 1.59 $\cdot 10^{-5}$ & 2.53 & 3.60 $\cdot 10^{-5}$ & 3.28 \\ 
87 $\times$ 87 & 1.22 $\cdot 10^{-5}$ & 1.28 & 1.24 $\cdot 10^{-5}$ & 1.30 & 2.40 $\cdot 10^{-5}$ & 2.15 \\ 
104 $\times$ 104 & 7.19 $\cdot 10^{-6}$ & 2.94 & 7.06 $\cdot 10^{-6}$ & 3.16 & 1.51 $\cdot 10^{-5}$ & 2.59 \\ 
125 $\times$ 125 & 5.41 $\cdot 10^{-6}$ & 1.55 & 6.04 $\cdot 10^{-6}$ & 0.84 & 1.08 $\cdot 10^{-5}$ & 1.84 \\ 
150 $\times$ 150 & 3.65 $\cdot 10^{-6}$ & 2.15 & 4.07 $\cdot 10^{-6}$ & 2.16 & 9.91 $\cdot 10^{-6}$ & 0.46 \\ 
180 $\times$ 180 & 2.35 $\cdot 10^{-6}$ & 2.43 & 2.08 $\cdot 10^{-6}$ & 3.68 & 5.04 $\cdot 10^{-6}$ & 3.71 \\ 
\hline \hline 
No.~of points & $L^\infty$ error of $u$ & order & $L^\infty$ error of $v$ & order & $L^\infty$ error of $\nabla \cdot \vec{u}$ & order \\ 
\hline 
60 $\times$ 60 & 8.42 $\cdot 10^{-4}$ & - & 1.23 $\cdot 10^{-3}$ & - & 6.55 $\cdot 10^{-5}$ & - \\ 
72 $\times$ 72 & 7.12 $\cdot 10^{-4}$ & 0.92 & 5.09 $\cdot 10^{-4}$ & 4.84 & 3.60 $\cdot 10^{-5}$ & 3.28 \\ 
87 $\times$ 87 & 7.92 $\cdot 10^{-4}$ & -0.56 & 8.08 $\cdot 10^{-4}$ & -2.45 & 2.40 $\cdot 10^{-5}$ & 2.15 \\ 
104 $\times$ 104 & 3.78 $\cdot 10^{-4}$ & 4.14 & 2.92 $\cdot 10^{-4}$ & 5.70 & 1.51 $\cdot 10^{-5}$ & 2.59 \\ 
125 $\times$ 125 & 2.19 $\cdot 10^{-4}$ & 2.97 & 3.52 $\cdot 10^{-4}$ & -1.03 & 1.08 $\cdot 10^{-5}$ & 1.84 \\ 
150 $\times$ 150 & 1.51 $\cdot 10^{-4}$ & 2.03 & 2.32 $\cdot 10^{-4}$ & 2.29 & 9.91 $\cdot 10^{-6}$ & 0.46 \\ 
180 $\times$ 180 & 2.68 $\cdot 10^{-4}$ & -3.14 & 1.95 $\cdot 10^{-4}$ & 0.95 & 5.04 $\cdot 10^{-6}$ & 3.71 \\ 
\hline \hline 
\end{tabular}
\caption{ \footnotesize{Test \textsc{Moving Flower}: errors and accuracy orders in the $L^1$ (top) and $L^\infty$ (bottom) norms, for $u$ (left), $v$ (middle) and $\nabla \cdot \vec{u}$ (right). The exact solution is \eqref{sol:exa}; the angular velocity is $\omega = 2 \pi / 5$.}}
    \label{table:flower_moving}  
 \end{table}

\subsection{Driven cavity} 
In this section we perform the driven cavity test around a steady flower-shaped object. No-slip boundary conditions $\vec{u}=(0,0)$ are prescribed on the boundary of the object and on the walls of the computational domain $(-1,1)^2$, except for the top wall, that models a lid moving horizontally and at a constant velocity, i.e.~$\vec{u}=(1,0)$. External forces, such as gravity, are neglected ($\vec{f}=0$).

\begin{figure}[H]
 \begin{minipage}[c]{0.49\textwidth}
   	\centering
   	\captionsetup{width=0.80\textwidth}
		\includegraphics[width=0.99\textwidth]{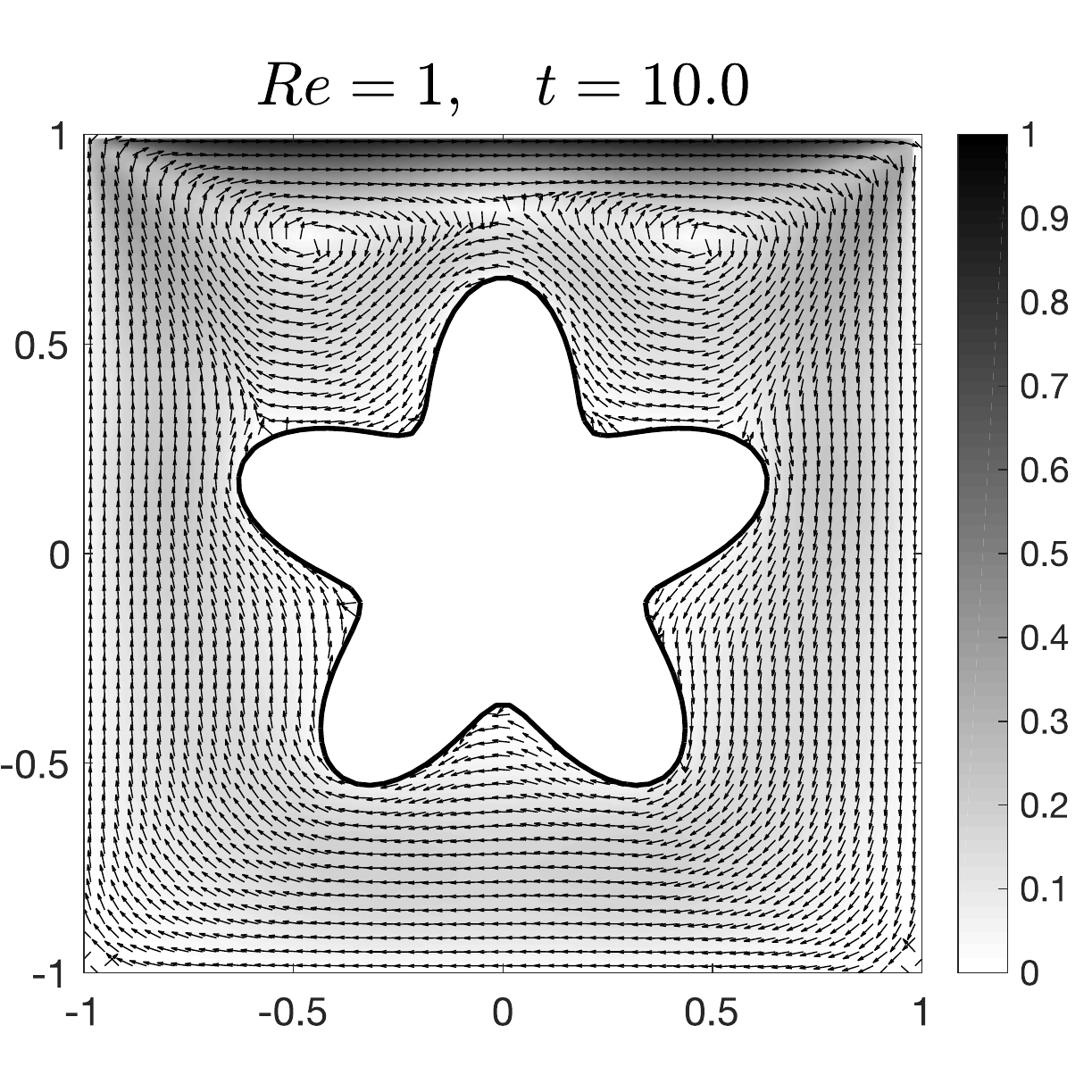}
 \end{minipage}
 \begin{minipage}[c]{0.49\textwidth}
   	\centering
   	\captionsetup{width=0.80\textwidth}
		\includegraphics[width=0.99\textwidth]{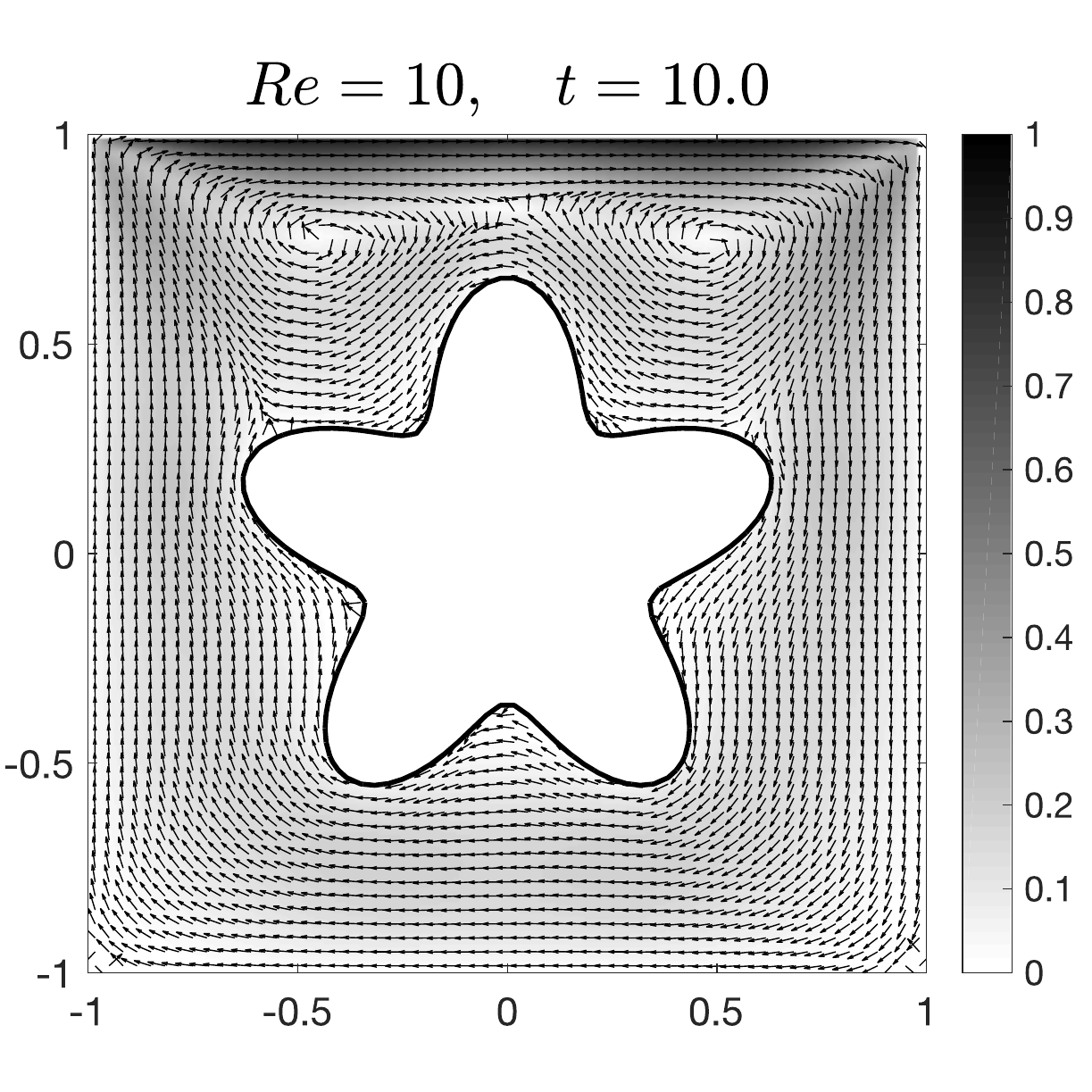}
 \end{minipage}
 \begin{minipage}[c]{0.49\textwidth}
   	\centering
   	\captionsetup{width=0.80\textwidth}
		\includegraphics[width=0.99\textwidth]{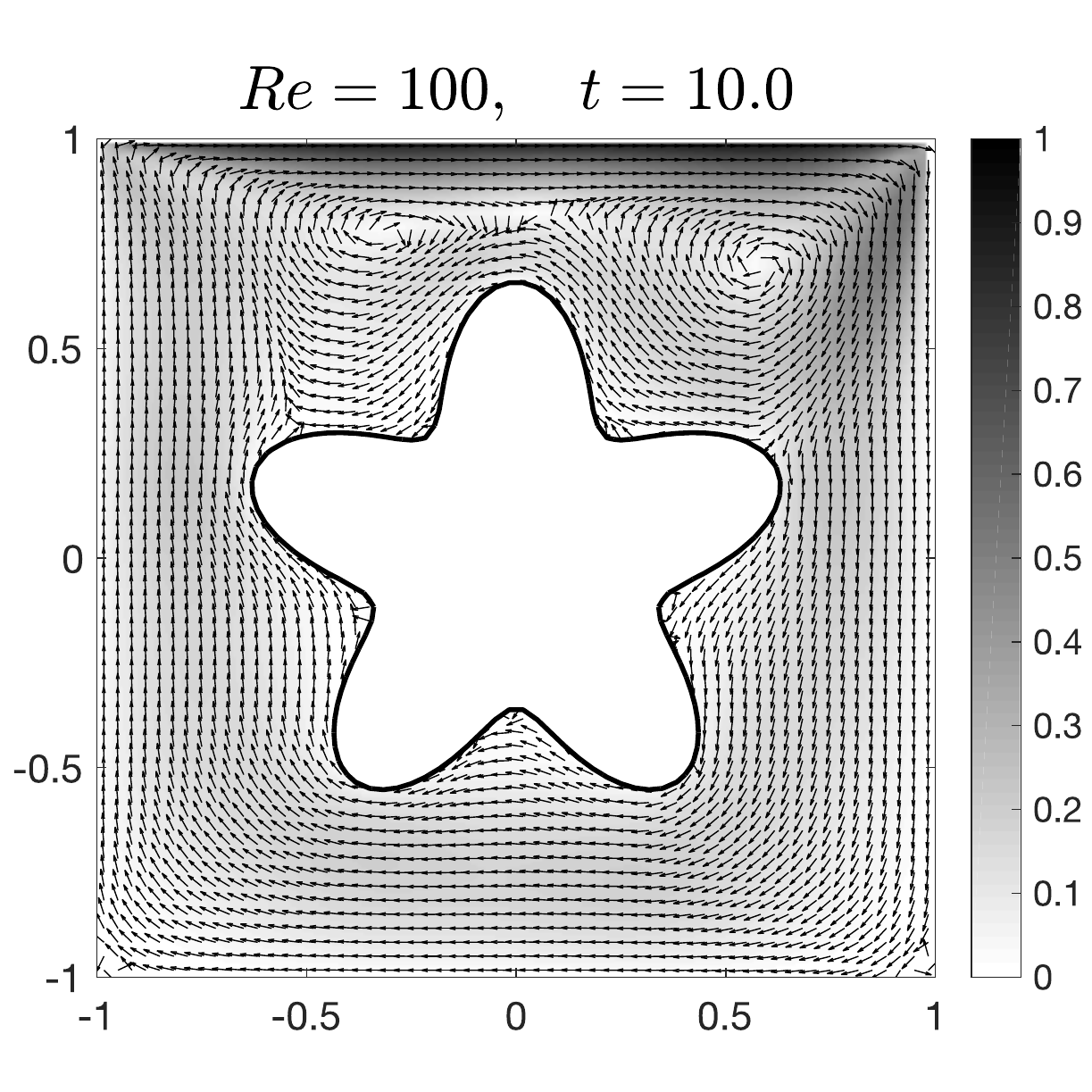}
 \end{minipage}
 \begin{minipage}[c]{0.49\textwidth}
   	\centering
   	\captionsetup{width=0.80\textwidth}
		\includegraphics[width=0.99\textwidth]{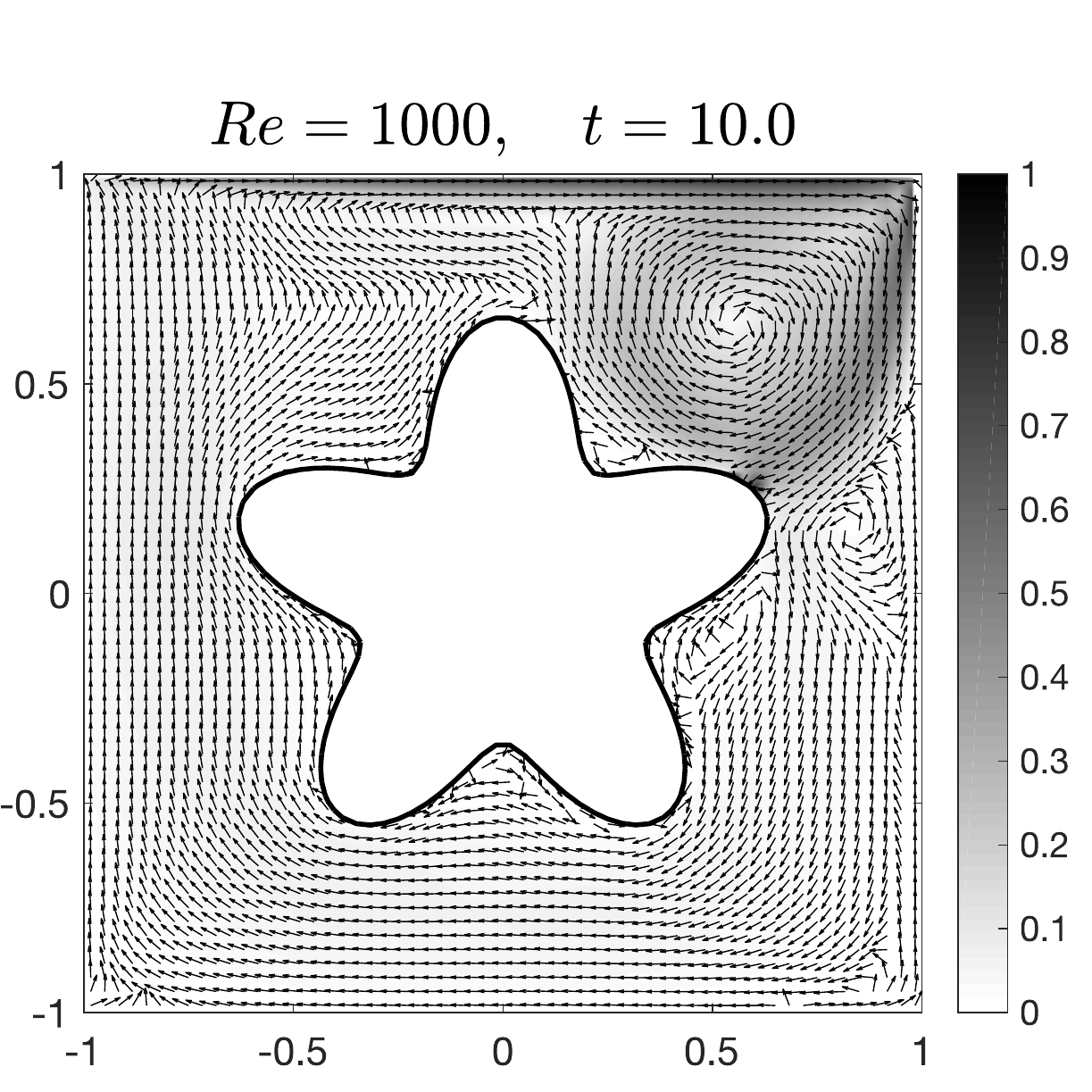}
 \end{minipage}
\caption{ \footnotesize{Velocity field for the driven cavity test with $N=60$. The velocity vectors are scaled in order to have constant length. Their magnitude can be inferred from the color map.} }
	\label{fig:flower_drivencavity}
\end{figure}

In Fig.~\ref{fig:flower_drivencavity} we plot the velocity field at time $t=10$, $N=60$ and with four different values of the Reynolds number: $\text{Re}=1,10,100$ and $1000$. The velocity vectors are scaled in order to have the same length for graphical purposes. Their magnitude can be inferred from the color map.

The presence of the object prevents the formation of anticlockwise vortices in the bottom corners, which are instead observed in classical driven cavity tests (i.e.~without objects immersed in the fluid). Moreover, the central large clockwise vortex observed in classical tests is here split into two vortices for low Reynolds numbers $\text{Re}=1$ and $10$, mainly driven by the reentrant boundaries of the flower-shaped object. The $x$-coordinates of the centres of the vortices are about $x=-0.5$ and $x=0.5$. For $\text{Re}=100$, the vortex on the top-left is more eccentric and slightly moved to the right. For $\text{Re}=1000$ the vortex on the top-left is not observed, while another vortex is forming in the middle-right region.

\subsection{Rotating flower-shaped object}\label{sect:movingFlower}
\begin{figure}[H]
\begin{minipage}[c]{0.33\textwidth}
   	\centering
		\includegraphics[width=0.99\textwidth]{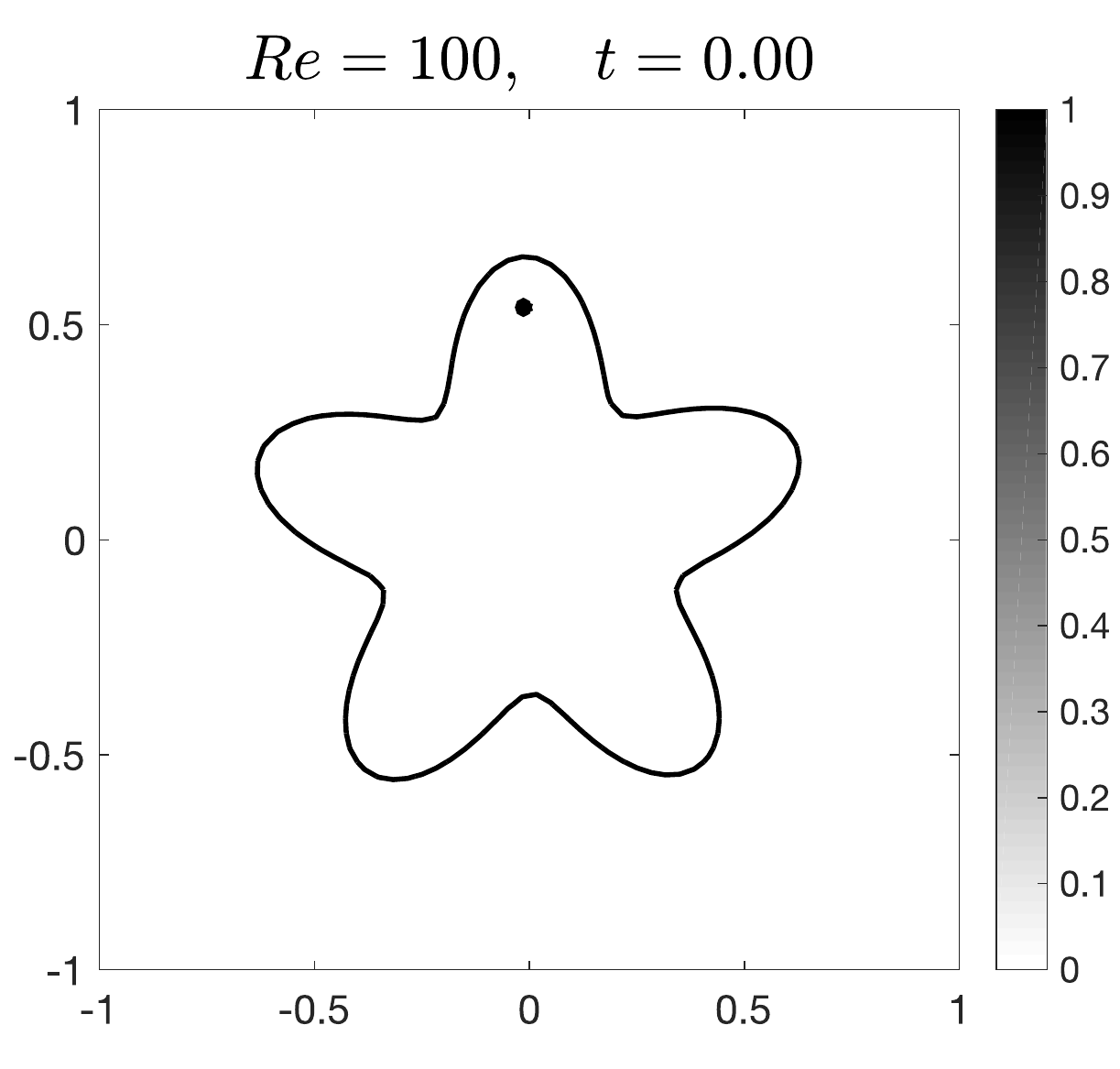}
 \end{minipage}
 \begin{minipage}[c]{0.33\textwidth}
   	\centering
		\includegraphics[width=0.99\textwidth]{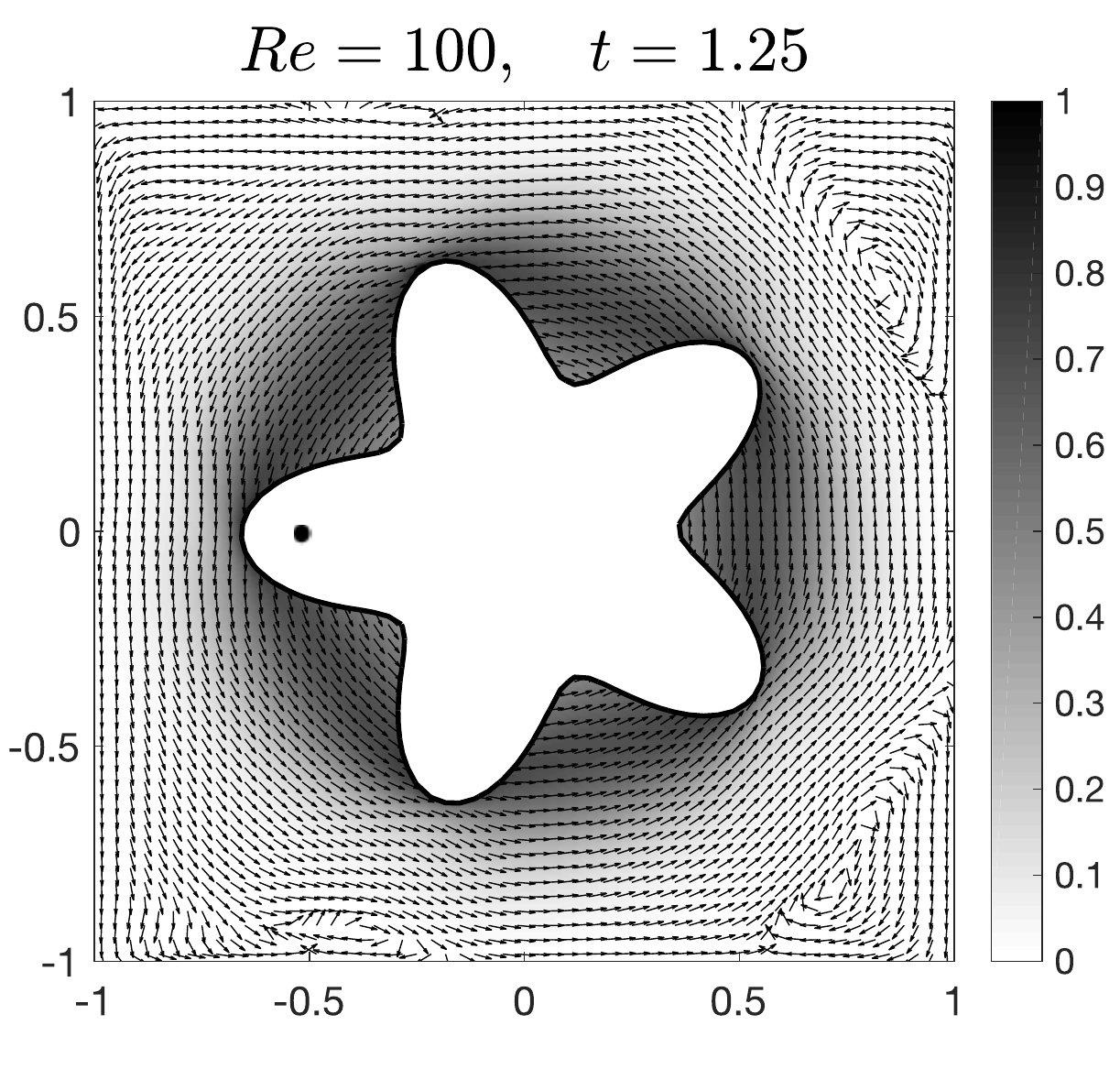}
 \end{minipage}
 \begin{minipage}[c]{0.33\textwidth}
   	\centering
		\includegraphics[width=0.99\textwidth]{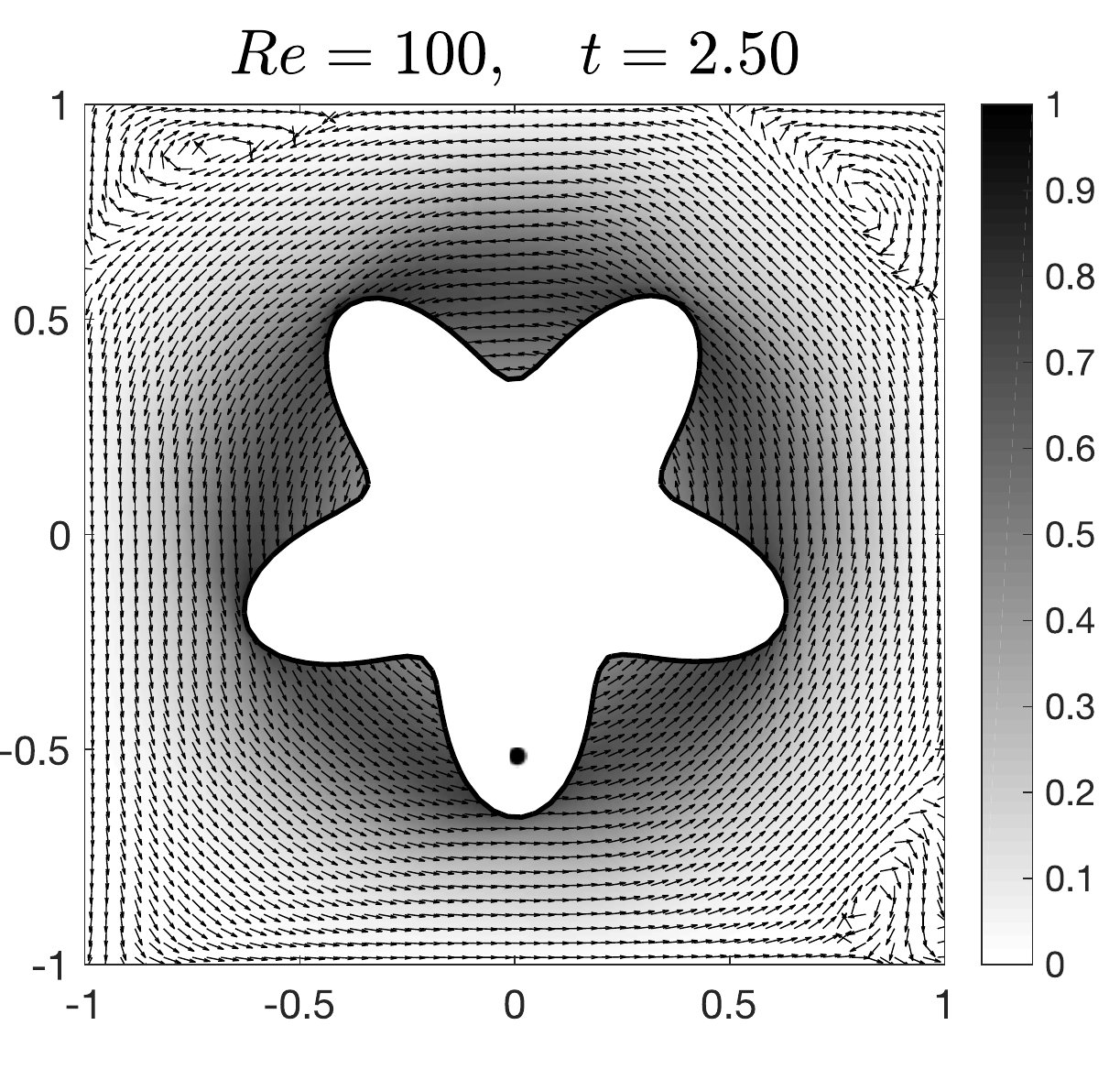}
 \end{minipage}
 \begin{minipage}[c]{0.33\textwidth}
   	\centering
		\includegraphics[width=0.99\textwidth]{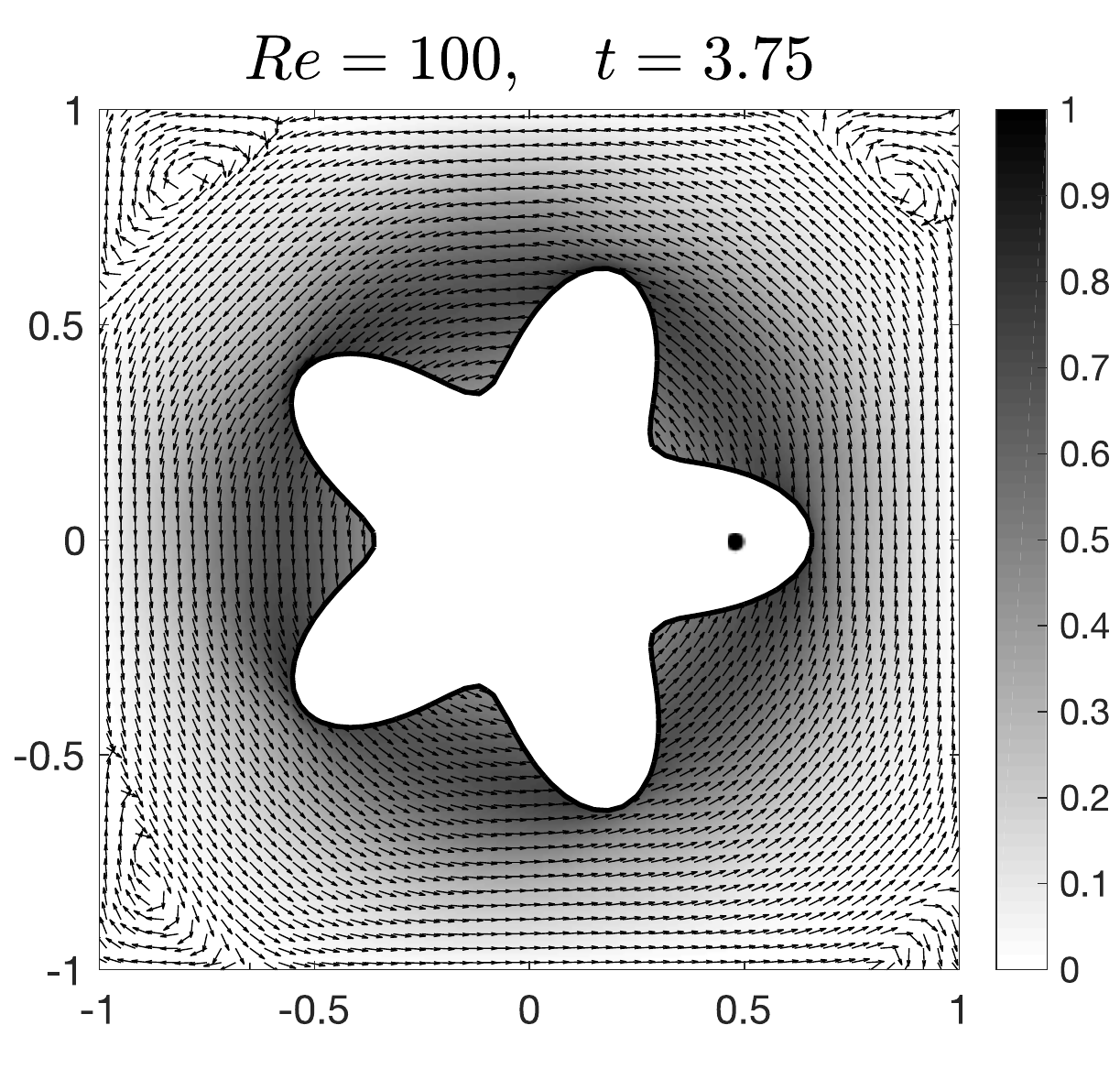}
 \end{minipage}
 \begin{minipage}[c]{0.33\textwidth}
   	\centering
		\includegraphics[width=0.99\textwidth]{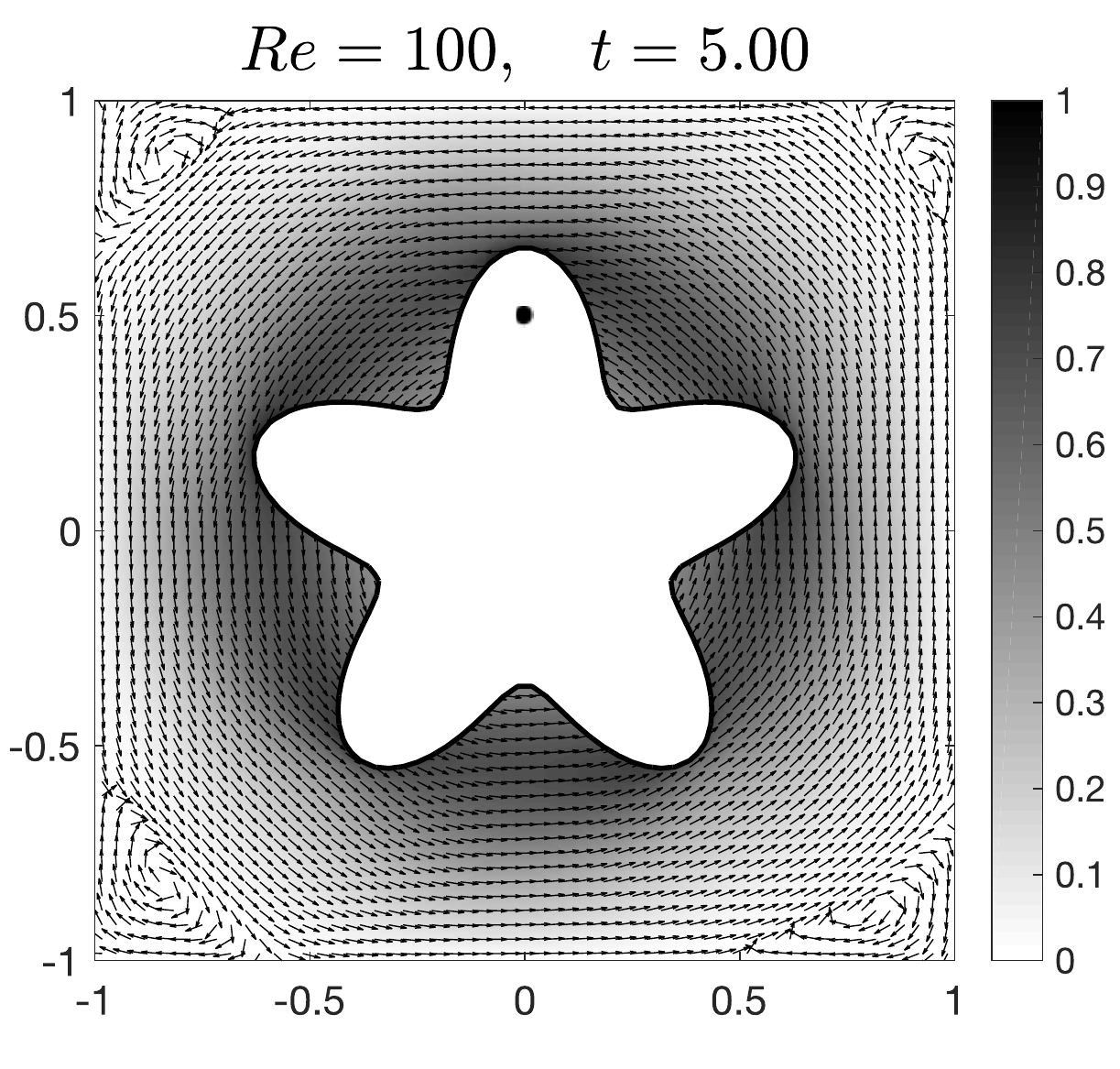}
 \end{minipage}
\begin{minipage}[c]{0.33\textwidth}
   	\centering
		\includegraphics[width=0.99\textwidth]{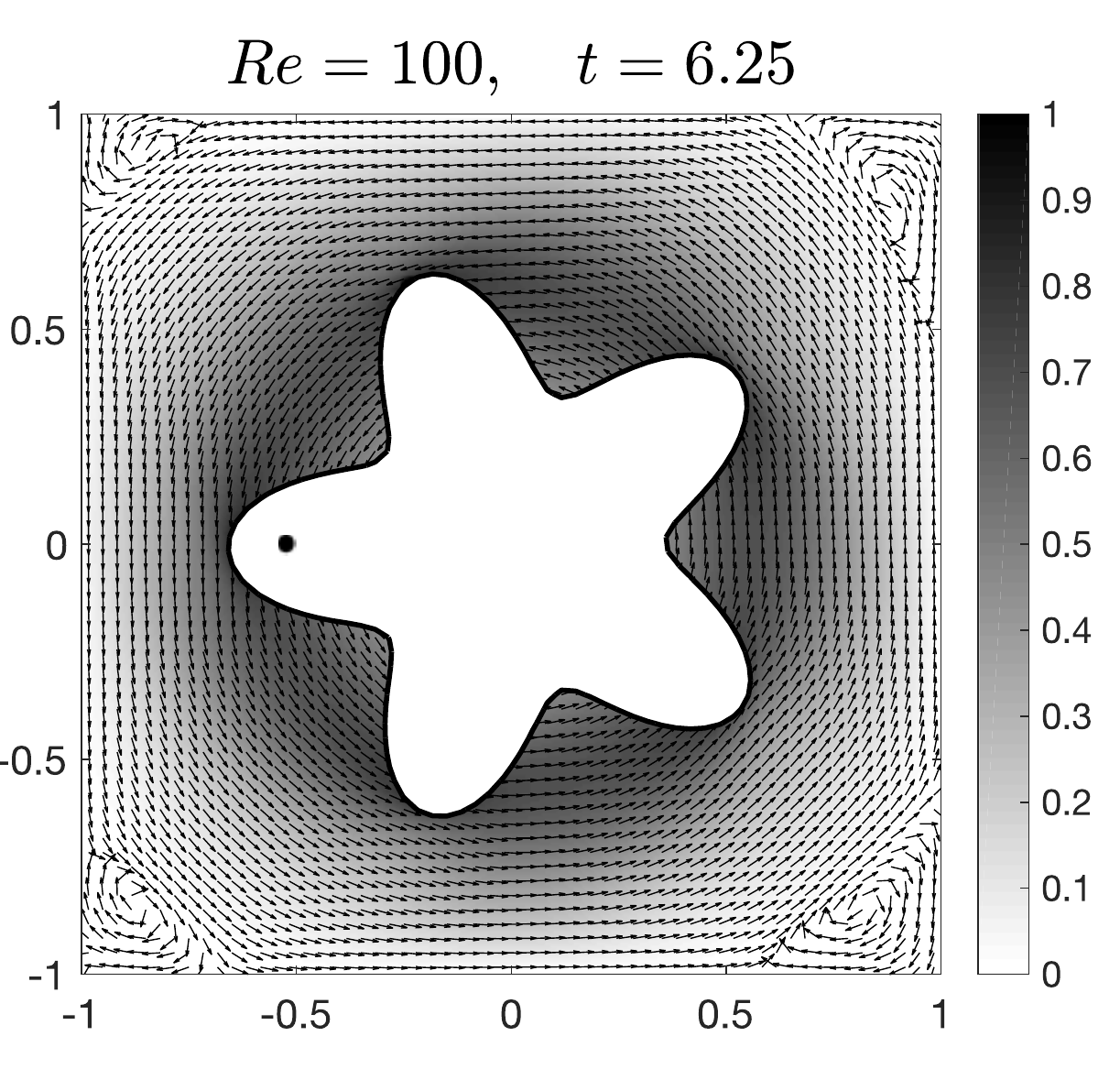}
 \end{minipage}
 \begin{minipage}[c]{0.33\textwidth}
   	\centering
		\includegraphics[width=0.99\textwidth]{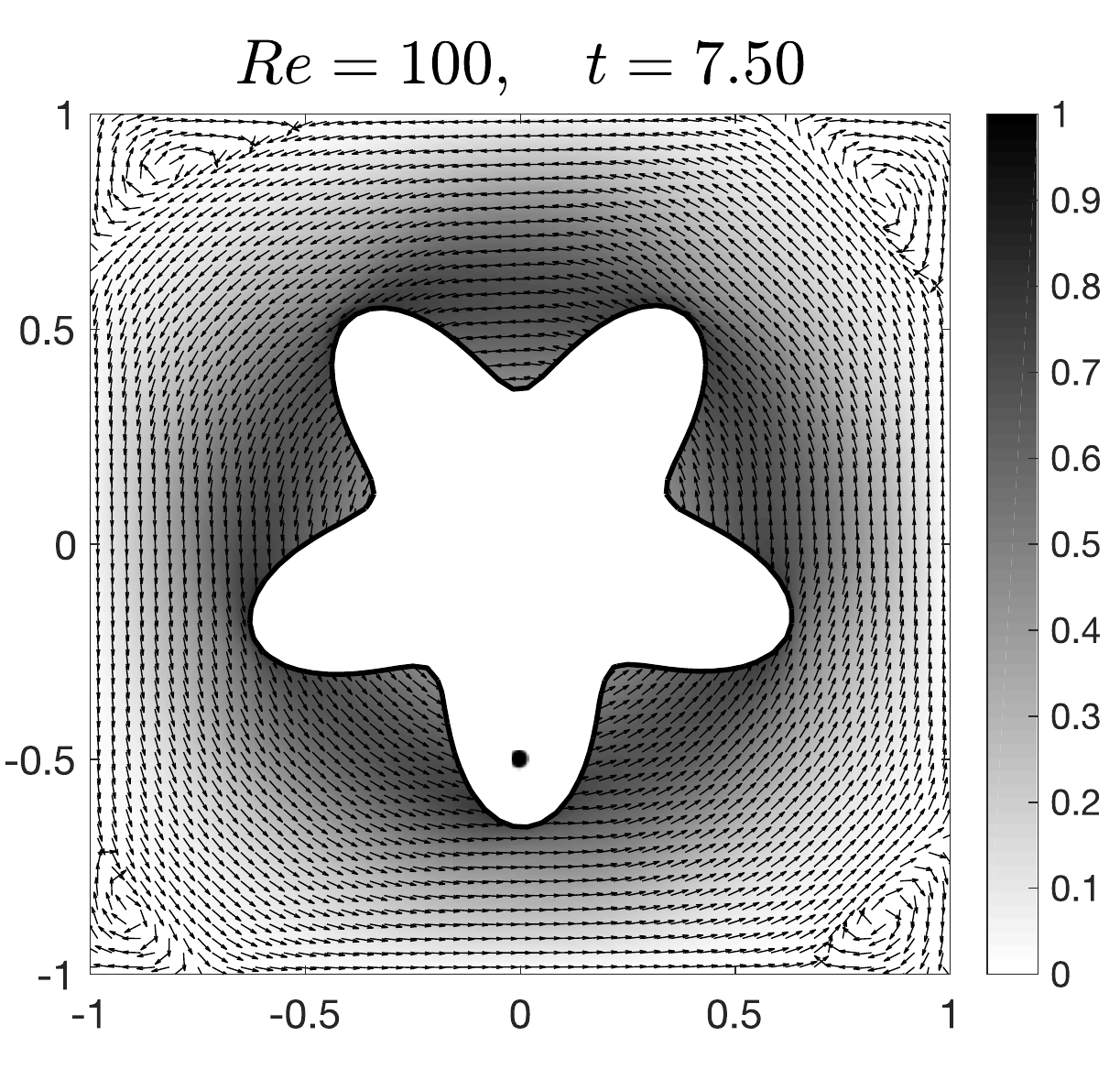}
 \end{minipage}
\begin{minipage}[c]{0.33\textwidth}
   	\centering
		\includegraphics[width=0.99\textwidth]{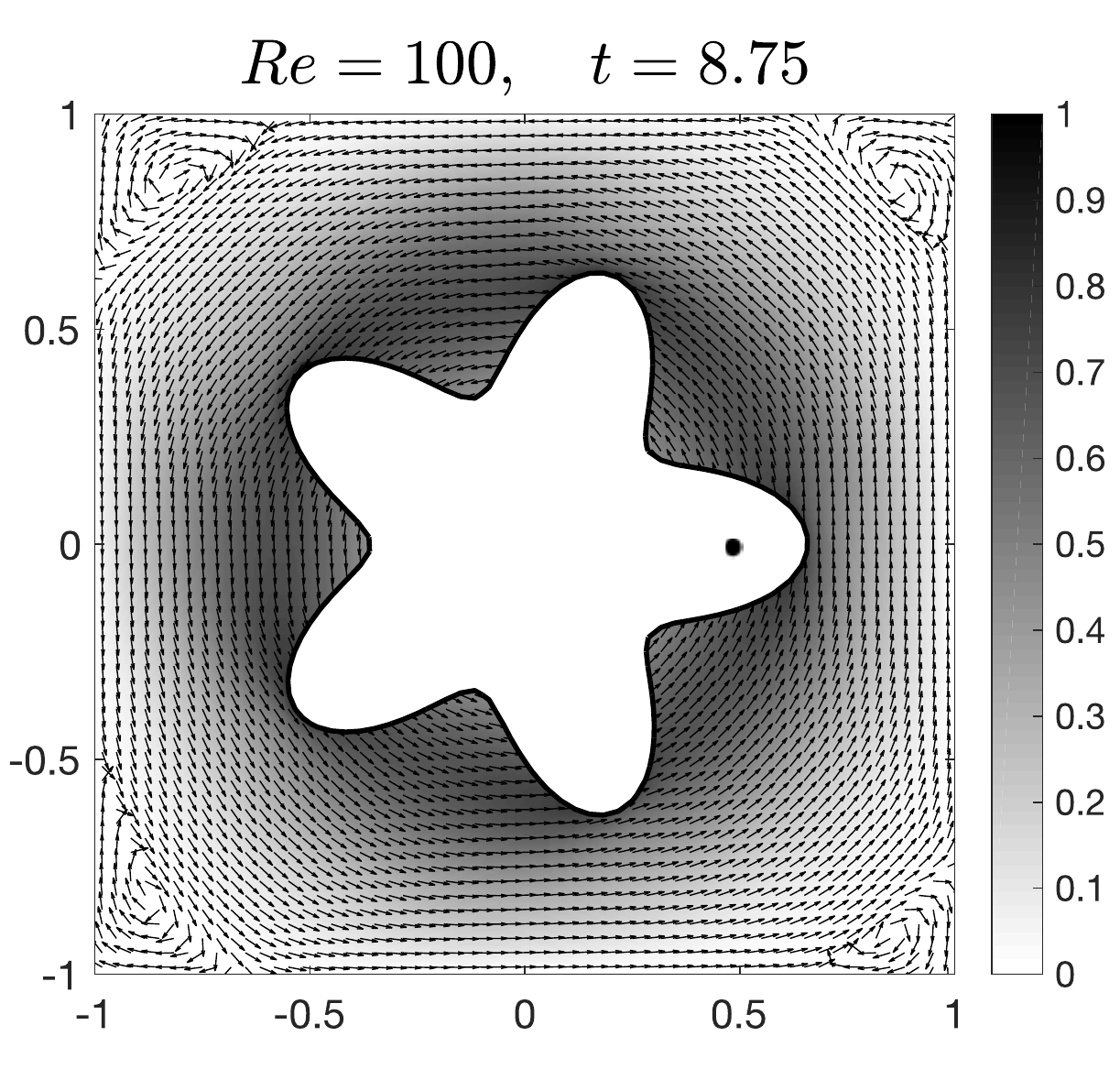}
 \end{minipage}
 \begin{minipage}[c]{0.33\textwidth}
   	\centering
		\includegraphics[width=0.99\textwidth]{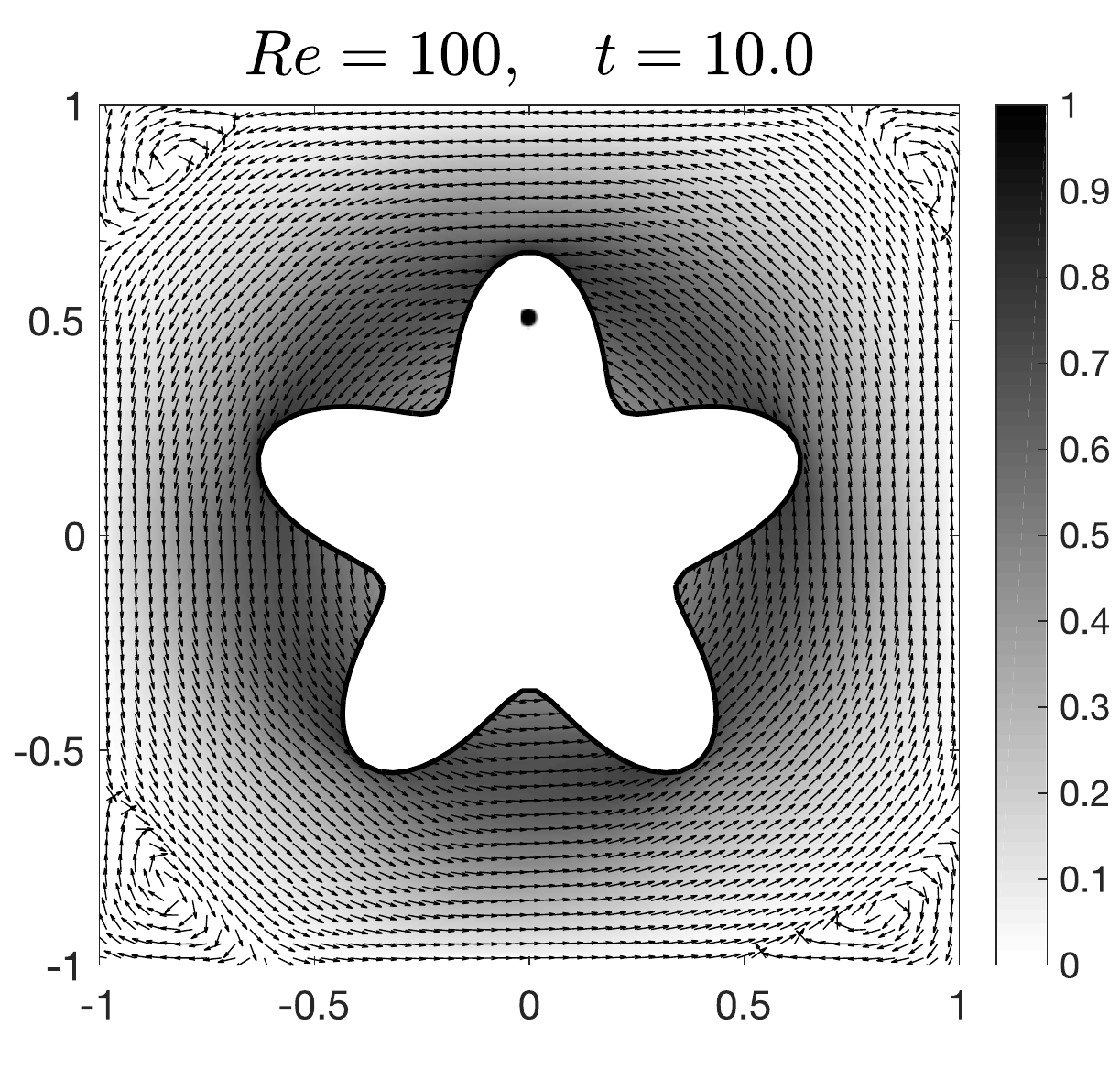}
 \end{minipage}
\caption{ \footnotesize{Velocity field at nine time steps and $N=60$. The black dot is plotted on the same petal and it is meant to help the visualisation of the rotation. At time $t=1.25$ two clockwise vortices occupy a large region of the top-right and bottom-right corners. At $t=2.50$ another clockwise vortex forms at the top-left corner. At $t=3.75$ the vortex on the bottom-right corner almost disappears, while at $t=5.00$ (i.e.~after a complete revolution) there is a clockwise vortex in each of the four corners. We observe that at time $t=5.00$ the vortices on the right occupy a smaller region than at $t=1.25$. From $t=5.00$ to $t=10.00$, i.e.~during the second revolution, the velocity field is almost constant in time away from the object, showing that it is approaching a steady-state solution.} }
	\label{fig:flower_moving_sol}
\end{figure}

The flower shaped object is rotated around the origin at a constant anticlockwise rotation velocity of $\omega = 2 \pi / 5$. The fluid is initially steady and no-slip boundary conditions on the walls ($\vec{u}=\vec{0}$) and on the boundary of the object ($(u(x,y),v(x,y),0)=(0,0,\omega) \times (x,y,0)$) are implemented for the entire simulation.
Reynolds number is $\text{Re}=100$.
Fig.\ref{fig:flower_moving_sol} shows the velocity field at nine time steps and $N=60$. The black dot is plotted on the same petal and it is meant to help the visualisation of the rotation. At time $t=1.25$ two clockwise vortices occupy a large region of the top-right and bottom-right corners. At $t=2.50$ another clockwise vortex forms at the top-left corner. At $t=3.75$ the vortex on the bottom-right corner almost disappears, while at $t=5.00$ (i.e.~after a complete revolution) there is a clockwise vortex in each of the four corners. We observe that at time $t=5.00$ the vortices on the right occupy a smaller region than at $t=1.25$. From $t=5.00$ to $t=10.00$, i.e.~during the second revolution, the velocity field is almost constant in time away from the object, showing that it is approaching a steady-state solution.

\subsubsection{Multigrid efficiency}
\begin{figure}[H]
   	\centering
		\includegraphics[width=0.69\textwidth]{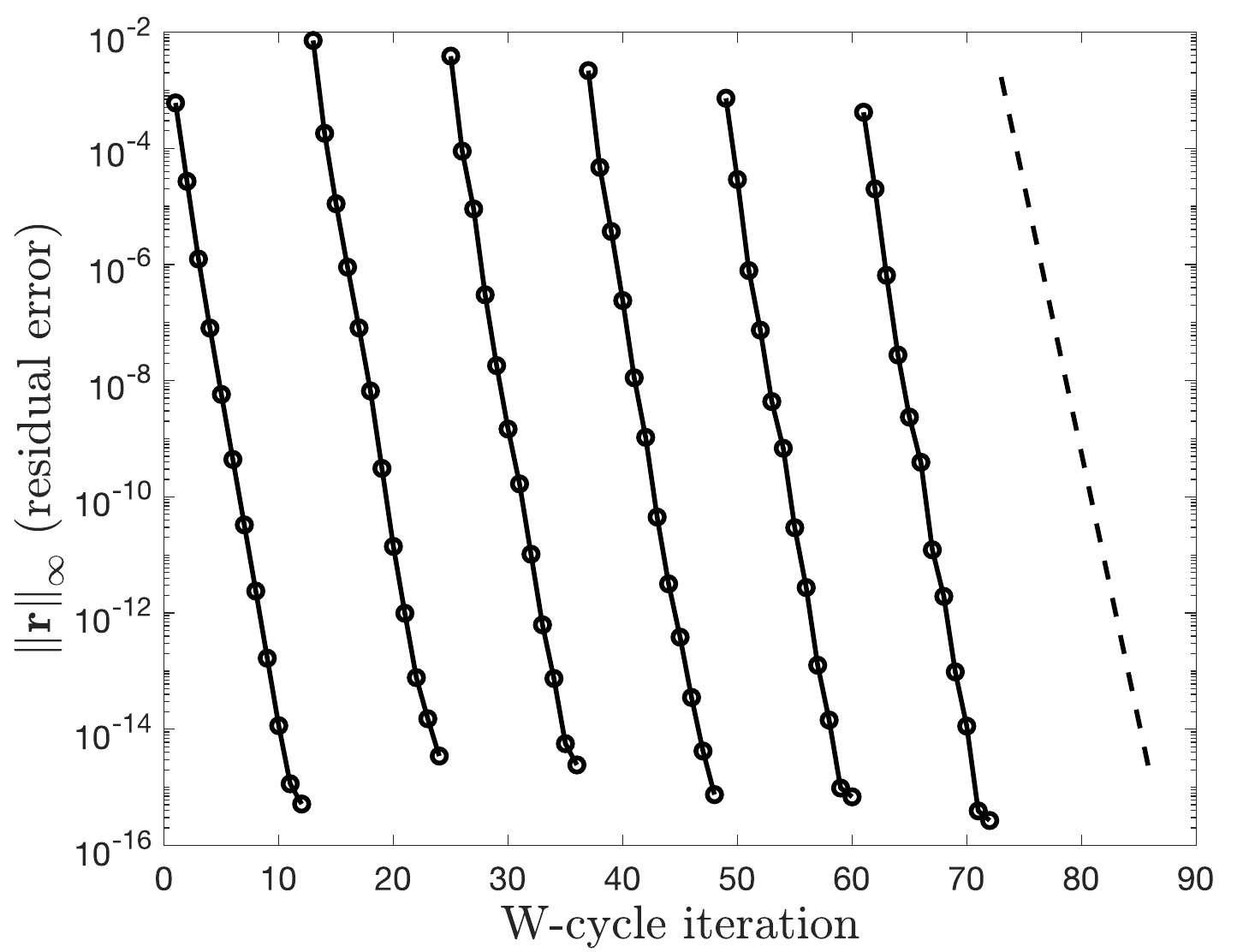}
\caption{ \footnotesize{residual error for each $W-$cycle. Each line refers to a single time step. The dashed line is the reference decay of $\rho=0.119$ (see Eq.~\ref{eq:convfact}), as expected by the Local Fourier Analysis for elliptic equations in rectangular domains.} }
	\label{fig:convfact}
\end{figure}

To investigate the efficiency of the multigrid, we compute the convergence factor 
\begin{equation}\label{eq:convfact}
\rho = \frac{\left\| \vec{r}_h^{(s+1)} \right\|_\infty}{\left\| \vec{r}_h^{(s)} \right\|_\infty},
\end{equation}
where $\vec{r}_h^{(s)}$ is the residual after $s$ W-cycles.
We use $\nu_1=2$ pre-relaxations and $\nu+1$ post-relaxations. For these values of $\nu_1$ and $\nu_2$, the Local Fourier Analysis predicts an asymptotic convergence factor of $\rho=0.119$ for standard multigrid applied to elliptic equations in rectangular domains (see~\cite[Sect.~4]{Trottemberg:MG}).
In Fig.~\ref{fig:convfact} we plot residuals obtained in the most complex problem of this paper, namely the rotating flower-shaped object of Sect.~\ref{sect:movingFlower}.
We observe that the convergence factor for the tests proposed in this paper is not degraded by the curved boundary and it is in agreement with the expected one $\rho=0.119$ (dashed line).

\section{Conclusion}\label{sect:concl}
We have presented a numerical method to solve the Navier-Stokes equations in a moving domain. The equations are discretized by Crank-Nicholson in time and finite differences in space. The moving domain is described by a level-set function and the curved boundaries are discretized by a ghost-point method. The resulting linear system is solved by a multigrid approach properly developed for curved boundaries.
The method is second order accurate in space and time. Several tests are provided, varying the Reynolds number and/or the shape of the domain.
The method falls under the class of unfitted boundary methods, since the curved boundaries are embedded in a fixed grid. These methods are preferred for moving domain problems, since fitted boundary methods require an additional computational cost for the re-meshing procedures (needed at each time step as the domain moves).

In fact, one of the bottlenecks of computer simulation is the computational time needed to represent complex-shaped objects in 3D. Most existing commercial software used by research centres and industries are based on fitted boundary method approaches, such as FEM. On the other hand, numerical approaches embedded in fixed grids have recently gained an increased interest in scientific computing, as simulation time is drastically reduced, because computationally expensive algorithms to generate the mesh are not required. This aspect is more exasperated if the domain moves and changes shape repetitively.
In addition, the increase of computational power of modern massively parallel supercomputers such as exascale computing demands for more efficient numerical methods to solve complex problems with higher accuracy.

The numerical approaches adopted in this paper, namely the ghost-point method on Cartesian meshes and the multigrid solver, are also designed to allow a more natural implementation in parallel computer architectures.
 
Finally, the method is suitable to model non turbulent viscous flows of incompressible fluids around moving objects, 
but can be applied in several other contexts where a moving/deforming domain is present, such as fluid-structure interaction or free-boundary problems. The latter is part of the author ongoing efforts.

\section*{Acknowledgments} The work has been partially supported by the London Mathematical Society Computer Science Small Grants -- Scheme 7 (Ref.~SC7-1617-02) and the Research in Pairs -- Scheme 4 (Ref.~41739).

\addcontentsline{toc}{chapter}{References}
\bibliographystyle{jabbrv_abbrv}
\bibliography{bibliography}

\end{document}